\documentclass[a4paper,12pt]{amsart}
\usepackage{amssymb,amsmath,latexsym,color,graphicx}%,stmaryrd}
\usepackage{hyperref}
\usepackage{amsfonts}
\usepackage{amsthm}
 \usepackage[color]{showkeys}
\usepackage{tikz}
\usetikzlibrary{positioning}

\newtheorem{dref}{Definition}[section] 
\newtheorem{theo}[dref]{Theorem} \newtheorem{prop}[dref]{Proposition}
\newtheorem{remark}[dref]{Remark}

\begin{document}

\title[]{Large $|k|$ behavior of complex geometric optics solutions 
to d-bar problems}
\author{Christian Klein}
\address{Institut de Math\'ematiques de Bourgogne, UMR 5584\\
                Universit\'e de Bourgogne-Franche-Comt\'e, 9 avenue Alain Savary, 21078 Dijon
                Cedex, France\\
    E-mail Christian.Klein@u-bourgogne.fr}
    
\author{Johannes Sj\"ostrand}
\address{Institut de Math\'ematiques de Bourgogne, UMR 5584\\
                Universit\'e de Bourgogne-Franche-Comt\'e, 9 avenue Alain Savary, 21078 Dijon
                Cedex, France\\
    E-mail Johannes.Sjostrand@u-bourgogne.fr}
    
\author{Nikola Stoilov}
\address{Institut de Math\'ematiques de Bourgogne, UMR 5584\\
                Universit\'e de Bourgogne-Franche-Comt\'e, 9 avenue Alain Savary, 21078 Dijon
                Cedex, France\\
    E-mail Nikola.Stoilov@u-bourgogne.fr}
\date{\today}
\begin{abstract}
	Complex geometric optics solutions to a system of d-bar equations 
	appearing in the context of electrical impedance tomography and 
	the scattering theory of the integrable Davey-Stewartson II 
	equations are studied for large values of the spectral parameter 
	$k$. For potentials \( q\in \langle \cdot \rangle^{-2} H^{s}(\mathbb{C}) \) for 
	some $s \in]1,2]$, it is shown that the solution converges as the 
	geometric series in $1/|k|^{s-1}$. For potentials $q$ being the 
	characteristic function of a strictly convex open set with smooth 
	boundary, this still holds with $s=3/2$ i.e., with  
	$1/\sqrt{|k|}$ instead of $1/|k|^{s-1}$.  The leading order contributions are computed 
	explicitly. Numerical simulations show the applicability of the 
	asymptotic formulae for the example of the characteristic 
	function of the disk. 
\end{abstract}

\thanks{This work is partially supported by 
the ANR-FWF project ANuI - ANR-17-CE40-0035, the isite BFC project 
NAANoD, the EIPHI Graduate School (contract ANR-17-EURE-0002) and by the 
European Union Horizon 2020 research and innovation program under the 
Marie Sklodowska-Curie RISE 2017 grant agreement no. 778010 IPaDEGAN}

\maketitle
\tableofcontents

\section{Introduction}
This paper is concerned with the large $|k|$ behavior of solutions to 
the Dirac system 
\begin{equation}\label{dbarphi}
  \begin{cases}
    \bar{\partial}\phi_{1}=\frac{1}{2}q\mathrm{e}^{\bar{k}\bar{z}-kz}\phi_{2},\\
    \partial\phi_{2}=\sigma\frac{1}{2}\bar{q}\mathrm{e}^{kz-\bar{k}\bar{z}}\phi_{1},\quad \sigma=\pm1,
\end{cases}   
\end{equation}
subject to the asymptotic conditions 
\begin{equation}
    \lim_{|z|\to\infty}\phi_{1}=1,\quad \lim_{|z|\to\infty}\phi_{2}=0;
    \label{Phisasym}
\end{equation}
here  $q=q(x,y)$ is a complex-valued field, the \emph{spectral 
parameter} $k\in\mathbb{C}$ is independent of $z=x+\mathrm{i} y$, and
\begin{equation*}
\partial:=\frac{1}{2}\left(\frac{\partial}{\partial x}-\mathrm{i}\frac{\partial}{\partial y}\right)\quad\text{and}\quad
\bar{\partial}:=\frac{1}{2}\left(\frac{\partial}{\partial 
x}+\mathrm{i}\frac{\partial}{\partial y}\right).
%\label{eq:d-dbar}
\end{equation*}
The functions $\phi_{i}(z;k)$, $i=1,2$  depend on $z$ and $k$ where it 
is understood that they need not be holomorphic in either variable.

As discussed in more detail in Subsection 
\ref{subdirac}, the solutions to the system (\ref{dbarphi}), subject to (\ref{Phisasym}) are  \emph{complex 
geometric optics} (CGO) solutions to a d-bar problem. The latter 
appears in the scattering theory of two-dimensional integrable 
equations as the Davey-Stewartson (DS) equation (\ref{eq:DSII}), in 
electrical impedance tomography (EIT) and in the theory of random matrix 
models. 
Of special interest in the context of the DS equation is the 
\emph{reflection coefficient} $R$, where
\begin{equation}
    \bar{R} = 
    \frac{2\sigma}{\pi}\int_{\mathbb{C}}^{}\mathrm{e}^{kz-\bar{k}\bar{z}}\bar{q}(z)\phi_{1}(z;k)d^{2}z
    \label{Rintro},
\end{equation}
which can be seen as a nonlinear analog to the Fourier transform of 
the potential $q$. 

The existence and uniqueness of CGO solutions to system 
(\ref{dbarphi}) with $\sigma=1$ was studied in \cite{BC} for Schwartz class 
potentials and in \cite{Su1,Su2,Su3} for potentials  
$q\in L^{\infty}(\mathbb{C})\cap L^{1}(\mathbb{C})$ such that 
also $\hat{q}\in L^{\infty}(\mathbb{C})\cap L^{1}(\mathbb{C})$ 
where $\hat{q}$ is the Fourier transform of $q$ (the potentials have 
to satisfy a smallness condition in the focusing case).
The results for system (\ref{dbarphi}) for Schwartz class potentials were generalized respectively to real-valued, compactly supported potentials 
in $L^{p}(\mathbb{C})$ \cite{BU} and to potentials in 
$H^{1,1}(\mathbb{C})$ \cite{Perry2012}, and in \cite{NRT} to potentials 
in  $L^{2}(\mathbb{C})$.

For  potentials $q$ in the Schwartz class of rapidly decreasing smooth potentials, the reflection coefficient (\ref{Rintro}) is also 
   in the Schwartz class. However, the dependence of $\phi_{i}$, $i=1,2$ 
   and $R$ on $k$ is much less clear for  potentials $q$ of lower 
   regularity, such as potentials with compact support, or 
%    potentialss of slow 
   slow decrease towards infinity, which 
   are both important in EIT. It is the purpose of the present paper 
   to give the large $|k|$ behaviour of the solutions to the system 
   \ref{dbarphi} for such cases.

The importance of d-bar problems in applications has led to many 
numerical approaches to solve them. The standard method is to use 
that the inverse of the d-bar operator is given by the solid Cauchy 
transform, a weakly singular integral which was first computed in the current context in  \cite{KnMuSi2004}, see \cite{MS} for a review of more recent 
developments, with Fourier methods and a simple regularization of the 
integrand. These approaches are of first order which means the 
numerical error decreases linearly with the number of Fourier modes.  
The first Fourier approach with  an exponential decrease of the numerical error with the 
number of Fourier modes, or \emph{spectral convergence}, was presented in \cite{KM,KMS} for Schwartz  class potentials via an analytical regularization of the integrand in  the solid Cauchy transform. For potentials with compact support on a 
disk, a numerical approach of formally infinite order was presented 
in \cite{KS19} based on a formulation of the problem in 
polar coordinates and the solution of the resulting system by a 
Chebyshev-Fourier method. As will be shown in this paper, the 
reflection coefficient (\ref{Rintro}) decreases algebraically in 
$1/|k|$ in this case as $|k|\to\infty$. Thus in contrast to the case 
of Schwartz potentials, a purely numerical approach cannot be of the 
order of machine precision (here $10^{-16}$) for all values of
$k$. One of the
motivations for the present paper is to present, in the concrete example 
of the characteristic function of the disk as the potential $q$, 
formulae for the large $|k|$ behaviour which together with the 
numerical approach \cite{KS19} yield %for a given accuracy 
a complete description of the solutions in the whole complex $k$ plane
within a predetermined accuracy. This gives an upper
  bound on the reflection coefficient, and  we hope to be able to push
  the methods further for the leading asymptotics.

An interesting question in the context of the DS equation is the 
appearance of \emph{dispersive shock waves} (DSWs), zones of rapid 
modulated oscillations in the vicinity of shocks in the 
semi-classical DS system for the same initial data. Such DSWs were 
studied numerically in \cite{KR14}. A first attempt towards an 
asymptotic description of DSWs based on inverse scattering techniques 
was presented in \cite{AKMM}. Note that the first system for which a 
rather complete asymptotic description of DSWs exists is the 
completely integrable Korteweg-de Vries (KdV) equation. Historically the 
Gurevitch-Pitaevskii (GP) work \cite{GP} on solutions to the KdV equation for 
steplike initial data was very influential. It was one of the 
motivation of the numerical work \cite{KS19} to provide numerical 
tools for the study of the corresponding GP problem for DS II, 
initial data given by the characteristic function of the disk. 
In the present work, this case is addressed in some 
detail. 
   
\subsection{Applications of the Dirac system}
\label{subdirac}

The system (\ref{dbarphi}) and the conditions (\ref{Phisasym}) are
equivalent to the Dirac system (simply put $\phi_{1}:=\psi_{1}e^{-kz}$, 
$\phi_{2}:=\psi_{2}e^{-\overline{kz}}$)
\begin{equation}
\begin{array}{l}
    \bar{\partial}\psi_{1}=\frac{1}{2}q\psi_{2},\\
    ~\\
    \partial\psi_{2}=\sigma\frac{1}{2}\bar{q}\psi_{1},\quad 
	\sigma=\pm1
    \label{dbarpsi}
\end{array}
\end{equation}
where the scalar functions $\psi_{1}$ and $\psi_{2}$ satisfy the 
CGO asymptotic conditions
\begin{equation}
\begin{array}{l}
    \lim_{|z|\to\infty}\psi_{1}\mathrm{e}^{-kz}=1,\\
    ~\\
    \lim_{|z|\to\infty}\psi_{2}\mathrm{e}^{-\bar{k}\bar{z}}=0.
    \label{dbarpsiasym}
\end{array}
\end{equation}
Putting $\Psi_{\pm}:=\psi_{1}\pm\bar{\psi}_{2}$, the system (\ref{dbarpsi}) 
is diagonalized,
\begin{equation}
    \bar{\partial}\Psi_{\pm}=\pm\frac{1}{2}q\bar{\Psi}_{\pm}
    \label{scalar}
\end{equation}
subject to the asymptotic condition 
$\lim_{|z|\to\infty}\Psi_{\pm}\mathrm{e}^{-kz}=1$. The disadvantage of 
equation (\ref{scalar}) is that it is not complex linear in $\Psi_{\pm}$, 
which is why we use the system (\ref{dbarpsi}). 

The system (\ref{dbarpsi}) has many applications, the first being in 
completely integrable equations in two dimensions. As was shown in 
\cite{Fok,FA} system 
(\ref{dbarpsi}) gives both the scattering and 
inverse scattering map for the  Davey-Stewartson II 
equation
\begin{equation}
\begin{split}
i q_t + (q_{xx}-q_{yy}) + 2\sigma(\Phi + |q|^2)q&=0,\\
\Phi_{xx}+\Phi_{yy} +2(|q|^2)_{xx}&=0,
\end{split}
    \label{eq:DSII}
\end{equation}
a two-dimensional nonlinear Schr\"odinger equation; here the 
parameter $\sigma=1$ in 
the \emph{defocusing} case, and $\sigma=-1$
in the \emph{focusing} case. DS systems, which are in general not 
integrable, appear in the modulational regime of many dispersive 
equations as for instance the water wave systems, see e.g.,  \cite{KS15} 
for a review on DS equations and a 
comprehensive list of references. 

In the smooth case the scattering data  are given in terms of the reflection coefficient $R=R(k)$  in terms of 
$\psi_2(z;k)$ as follows (in the general case, one has to consider 
(\ref{Rintro})):
\begin{equation}
\mathrm{e}^{-kz}\overline{\psi_2(z;k)}=\frac{1}{2}R(k) z^{-1}+ 
O(|z|^{-2}),\quad |z|\to\infty,
\label{eq:r-def}
\end{equation}
which is equivalent to (\ref{Rintro}) after writing $\phi_{2}$  as 
the solid Cauchy transform of the right hand side of the second equation 
of (\ref{dbarphi}) and taking the limit $z\to\infty$. 
Note that the understanding of the  Dirac system (\ref{dbarphi}) 
with $\sigma=-1$
is much less complete than in the 
defocusing case $\sigma=1$.   In the former case the system no longer 
has generically a unique solution for large classes of potentials $q$ 
for all $k\in\mathbb{C}$. There can be special values of the spectral parameter, 
called exceptional points, where the system is not uniquely solvable. 
Therefore we concentrate for the examples (which cover all values 
of $k\in\mathbb{C}$) on the defocusing
case. But the results for large $|k|$ we present in this paper hold 
for both cases. Note that this implies that the exceptional points 
can only occur in a bounded set. Perry \cite{perrysol} gave a bound 
for the radius of this set based on the $H^{1,1}$ norm of the 
potential. It is beyond the scope of the current paper to establish 
similar bounds based on our approach.

As $q$ in (\ref{eq:DSII}) evolves in time $t$, the reflection coefficient evolves by a trivial phase factor:
\begin{equation}
R(k;t)=R(k,0)e^{4i t\Re(k^2)}.
\end{equation}
The inverse scattering transform for DS II is then given by 
(\ref{dbarpsi}) and
(\ref{dbarpsiasym}) after replacing $q$ by $R$ and vice versa, the derivatives with 
respect to $z$ by the corresponding derivatives with respect to $k$, 
and asymptotic conditions for $k\to\infty$ instead of $z\to\infty$, 
see \cite{AF}.

Systems of the form (\ref{dbarpsi}) also appear in electrical impedance tomography 
(EIT) in 2d, the reconstruction of 
the conductivity in a given domain $\Omega$ from measurements of the electrical 
current through its boundary, induced by an applied voltage, i.e., from 
the Dirichlet-to-Neumann map. This problem was first posed by Calder\'on
\cite{calderon} and bears his name. For a comprehensive review of the 
mathematical aspects and advances see \cite{Uhl,MS}. The basic idea 
of EIT is to construct CGO solutions to the conductivity equation in 
some domain $\Omega\in\mathbb{R}^{2}$,
\begin{equation}
    \nabla \cdot(\sigma(x,y)\nabla u(x,y)) = 0, \quad (x,y)\in \Omega 
    \label{conduc}.
\end{equation}
In \cite{BU}, this was done in form of the system (\ref{dbarpsi}) by 
putting $q=-(1/2)\partial \ln \sigma$ for conductivities $\sigma\in 
C^{1}(\Omega)$. The CGO solutions satisfy slightly different 
asymptotic  conditions in this case which is why we study  larger 
classes of conditions than (\ref{dbarpsiasym}). 
The reconstruction of the conductivity from the Dirichlet-to-Neumann 
map is also achieved via a d-bar problem, see \cite{Uhl}. 

D-bar problems also appear
 in the context of 2d orthogonal polynomials, and of 
Normal Matrix Models in Random Matrix Theory, see e.g.\ \cite{KM}.

\subsection{Main results}

\medskip
\paragraph{\underline{\it 1.\ The d-bar equation on $\mathbb{C}$}}\ \\
For $s\in\mathbb{R}$, we write $\langle \cdot 
\rangle^{s}L^{2}=\left\{\langle z \rangle^{s}u(z); u\in 
L^{2}(\mathbb{C})\right\}$, where $\langle z 
\rangle^{s}=(1+|z|^{2})^{1/2}$. Thus 
$||\tilde{u}||_{\langle z \rangle^{s}L^{2}}=||u||_{L^{2}}$ if 
$\tilde{u}=\langle z 
\rangle^{s}u\in \langle \cdot \rangle^{s}L^{2}(\mathbb{C})$.\\
{\it Reviewing
  H\"ormander's approach with Carleman estimates, we show in the
  propositions \ref{db1}, \ref{db2} that if $0<\epsilon \le 1$, then
  for every $v\in \langle \cdot \rangle^{\epsilon -2}L^2(\mathbb{C})$ 
  the equation $\overline{\partial }u=v$ in (\ref{db.10}) has a unique
  solution $u\in \langle \cdot \rangle^\epsilon L^2$. When $\epsilon
  =1$ we show in Proposition \ref{db3} that  when $v\in \langle \cdot 
  \rangle^{-1}L^2$ the unique solution $u\in\langle \cdot \rangle L^2$
  is given by the
  standard  formula (\ref{db.16}),
  $$u(z)=\frac{1}{\pi }\int \frac{1}{z-w}v(w)L(dw),$$
  where $L(dw)$ denotes the Lebesgue measure on $\mathbb{C}\simeq
  \mathbb{R}^2$. Of course we have the same results for the complex
  conjugate equation $\partial u=v$ and we then have to replace
  $1/(z-w)$ with $1/(\overline{z}-\overline{w})$ in the integral above.
}

\medskip
\par
We are interested in the case when $|k|$ is large and notice
  that
  \begin{equation}\label{omega}
kz-\overline{kz}=i|k|\Re (z\overline{\omega })=i|k|\langle z,\omega
\rangle_{\mathbb{R}^2},\ \ \omega= \frac{2i\overline{k}}{|k|},
  \end{equation}
so that $ |\omega |=2$. Here we identify $\mathbb{C}\simeq \mathbb{R}^2$ in
the usual way.
Introduce the
semi-classical parameter $h=1/|k|$, $0<h\ll 1$
and put 
$$
\widehat{\tau }_\omega u(z)=e^{\frac{i}{h}\langle z,
  \omega\rangle_{\mathbb{R}^2} }u(z).
$$
$\widehat{\tau }_\omega $ is translation by $\omega \in\mathbb{R}^2$ on the
$h$-Fourier transform side, see (\ref{sy.7.5}). 

\par For $v\in \langle \cdot \rangle^{-1}L^2$, let $u$,
$\widetilde{u}$ be the unique solutions in $\langle \cdot
  \rangle L^2$ (see Section \ref{db})
of the equations
\begin{equation}\label{sy.9int}
h\overline{\partial }u=v,\ h\partial \widetilde{u}=v,
\end{equation}
and write $u=Ev$, $\widetilde{u}=Fv$.  

The goal is to solve the system (\ref{dbarphi}) iteratively for small 
\( h \). With the above notation this implies that we can write  (\ref{dbarphi}) in the 
form ($\phi_{i}\in \langle \cdot \rangle^{-1}L^2$,
$i=1,2$)
\begin{equation}\label{sy.14int}
\begin{cases}
\phi _1-E\widehat{\tau }_{-\omega }\frac{hq}{2}\phi _2=E\Psi _1,\\
\phi _2-\sigma F\widehat{\tau }_\omega \frac{h\overline{q}}{2}\phi _1=F\Psi _2,
\end{cases}
\end{equation}
with $\Psi_{1}=\Psi_{2}=0$, 
or 
\begin{equation}\label{sy.15int}
(1-\mathcal{K})\begin{pmatrix}\phi _1\\\phi
  _2\end{pmatrix}=\begin{pmatrix}E\Psi _1\\F\Psi _2\end{pmatrix},
\end{equation}
where
\begin{equation}\label{sy.16int}
\mathcal{K}=\begin{pmatrix}0 &A\\ B &0\end{pmatrix},\ \ 
\begin{cases}A=E\widehat{\tau }_{-\omega }\frac{hq}{2},\\
B=\sigma F\widehat{\tau }_\omega \frac{h\overline{q}}{2}\ (=\sigma\overline{A})
\end{cases} .
\end{equation}
The inhomogeneities \( \Psi_{1} \), \( \Psi_{2} \) are introduced in 
the above system to allow for an iterative solution in $h$. In 
leading order they depend on the asymptotic conditions, e.g.\
(\ref{Phisasym}), which will not be specified for the moment in order 
to allow for rather general conditions.

\medskip
\paragraph{\underline{\it 2.\ Inverting the system}}\ \\ 
	{\it Let \( q\in \langle \cdot \rangle^{-2} H^{s}(\mathbb{C}) \) for 
	some $s \in]1,2]$ and fix $\epsilon \in ]0,1]$. Then 
	\( \mathcal{K}^{2} =\mathcal{O}(h^{s-1}): 
	(\langle\cdot\rangle^{\epsilon} L^{2})^{2} \mapsto 
	(\langle\cdot\rangle^{\epsilon} L^{2})^{2} \), see 
	Proposition \ref{sy1}.
	It follows that \( 1-\mathcal{K} \) is bijective with inverse \( 
	(1+\mathcal{K})(1-\mathcal{K}^{2})^{-1} 
	=(1-\mathcal{K}^{2})^{-1}(1+\mathcal{K})=1+\mathcal{K}+\mathcal{O}(h) \).\\
	This implies that the  solution of (\ref{sy.15int}) in terms of a 
	Neumann series in \( 
	h \) converges like the geometric series, see Proposition
        \ref{sy2}.\\
        Let $\phi _1=\phi _1^0+\phi _1^1$, $\phi _2=\phi _2^0+\phi
        _2^1$ be the solution, where $\phi _1^0=1$, $\phi _2^0=0$ is
        an approximate solution (cf.\ (\ref{sy.11}), (\ref{sy.12}))
        and $(\phi _1^1,\phi _2^1)$ the correction. In (\ref{sy.41}),
        (\ref{sy.42}), (\ref{sy.44}) we show that $\phi
        _2^1-F\widehat{\tau }_\omega h\overline{q}/2$ and $\phi
        _1^1-AF\widehat{\tau }_\omega h\overline{q}/2$ are small in
        $\langle \cdot \rangle ^\epsilon L^2$ for every fixed
        $\epsilon \in ]0,1]$. 
      }

      \medskip
      \par Apriori, we only have $\mathcal{ K}=\mathcal{ O}(1): (\langle \cdot
      \rangle^\epsilon L^2)^2\to \langle \cdot
      \rangle^\epsilon L^2)^2$ which is not enough for the convergence
      of the Neumann series. The improvement for $\mathcal{ K}^2$ comes
      from phase space truncations (microlocal analysis).

It appears difficult to extend this result to general lower regularity of \( 
q \). Therefore we concentrate on the case of potentials with 
compact support, in particular on potentials being the characteristic 
function of a simply connected \( \Omega\Subset\mathbb{C} \) with 
smooth boundary.

\medskip\paragraph{\underline{\it 3. Inversion when $q$ is a characteristic
  function}}\ \\
{\it Let \( q= \mathbf{1}_{\Omega} \) where
  $\Omega \Subset \mathbb{C}^2$ is a strictly convex open set with
  smooth boundary. Then the results in Item 2 are valid with
  $s=3/2$. In particular,
  \( \mathcal{K}^{2} =\mathcal{O}(h^{1/2}):
  (\langle\cdot\rangle^{\epsilon_{0}} L^{2})^{2} \mapsto
  (\langle\cdot\rangle^{\epsilon_{0}} L^{2})^{2} \). See Proposition
  \ref{qcf1}.}

\medskip\paragraph{\underline{\it 4. Asymptotics when $q$ is a
    characteristic
    function}}\ \\
{\it Theorem \ref{ic2} gives a detailed asymptotic description in
  various regions of the function
  $f(z,k)=\overline{F\widehat{\tau }_\omega h\overline{q}/2}$ (and
  hence also of $F\widehat{\tau }_\omega h\overline{q}/2$). Combining this with the
  approximation results in the items 2, 3, we get as a consequence the
  estimates (\ref{wn.12})
	$$\phi _2^1=\frac{1}{2k}e^{i|k|\Re (\cdot \overline{\omega })}1_\Omega +\mathcal{
  O}(1)h^{3/2}(\ln (1/h))^{1/2}\hbox{ in }\langle \cdot
\rangle^\epsilon L^2,
$$ and (\ref{wn.14})
	$$\phi _1^1=\frac{h}{4k}E(1_\Omega )+\mathcal{ O}(1)h^{3/2}(\ln (1/h))^{1/2}
\hbox{ in }\langle \cdot \rangle^\epsilon L^2. 
$$} This allows us to estimate the 
reflection coefficient via (\ref{Rintro}). To get an asymptotic
formula with a non-vanishing leading term seems to require more work
however, because of possible cancellations.

\subsection{Outline of the paper}

The paper is organized as follows: in Section \ref{db} we introduce 
some notation and summarize H\"ormander's approach to d-bar problems 
and weighted Carleman estimates. In Section \ref{sy} we present a 
proof of the first part of the main results. The case of the
characteristic function of a simply connected domain in \( \mathbb{C} 
\) with smooth boundary is addressed in Section \ref{qcf}. In Section 
\ref{ic} we give explicit formulae for an integral appearing in the 
special case of the characteristic function of a simply connected 
compact domain with smooth boundary. In Section 
\ref{num} we present a numerical study of the Dirac system for the 
example of the characteristic function of the disk and address the 
question when the asymptotic formulae for small \( h \)  can be 
applied in practice.

\section{The $\overline{\partial }$-operator on $\mathbb{C}$ with
  polynomial weights}\label{db}
\setcounter{equation}{0}
In this section we introduce basic notation and 
review H\"ormander's solution of the d-bar equations with 
Carleman estimates, see \cite{Hor}.

Consider $$
\Phi (z)=\ln \langle z\rangle^2,\ \langle z\rangle =
(1+|z|^2)^{1/2}=(1+z\overline{z}) ^{1/2},\ z\in \mathbb{C}.$$ 
Then 
\begin{equation}\label{db.1}
\partial _{\overline{z}}\partial _z\Phi
=\frac{1}{(1+z\overline{z})^2}=\frac{1}{\langle z\rangle^4},
\end{equation}
where we use the standard notation $\ z=x+iy,\ x=\Re z,\ y=\Im z$ and 
\begin{equation*}
\begin{array}{c}
\partial =\partial _z=\frac{1}{2}(\partial _x+\frac{1}{i}\partial
_y),\ \overline{\partial }=\partial _{\overline{z}}=\frac{1}{2}(\partial _x-\frac{1}{i}\partial
_y).
\end{array}
\end{equation*}
\par For $\epsilon \in \mathbb{R}$, put
\begin{equation}\label{db.2}
P_\epsilon =\langle \cdot \rangle^{-\epsilon }\circ \overline{\partial
}\circ \langle \cdot \rangle^\epsilon =e^{-\epsilon \Phi /2}\circ
\overline{\partial }\circ e^{\epsilon \Phi /2}=\overline{\partial }+\epsilon
\overline{\partial }\Phi /2,
\end{equation}

\begin{equation}\label{db.3}
P_\epsilon^* =-\langle \cdot \rangle^{\epsilon }\circ \partial
\circ \langle \cdot \rangle^{-\epsilon } =-e^{\epsilon \Phi /2}\circ
\partial \circ e^{-\epsilon \Phi /2}=-\partial +\epsilon
\partial \Phi /2.
\end{equation}
Here $^*$ denotes the complex adjoint in $L^2(\mathbb{C})$.

\par We have for the commutator:
\begin{equation}\label{db.4}
[P_\epsilon ,P_\epsilon ^*]=\epsilon \partial \overline{\partial }\Phi
=\frac{\epsilon }{\langle z\rangle^4}.
\end{equation}
When $\epsilon >0$, we get, using a standard trick:
\begin{equation}\label{db.5}
\| P_\epsilon ^*\phi \|^2\ge \| P_\epsilon ^*\phi \|^2-\| P_\epsilon
\phi \|^2 = ([P_\epsilon ,P_\epsilon ^*]\phi |\phi )=\epsilon \|
\langle \cdot \rangle^{-2}\phi \|^2,
\end{equation}
for every $\phi \in C_0^\infty (\mathbb{C})$, where $\|\cdot \|$ and
$(\cdot |\cdot )$ denote the norm and scalar product on
$L^2$. Hence, we have the apriori estimate,
\begin{equation}\label{db.6}
\epsilon ^{1/2}\| \langle \cdot \rangle^{-2}\phi \|\le \| P_\epsilon
^*\phi \| .
\end{equation}

\par This leads to an existence result for $P_\epsilon $ in the usual
way. Assume that $v\in \langle \cdot \rangle^{-2}L^2(\mathbb{C})$, i.e.\
$\langle \cdot \rangle^2v\in L^2$. Then for every $\phi \in C_0^\infty
(\mathbb{C})$:
$$
|(\phi | v )|=|( \langle \cdot
\rangle^{-2}\phi | \langle \cdot \rangle^2 v )|\le \epsilon ^{-1/2}\| \langle \cdot \rangle^2 v\|
\| P_\epsilon ^* \phi \|.
$$
Hence $\phi \mapsto (\phi |v )$ is a bounded linear form acting on
$P_\epsilon ^*\phi $, so $\exists\, u\in L^2$ with
\begin{equation}\label{db.7}
\| u\|\le \epsilon ^{-1/2}\| \langle \cdot \rangle^2v\| ,
\end{equation}
such that
$$
(\phi | v )=(P_\epsilon ^*\phi |u ),\ \forall\, \phi \in C_0^\infty 
(\mathbb{C}),
$$
i.e.
\begin{equation}\label{db.8}
P_\epsilon u=v .
\end{equation}
This can be written
\begin{equation}\label{db.9}
\overline{\partial }\underbrace{\langle \cdot \rangle^\epsilon
  u}_{\widetilde{u}}=\underbrace{\langle \cdot \rangle^\epsilon v}_{\widetilde{v}},
\end{equation}
where $\langle \cdot \rangle^{2-\epsilon }\widetilde{v}=\langle \cdot
\rangle^2v\in L^2$, $\langle \cdot \rangle^{-\epsilon
}\widetilde{u}=u\in L^2$ and after dropping the tildes, we get:
\begin{prop}\label{db1}
Let $\epsilon >0$. For every $v\in \langle \cdot \rangle^{\epsilon
  -2}L^2$, there exists $u\in \langle \cdot \rangle^{\epsilon }L^2$
such that
\begin{equation}\label{db.10}
\overline{\partial }u=v,
\end{equation}
and 
\begin{equation}\label{db.11}
  \| \langle \cdot \rangle^{-\epsilon }u\|\le \epsilon ^{-1/2}\| \langle
  \cdot   \rangle^{2-\epsilon }v\| .
\end{equation}
\end{prop}

\par In order to consider the uniqueness in the proposition, let $u\in
\langle \cdot \rangle^\epsilon L^2$ satisfy (\ref{db.10}) with $v=0$,
so that $u$ is an entire function. Using the mean-value property,
$$
u(z)=\frac{1}{\pi }\int_{D(z,1)} u(w)L(dw),$$ 
where 
$$D(z,r)=\{ w\in \mathbb{C};\, |w-z|<r \},\ L(dw)=d\Re w\, d\Im w,
$$
we get by the Cauchy-Schwarz inequality,
$$
|u(z)|\le \pi ^{-1/2}\| u\|_{L^2(D(z,1))}\le o(1)\langle
z\rangle^\epsilon \|\langle \cdot \rangle^{-\epsilon }u\|_{L^2},
$$
and hence $u$ is a polynomial. A polynomial of
degree $\le m\in \mathbb {N}$ belongs to $\langle \cdot \rangle^\epsilon
L^2$ iff $2(m-\epsilon )<-2$, i.e.\ iff $m<\epsilon -1$. Thus, with
$\mathcal{N}$ denoting the nullspace of an operator, $\mathcal{
  N}(\overline{\partial })\cap \langle \cdot \rangle^\epsilon L^2$ is
equal to the space of polynomials of degree $<\epsilon -1$. In
particular this space is reduced to $0$ when $\epsilon \in ]0,1]$.
\begin{prop}\label{db2}
Let $0<\epsilon \le 1$. Then for every $v\in \langle \cdot
\rangle^{\epsilon -2}L^2$, the solution $u\in \langle \cdot
\rangle^\epsilon L^2$ of (\ref{db.10})(cf.\ Proposition \ref{db1}) is unique.
\end{prop}

\par Notice that it would have sufficed to state the last proposition
in the case of the largest spaces, i.e. in the case $\epsilon
=1$. With this value of $\epsilon $, let $v\in \langle \cdot 
\rangle^{-1}L^2$ and consider
\begin{equation}\label{db.12}
\widetilde{u}(z)=\frac{1}{\pi }\int \frac{1}{z-w}v(w)L(dw)=\widetilde{u}_1(z)+\widetilde{u}_2(z),
\end{equation}
where $\widetilde{u}_1(z)$ and $\widetilde{u}_2(z)$ are obtained by
inserting the factors $\chi (z-w)$ and $1-\chi (z-w)$ respectively in
the integral in (\ref{db.12}) and $\chi \in C_0^\infty (\mathbb{C})$ is
equal to 1 near $0$. We see that $\widetilde{u}_1(z)$ is well defined
since
$$
\langle z\rangle\widetilde{u}_1(z)=\frac{1}{\pi }\int \frac{\langle
  z\rangle}{\langle w\rangle}\frac{\chi (z-w)}{z-w}\underbrace{\langle
  w\rangle v(w)}_{\in L^2}L(dw)
$$
and
$$
\frac{\langle
  z\rangle}{\langle w\rangle}\frac{\chi (z-w)}{z-w}=\mathcal{
  O}(1)\frac{\chi (z-w)}{z-w},\hbox{ and }\frac{\chi (\cdot )}{\cdot
}\in L^1.
$$ 
It follows that $\widetilde{u}_1$ is well defined and
\begin{equation}\label{db.13}
\| \langle \cdot \rangle \widetilde{u}_1\| \le \mathcal{O}(1)\|\langle
\cdot \rangle v\|.
\end{equation}
$\widetilde{u}_2$ is well defined because of the Cauchy-Schwarz
inequality,
\begin{equation}\label{db.14}
\begin{split}
|\widetilde{u}_2(z)|&\le \frac{1}{\pi }\int \frac{|1-\chi
  (z-w)|}{|z-w|\langle w\rangle}\langle w \rangle |v(w)| L(dw)\\
&\le \frac{1}{\pi }\|\frac{|1-\chi
  (z-\cdot )|}{|z-\cdot |\langle \cdot \rangle}\|\|\langle \cdot
\rangle v\|\\
&\le \mathcal{O}(1)\frac{\langle \ln \langle z\rangle\rangle^{1/2}}{\langle
  z\rangle} \|\langle \cdot
\rangle v\| .
\end{split}
\end{equation}

Using the same decomposition, we can show that 
\begin{equation}\label{db.15}
\overline{\partial }\widetilde{u}=v
\end{equation}
and we have seen that we have (\ref{db.12}) with
$\widetilde{u}_1\in \langle \cdot \rangle^{-1}L^2$ and
$\widetilde{u}_2$ satisfying (\ref{db.14}).

\par For the same $v$, let $u\in \langle \cdot \rangle L^2$ be the
solution of $\overline{\partial }u=v$ in Propositions \ref{db1},
\ref{db2}. Then $\overline{\partial }(\widetilde{u}-u)=0$ and 
$$
\widetilde{u}-u=u_1+u_2,\hbox{ where }u_1=\widetilde{u}_1-u\in \langle
\cdot \rangle L^2,\ u_2=\widetilde{u}_2\to 0,\ z\to \infty .
$$
Following the discussion before Proposition \ref{db2}, we first see that
$\widetilde{u}-u$ is a polynomial, then that $\widetilde{u}-u=0$, 
and we obtain
\begin{prop}\label{db3}
Let $v\in \langle \cdot \rangle^{-1}L^2$ and let $u$ be the unique
solution in $\langle \cdot
\rangle L^2$ of $\overline{\partial }u=v$. Then
\begin{equation}\label{db.16}
u(z)=\frac{1}{\pi }\int \frac{1}{z-w}v(w) L(dw),
\end{equation}
where the integral is well defined according to the above discussion. 
\end{prop}

In the following we shall use the semi-classical calculus of
pseudodifferential operators, see e.g. \cite{DiSj99}. Let $\chi \in
C_0^\infty (\mathbb{R}^2_{\xi ,\eta })$, the space of smooth function 
with compact support on $\mathbb{R}^2_{\xi ,\eta }$. Then for $0<h\le 1$ we put 
$$\chi ^w=\mathrm{Op}(\chi )=\chi (hD),\ D=(D_x,D_y),\ z=x+iy$$
with the usual convention that $D_x=i^{-1}\partial _x$,
$D_y=i^{-1}\partial _y$. The exponent $w$ indicates that we use Weyl
quantization. The pseudodifferential operator calculus shows that if 
$\theta\in\mathbb{R}$, then
$\langle \cdot \rangle^{-\theta }\chi (hD)\langle \cdot \rangle^\theta
$ is again a pseudodifferential operator with a symbol of class
$S(1)$:
$$
\langle \cdot \rangle^{-\theta }\chi (hD)\langle \cdot
\rangle^\theta=a^w(x,hD),\ a\in S(1).
$$
Here, if $m>0$ denotes an order function (see \cite{DiSj99}), we let
$S(m)\subset C^\infty (\mathbb{R}^4)$ be the Fr\'echet space of all smooth
functions $a(x,y;\xi ,\eta   ;h)$ of $(x,y;\xi ,\eta )\in 
\mathbb{R}^2\times \mathbb{R}^2$ such that for all $\alpha ,\, \beta \in 
  \mathbb{N}^2$, there is a constant, $C=C_{\alpha ,\beta }$ such that
$$
|\partial _{x,y}^\alpha \partial _{\xi ,\eta }^\beta a(x,y;\xi ,\eta
;h)|\le Cm(x,y;\xi ,\eta), 
$$ 
uniformly with respect to $h$. (We may also need this definition for a
fixed value of $h$.) Here we use standard multiindex notation,
$$
\partial _{x,y}^\alpha =\partial _x^{\alpha _1}\partial _y^{\alpha
  _2},\ \alpha =(\alpha _1,\alpha _2)\in \mathbb{N}^2,
$$
and similarly for $\partial _{\xi ,\eta }^\beta $.
From the standard $L^2$ boundedness result for
pseudodifferential operators (here basically the
Cald\'eron-Vaillancourt theorem) we conclude that $\langle \cdot \rangle^{-\theta }\chi (hD)\langle \cdot
\rangle^\theta=\mathcal{O}(1):\, L^2\to L^2$ uniformly for $(\theta
,h)\in K\times ]0,1]$ for every fixed bounded interval $K$. It follows
that 
$$
\chi ^w=\chi (hD)=\mathcal{O}(1):\, \langle \cdot \rangle^\theta L^2\to
\langle \cdot \rangle^\theta L^2  
$$
with the same uniformity. 

\par In the situation of Propositions \ref{db1}, \ref{db2} we can
apply $\chi ^w$ to (\ref{db.10}) and get 
\begin{equation}\label{db.17}
h\overline{\partial }\chi ^wu=\chi ^wv,\ \chi^wu\in \langle
\cdot \rangle^\epsilon L^2,\ \chi ^wv\in \langle \cdot
\rangle^{\epsilon -2}L^2, 
\end{equation}
and from (\ref{db.11}) and uniqueness, we get
\begin{equation}\label{db.18}
\| \langle \cdot \rangle^{-\epsilon }\chi ^wu\|\le \epsilon ^{-1/2}\|
\langle \cdot \rangle^{2-\epsilon }\chi ^wv \|,
\end{equation}
for $u,v$ as in Propositions \ref{db1}, \ref{db2}.

Let $\xi ,\eta $ denote the dual variables to $x,y$ and notice that
the semi-classical symbols of $h\partial _{\overline{z}}$ and
$h\partial _{z}$ are equal to $i\overline{\zeta }$ and $i\zeta $
respectively, where \begin{equation}\label{db.18.5}\zeta :=(\xi -i\eta )/2,\ \overline{\zeta }=(\xi
+i\eta )/2. \end{equation} Notice that $x\cdot
\xi +y\cdot \eta =z\cdot \zeta +\overline{z}\cdot \overline{\zeta }$,
when $z=x+iy$.

\par If 
\begin{equation}\label{db.19}
0\notin\mathrm{supp\,}(1-\chi ),
\end{equation}
then $q=(1-\chi (\xi ,\eta ))/(i\overline{\zeta })$ is a smooth
function with $\partial _{\xi ,\eta }^\alpha q=\mathcal{O}(\langle \xi
,\eta \rangle^{-1-|\alpha |})$ (with the usual convention for
multiindices, that $|\alpha |=\|\alpha \|_{\ell^1}$), and from the
equation
\begin{equation}\label{db.20}
h\partial _{\overline{z}}(1-\chi ^w)u=(1-\chi ^w )v
\end{equation}
(still for $u,v$ as in Propositions \ref{db1}, \ref{db2}), we get
\begin{equation}\label{db.21}
(1-\chi ^w)u=\mathrm{Op}((1-\chi (\xi ,\eta ))/(i\overline{\zeta }))v.
\end{equation}
The pseudodifferential operator to the right is $\mathcal{O}(1):\ \langle
\cdot \rangle^\theta L^2\to \langle
\cdot \rangle^\theta L^2$ for every $\theta \in \mathbb{R}$. The apriori
estimate 
$$
\|\langle \cdot \rangle ^{-\epsilon }(1-\chi ^w)u\|\le \epsilon
^{-1/2}
\|\langle \cdot \rangle ^{2-\epsilon }(1-\chi ^w)v\|,
$$
(following from (\ref{db.20}),\ Propositions \ref{db1}, \ref{db2} and
the fact that $(1-\chi ^w)u\in \langle \cdot \rangle^\epsilon L^2$,
$(1-\chi ^w)v\in \langle \cdot \rangle^{\epsilon -2}L^2$,) improves
partially to
\begin{equation}\label{db.22}
\|\langle \cdot \rangle ^{2-\epsilon }(1-\chi ^w)u\|\le \mathcal{O}(1)
\|\langle \cdot \rangle ^{2-\epsilon }(1-\chi ^w)v\|\le \mathcal{O}(1)\|
\langle \cdot \rangle^{2-\epsilon }v\|.
\end{equation}
 
\section{Application to a $2\times 2$ system.}\label{sy}
\setcounter{equation}{0}

Let $q:\mathbb{C}\to \mathbb{C}$ and assume that for some $s\in ]1,2]$,
\begin{equation}\label{sy.1}
\langle z\rangle^2q\in H^s(\mathbb{C}). 
\end{equation}
This implies that $|\langle z\rangle^2q(z)|\le C_s\|\langle \cdot
\rangle q\|_{H^s}$ so that
\begin{equation}\label{sy.1.5}
|q(z)|\le \mathcal{ O}(1)\langle z\rangle^{-2}.
\end{equation}
We study the system (\ref{dbarpsi}),
\begin{equation}\label{sy.2}
\begin{cases}
\overline{\partial }\psi _1=(q/2)\psi _2,\\
\partial \psi _2=\sigma(\overline{q}/2)\psi _1,
\end{cases}
\end{equation}
with the condition that for some $k\in 
  \mathbb{C}$,
\begin{equation}\label{sy.3}
\psi _1=e^{kz}\phi _1,\ \psi _2=e^{\overline{k}\overline{z}}\phi _2,
\end{equation}
where
\begin{equation}\label{sy.4}
\phi _1=1+o(1),\ \phi _2=o(1),\ z\to \infty .
\end{equation}
The system (\ref{sy.2}) is equivalent to (\ref{dbarphi}):
\begin{equation}\label{sy.5}
\begin{cases}
\overline{\partial }\phi
_1=(q/2)e^{\overline{k}\overline{z}-kz}\phi _2,\\
\partial \phi
_2=\sigma(\overline{q}/2)e^{kz-\overline{k}\overline{z}}\phi _1.
\end{cases}
\end{equation}
We are interested in the case when $|k|$ is large and introduce the
semi-classical parameter $h=1/|k|$, $0<h\ll 1$. Then as we have seen
in (\ref{omega}),
\begin{equation}\label{sy.6}
kz-\overline{k}\overline{z}=\frac{i}{h}\Re (z\overline{\omega
})=\frac{i}{h}\langle z,\omega \rangle_{\mathbb{R}^2},
\end{equation}
where 
$$
\omega =2i\frac{\overline{k}}{|k|},\ \ |\omega |=2.
$$
After multiplication with $h$, (\ref{sy.5}) takes the equivalent form:
\begin{equation}\label{sy.7}
\begin{cases}
h\overline{\partial }\phi
_1=h(q/2)e^{-\frac{i}{h}(x,y)\cdot \omega }\phi _2,\\
h\partial \phi
_2=\sigma h(\overline{q}/2)e^{\frac{i}{h}(x,y)\cdot \omega}\phi _1.
\end{cases}
\end{equation}
To shorten the notation, put 
$$
\widehat{\tau }_\omega u(z)=e^{\frac{i}{h}(x,y)\cdot \omega }u(z),\
\widehat{\tau }_{-\omega } u(z)=(\widehat{\tau }_\omega)^{-1} u(z)=e^{-\frac{i}{h}(x,y)\cdot \omega }u(z).
$$
$\widehat{\tau }_\omega $ is translation by $\omega $ on the
$h$-Fourier transform side. Here we use the $h$-Fourier
transform,
\begin{equation}
	\mathcal{ F}_hu(\xi,\eta )=\iint e^{-i(x,y)\cdot (\xi ,\eta )/h}u(x,y)dxdy
	\label{sy.7.5}
\end{equation}

so that $\mathcal{ F}_1$ is the usual Fourier transform.
$\widehat{\tau } _{\pm \omega }$ commute with the
multiplications in the right hand sides in (\ref{sy.7}) and this
system takes the form
\begin{equation}\label{sy.8}
\begin{cases}
h\overline{\partial }\phi
_1=\widehat{\tau }_{-\omega }h(q/2)\phi _2,\\
h\partial \phi
_2=\sigma\widehat{\tau }_\omega h(\overline{q}/2)\phi _1.
\end{cases}
\end{equation}

\par For $v\in \langle \cdot \rangle^{-1}L^2$, let $u$,
$\widetilde{u}$ be the unique solutions in $\langle \cdot \rangle L^2$
of the equations
\begin{equation}\label{sy.9}
h\overline{\partial }u=v,\ h\partial \widetilde{u}=v,
\end{equation}
and write $u=Ev$, $\widetilde{u}=Fv$. (Since $\partial $ is the
complex conjugate of $\overline{\partial }$, the results of Section
\ref{db} apply also to $\partial $.) Then $E,F:\, \langle \cdot
\rangle^{-1}L^2\to \langle \cdot \rangle L^2$ are bounded operators,
which are also bounded $\langle \cdot \rangle ^{\epsilon -2}L^2\to
\langle \cdot \rangle^{\epsilon }L^2$ for $0<\epsilon \le 1$ and by
Propositions \ref{db1}, \ref{db2}, we have
\begin{equation}\label{sy.10}
\| E\|_{\mathcal{L}(\langle \cdot \rangle^{\epsilon -2}L^2,\langle \cdot
  \rangle ^\epsilon L^2)},\ \| F\|_{\mathcal{L}(\langle \cdot \rangle^{\epsilon -2}L^2,\langle \cdot
  \rangle ^\epsilon L^2)}\le 1/(h\sqrt{\epsilon }).
\end{equation}

\par As a first approximate solution to (\ref{sy.8}), (\ref{sy.4}), we take
\begin{equation}\label{sy.11}
\phi _1^0=1,\ \phi _2^0=0.
\end{equation}
Then
\begin{equation}\label{sy.12}
\begin{cases}
h\overline{\partial }\phi _1^0-\widehat{\tau }_{-\omega }\frac{hq}{2}\phi
_2^0=0,\\
h\partial \phi _2^0-\sigma\widehat{\tau }_\omega \frac{h\overline{q}}{2}\phi
_1^0=-\sigma\widehat{\tau }_\omega \frac{h\overline{q}}{2}\phi _1^0,
\end{cases}
\end{equation}
where the term to the right in the second equation,
$\psi _2^0:=-\sigma\widehat{\tau }_\omega (h\overline{q}/2)\phi _1^0$ is
viewed as a remainder. By (\ref{sy.1}) we have
\begin{equation}\label{sy.12.5}
q\in \langle \cdot \rangle^{\epsilon -2}L^2,\hbox{ for }\epsilon
\in ]0,1].
\end{equation}

\par In order to correct for the remainder in (\ref{sy.12}), we
consider the inhomogeneous system
\begin{equation}\label{sy.13}
\begin{cases}
h\overline{\partial }\phi _1-\widehat{\tau }_{-\omega
}\frac{hq}{2}\phi _2=f _1,\\
h\partial \phi _2-\sigma\widehat{\tau }_{\omega
}\frac{h\overline{q}}{2}\phi _1=f _2,
\end{cases}
\end{equation}
for $f _j\in \langle \cdot \rangle^{\epsilon -2}L^2$, $\phi _j\in
\langle \cdot \rangle^\epsilon L^2$ and $\epsilon\in
]0,1]$ as in (\ref{sy.12.5}). This is equivalent to
\begin{equation}\label{sy.14}
\begin{cases}
\phi _1-E\widehat{\tau }_{-\omega }\frac{hq}{2}\phi _2=Ef _1,\\
\phi _2-\sigma F\widehat{\tau }_\omega \frac{h\overline{q}}{2}\phi _1=Ff _2,
\end{cases}
\end{equation}
or 
\begin{equation}\label{sy.15}
(1-\mathcal{K})\begin{pmatrix}\phi _1\\\phi
  _2\end{pmatrix}=\begin{pmatrix}Ef _1\\Ff _2\end{pmatrix},
\end{equation}
where
\begin{equation}\label{sy.16}
\mathcal{K}=\begin{pmatrix}0 &A\\ B &0\end{pmatrix},\ \ 
\begin{cases}A=E\widehat{\tau }_{-\omega }\frac{hq}{2},\\
B=\sigma F\widehat{\tau }_\omega \frac{h\overline{q}}{2}\ (=\sigma\overline{A})
\end{cases} .
\end{equation}

\par Since $Ef _1,\, Ff _2\in \langle \cdot \rangle^\epsilon
L^2$ (with $\epsilon$ as in (\ref{sy.12.5})), we want to
invert 
$$
1-\mathcal{K}:\ (\langle \cdot \rangle^\epsilon L^2)^2\to (\langle
\cdot \rangle^\epsilon L^2)^2.
$$
Using also (\ref{sy.1.5}), we see that $A,B=\mathcal{O}(1):\, \langle \cdot \rangle^\epsilon L^2\to
\langle \cdot \rangle^\epsilon L^2$ ($\epsilon$ is fixed) and hence $\mathcal{K}=\mathcal{O}(1):\, (\langle \cdot \rangle^\epsilon L^2)^2\to
(\langle \cdot \rangle^\epsilon L^2)^2$ which is not enough to imply
the invertibility of $1-\mathcal{K}$ without a smallness condition on
$q$. Instead, we shall show that $\mathcal{K}^2$ is of small norm and
obtain the inverse of $1-\mathcal{K}$ as
\begin{equation}\label{sy.17}
(1-\mathcal{K})^{-1}=(1+\mathcal{K})(1-\mathcal{K}^2)^{-1}=(1-\mathcal{K}^2)^{-1}(1+\mathcal{K}).
\end{equation}

\par We have 
\begin{equation}\label{sy.18}
\mathcal{K}^2=\begin{pmatrix}AB &0\\ 0 &BA\end{pmatrix},
\end{equation}
where we recall that $A$, $B$ are given in (\ref{sy.16}). Let $\chi
\in C_0^\infty (\mathbb{R}^2)$ have its support in a small neighborhood of
$0$ and satisfy (\ref{db.19}). From (\ref{db.22}) and the adjacent
discussion we know that 
$$
E(1-\chi ^w)=(1-\chi ^w)E=\mathcal{O}(1):\, \langle \cdot \rangle^\theta
L^2\to \langle \cdot \rangle^\theta L^2,
$$
for every $\theta \in \mathbb{R}$ and similarly with $E$ replaced by
$F$. Write
\begin{equation}\label{sy.19}
  AB=\frac{h^2}{4}E\widehat{\tau }_{-\omega }\sigma qF\widehat{\tau }_\omega \overline{q}=\mathrm{I}+\mathrm{II}+\mathrm{III},
\end{equation}
where 
\begin{equation}\label{sy.19.5}
\begin{split}
  \mathrm{I}&=\frac{h^2}{4}E\widehat{\tau }_{-\omega }\sigma qF(1-\chi
  ^w)\widehat{\tau }_\omega \overline{q},\\
  \mathrm{II}&=\frac{h^2}{4}E(1-\chi ^w)\widehat{\tau }_{-\omega
  }\sigma qF\chi ^w\widehat{\tau }_\omega \overline{q},\\
\mathrm{III}&=\frac{h^2}{4}E\chi ^w\widehat{\tau }_{-\omega
  }\sigma qF\chi ^w\widehat{\tau }_\omega \overline{q}.
\end{split}
\end{equation}
We decompose each of the three terms into factors:

\medskip
\par\noindent \hskip -1cm\begin{tabular}{llclclclcl}
    &$\langle \cdot \rangle^\epsilon
 L^2$
&$\leftarrow$
&$\langle \cdot \rangle^{\epsilon -2}
 L^2$
&$\leftarrow$
& $\langle \cdot \rangle^\epsilon  L^2$
&$\leftarrow $
&$\langle \cdot \rangle^{\epsilon -2}L^2$
&$\leftarrow$
&$\langle
\cdot \rangle^\epsilon  L^2$\\
$\mathrm{I}:$
&
&$\frac{h^2E}{4}$
&
&$\sigma\widehat{\tau }_{-\omega }q$
& 
&$F(1-\chi ^w) $
&
&$\widehat{\tau }_\omega \overline{q}$
&\\
$\mathrm{II}:$
&
&$\frac{h^2E}{4}(1-\chi ^w)$
&
&$\sigma\widehat{\tau }_{-\omega }q$
& 
&$F\chi ^w $
&
&$\widehat{\tau }_\omega \overline{q}$
&\\
$\mathrm{III}:$
&
&$\frac{h^2E}{4}\chi ^w$
&
&$\sigma\widehat{\tau }_{-\omega }q$
& 
&$F\chi ^w $
&
&$\widehat{\tau }_\omega \overline{q}$
&
\end{tabular}

\medskip
\par Correspondingly we get for the operator norms:

\medskip
\[
\begin{split}
  \mathrm{I}&=\mathcal{O}(1)\, h\times 1\times 1\times 1=\mathcal{
    O}(h):\langle \cdot \rangle^\epsilon L^2\to \langle \cdot
  \rangle^\epsilon L^2,\\
\mathrm{II}&=\mathcal{O}(1) \, h^2\times 1\times h^{-1}\times 1=\mathcal{
    O}(h):\langle \cdot \rangle^\epsilon L^2\to \langle \cdot
  \rangle^\epsilon L^2,\\
\mathrm{III}&=\mathcal{O}(1) \, h\times 1\times h^{-1}\times 1=\mathcal{
    O}(1):\langle \cdot \rangle^\epsilon L^2\to \langle \cdot
  \rangle^\epsilon L^2.
\end{split}
\]

\medskip So far we only used (\ref{sy.1.5}).
In order to improve the estimate for $\mathrm{III}$, we write
\begin{equation}\label{sy.20}
\mathrm{III}=\sigma\frac{h^2}{4}E\left(\chi ^w\widehat{\tau }_{-\omega
  }q\chi ^w \right) F\widehat{\tau }_\omega \overline{q},
\end{equation}
using also that $F\chi ^w=\chi ^wF$.

\par Here,
\begin{equation}\label{sy.21}
\chi ^w\widehat{\tau }_{-\omega }q\chi ^w=\widehat{\tau }_{-\omega }
(\widetilde{\chi }^wq\chi ^w),
\end{equation}
where $\widetilde{\chi }^w=\widehat{\tau }_\omega \chi ^w\widehat{\tau
}_{-\omega }$ has the symbol
$\widetilde{\chi }(\xi ,\eta )=\chi ((\xi ,\eta )-\omega )\in
C_0^\infty (\mathbb{R}^2)$. If $\mathrm{supp\,}\chi $ is contained in a
sufficiently small neighborhood of $0\in \mathbb{R}^n$, then
$\mathrm{supp\,}\widetilde{\chi }$ will be contained in a small
neighborhood of $\omega $ and $\mathrm{supp\,}\chi \cap
\mathrm{supp\,}\widetilde{\chi }=\emptyset $.

\par Using the full assumption (\ref{sy.1}) we shall see that the
operator (\ref{sy.21}) or equivalently the operator $\widetilde{\chi
}^wq\chi ^w$ is $\mathcal{ O}(h^{s-1}):\, \langle \cdot \rangle^\epsilon
L^2\to \langle \cdot \rangle^{\epsilon -2}L^2$ for $\epsilon \in
]0,1]$ fixed. ($\widehat{\tau }_{-\omega }$ is unitary in $\langle
\cdot \rangle^{\epsilon -2}L^2$.)

\par Equivalently, we shall see that
\begin{equation}\label{sy.22}
\langle \cdot \rangle^{2-\epsilon }\widetilde{\chi }^wq\chi ^w\langle \cdot
\rangle^\epsilon  = \mathcal{ O}(h^{s-1}):\, L^2\to L^2.
\end{equation}
This operator can be written
\begin{equation}\label{sy.23}
\langle \cdot  \rangle^{2-\epsilon }\widetilde{\chi }^w\langle \cdot
\rangle^{\epsilon -2}\langle \cdot
\rangle^{2-\epsilon } q\langle \cdot \rangle^\epsilon \langle \cdot \rangle^{-\epsilon} \chi ^w\langle \cdot
\rangle^\epsilon =\widetilde{\Lambda }_{\epsilon -2}q_0\Lambda
_\epsilon ,
\end{equation}
where
\begin{equation}\label{sy.24}
\widetilde{\Lambda }_{\epsilon -2}=\langle \cdot  \rangle^{2-\epsilon }\widetilde{\chi }^w\langle \cdot
\rangle^{\epsilon -2},\ \Lambda _\epsilon =\langle \cdot
\rangle^{-\epsilon }\chi ^w \langle \cdot \rangle^\epsilon ,
\end{equation}
\begin{equation}\label{sy.25}
q_0=\langle \cdot \rangle^2q,
\end{equation}
and $q_0\in H^s$ by (\ref{sy.1}). In particular (as we have already
seen), $q_0\in L^\infty $, since $s>1$.

By semi-classical calculus (\cite{DiSj99}) we know that
$\widetilde{\Lambda }$ and $\Lambda $ are $h$-pseudo-differential
operators with symbols of class $S(1)$ which are of class $h^\infty
S(1)$ for $\xi $ away from any fixed neighborhood of
$\mathrm{supp\,}\widetilde{\chi }$ and $\mathrm{supp\,}\chi $
respectively. It follows that if $\widehat{\widetilde{\chi }},\,
\widehat{\chi }\in C_0^\infty $ are equal to $1$
near $\mathrm{supp\,}\widetilde{\chi }$ and $\mathrm{supp\,}\chi $
respectively, then
$$
\widetilde{\Lambda }-  \widetilde{\Lambda }\widehat{\widetilde{\chi
  }}^w ,\ \Lambda -\widehat{\chi }^w\Lambda =\mathcal{ O}(h^{\infty
  }):\, L^2\to L^2.
  $$

  Using also that $\Lambda ,\widetilde{\Lambda }$ are uniformly
  bounded: $L^2\to L^2$, it then suffices to
  show that
  \begin{equation}\label{sy.26}
\widehat{\widetilde{\chi }}^w q_0 \widehat{\chi }^w=\mathcal{
  O}(h^{s-1}):\, L^2\to L^2.
\end{equation}
Taking $\widehat{\widetilde{\chi }}$, $\widehat{\chi }$ with support
in sufficiently small neighborhoods of $\mathrm{supp\,}\widetilde{\chi
}$ and $\mathrm{supp\,}\chi $ respectively, we can also assume that
$\mathrm{supp\,}\widehat{\widetilde{\chi }}\cap
\mathrm{supp\,}\widehat{\chi }=\emptyset $. Then
$(\widehat{\widetilde{\chi }},\widehat{\chi })$ has the same
properties as $(\widetilde{\chi },\chi )$ and to simplify the
notation we can drop the hats and show that
\begin{equation}\label{sy.27}
A:=\widetilde{\chi }^wq_0\chi ^w=\mathcal{ O}(h^{s-1}):\, L^2\to L^2.
\end{equation}

\par Write $\widehat{u}=\mathcal{ F }_1u$ for the standard Fourier
transform ($h=1$ in (\ref{sy.7.5})). Then
\begin{equation}\label{sy.28}
\widehat{Au}(\xi )=\widetilde{\chi }(h\xi )\int \widehat{q}_0(\xi
-\eta )\chi (h\eta )\widehat{u}(\eta )\frac{d\eta }{(2\pi )^2}.
\end{equation}
By Cauchy-Schwarz, we have
$$
\left| \int \widehat{q}_0(\xi -\eta )\chi (h\eta )\widehat{u}(\eta )
  \frac{d\eta }{(2\pi )^2} \right| \le \sup |\chi |
\| \widehat{q}_0(\xi -\cdot )\|_{L^2(\mathrm{supp\,}\chi (h\cdot ))}\|\widehat{u}\|.$$
We have
$\langle \xi -\eta \rangle\asymp 1/h$
for $(\xi ,\eta ) \in \mathrm{supp\,}\widetilde{\chi } (h\cdot )\times
\mathrm{supp\,}\chi (h\cdot )$, so for $\xi \in
\mathrm{supp\,}\widetilde{\chi }(h\cdot )$ we have uniformly,
$$
\| \widehat{q}_0(\xi -\cdot )\|_{L^2(\mathrm{supp\,}\chi (h\cdot
  ))}\le \mathcal{ O}(h^s)\|\widehat{q}_0(\xi -\cdot )\langle \xi -\cdot
\rangle^s\|_{L^2}=
\mathcal{ O}(h^s)\| q_0\|_{H^s}.
$$
Using this in (\ref{sy.28}) gives
\begin{equation}\label{sy.29}
|\widehat{Au}(\xi )|\le |\widetilde{\chi }(h\xi )|\mathcal{ O}(h^s)\| q_0\|_{H^s}\|\widehat{u}\|,
\end{equation}
which implies
\begin{equation}\label{sy.30}
\| \widehat{Au}\|\le \mathcal{ O}(h^{s-1})\| q_0\|_{H^s}\|\widehat{u}\| ,
\end{equation}
since $\|\widetilde{\chi }(h\cdot )\|=h^{-1}\|\widetilde{\chi
}\|=\mathcal{ O}(1/h)$, and (\ref{sy.27}) follows. We have then
established (\ref{sy.22}) and as for $\mathrm{I}$, $\mathrm{II}$ we
can use the estimates for $E$, $F$ to conclude that
\begin{equation}\label{sy.31}
\mathrm{III}=\mathcal{ O}(h^{s-1}):\, \langle \cdot \rangle^\epsilon
L^2\to \langle \cdot \rangle^\epsilon L^2.
\end{equation}

\par We conclude that $AB=\mathcal{O}(h^{s-1}):\langle \cdot \rangle^\epsilon
L^2\to \langle \cdot \rangle^\epsilon L^2$, for every fixed
$\epsilon\in ]0,1]$. The same conclusion holds for $BA$, so 
\begin{equation}\label{sy.32}\mathcal{K}^2=\mathcal{O}(h^{s-1}):\, (\langle \cdot
  \rangle^\epsilon L^2)^2\to (\langle \cdot
  \rangle^\epsilon L^2)^2 .\end{equation}
 Then by (\ref{sy.17}), $(1-\mathcal{K})^{-1}$
exists and is $\mathcal{O}(1): (\langle \cdot \rangle^\epsilon
L^2)^2\to (\langle \cdot \rangle^\epsilon L^2)^2$, when $h>0$ is small enough.

\par Summing up we have 
\begin{prop}\label{sy1}
Let $q\in \langle \cdot \rangle^{-2}H^s$ for some $s\in ]1,2]$ and fix
$\epsilon \in ]0,1]$. Define $\mathcal{ K}$ as in (\ref{sy.16}). Then
$\mathcal{ K}=\mathcal{ O}(1): (\langle \cdot \rangle^\epsilon  L^2)^2\to
(\langle \cdot \rangle^\epsilon  L^2)^2 $,
$$
\mathcal{ K}^2=\mathcal{ O}(h^{s-1}):\, (\langle \cdot \rangle^\epsilon
L^2)^2\to (\langle \cdot \rangle^\epsilon  L^2)^2.
$$
For $h_0>0$ small enough and $0<h\le h_0$,
$1-\mathcal{ K}: (\langle \cdot \rangle^\epsilon  L^2)^2\to (\langle \cdot \rangle^\epsilon  L^2)^2$
has a uniformly bounded inverse,
\begin{equation}\label{sy.32.5}
  (1-\mathcal{ K})^{-1}=(1-\mathcal{ K}^2)^{-1}(1+\mathcal{ K})=
  \begin{pmatrix} (1-AB)^{-1} &0\\ 0
    &(1-BA)^{-1}\end{pmatrix} \begin{pmatrix}
1 &A\\ B &1
  \end{pmatrix},
\end{equation}
and the system (\ref{sy.13}) has a
unique solution $(\phi _1,\phi _2)\in (\langle \cdot \rangle^\epsilon
L^2)^2$ for every $(f _1,f _2)\in (\langle \cdot
\rangle^{\epsilon -2}
L^2)^2$, which also satisfies (\ref{sy.14}) and the uniform apriori
estimate,
$$
\| (\phi _1,\phi _2)\|_{ (\langle \cdot \rangle^\epsilon
L^2)^2}\le \mathcal{ O}(1/h)\| (f _1,f _2)\|_{(\langle \cdot
\rangle^{\epsilon -2}
L^2)^2}
$$
\end{prop}

We return to the problem (\ref{sy.8}), (\ref{sy.4}) and recall that we
have the approximate solution $(\phi _1^0,\phi _2^0)$ in
(\ref{sy.11}), which satisfies (\ref{sy.12}). We look for the full
solution in the form,
\begin{equation}\label{sy.33}
\phi _1=\phi _1^0+\phi _1^1,\ \phi _2=\phi _2^0+\phi _2^1,
\end{equation}
where $\phi _1^1$, $\phi _2^1$ should fulfill
\begin{equation}\label{sy.34}
  \begin{cases}
h\overline{\partial }\phi _1^1-\widehat{\tau }_{-\omega
}\frac{hq}{2}\phi _2^1=0,\\
h{\partial }\phi _2^1-\sigma\widehat{\tau }_{\omega
}\frac{h\overline{q}}{2}\phi _1^1=\sigma\widehat{\tau }_\omega \frac{h\overline{q}}{2}. 
  \end{cases}
\end{equation}
Thus $\phi _1^1$, $\phi _2^1$ should satisfy (\ref{sy.13}) with $f_1=0$
and $f_2=\widehat{\tau }_\omega h\overline{q}/2$ which is $=\mathcal{
  O}(h)$ in $\langle \cdot \rangle^{-2}H^s$. We look for $\phi _j^1$
in $\langle \cdot \rangle^\epsilon L^2$ for $\epsilon \in ]0,1]$ and
get the equivalent system (cf.\ (\ref{sy.14})),
\begin{equation}\label{sy.35}
  \begin{cases}
    \phi _1^1-E\widehat{\tau }_{-\omega }\frac{hq}{2}\phi _2^1=E0,\\
    \phi _2^1-\sigma F\widehat{\tau }_\omega \frac{h\overline{q}}{2}\phi
    _1^1=\sigma F\widehat{\tau }_\omega \frac{h\overline{q}}{2},
  \end{cases}
\end{equation}
i.e.
\begin{equation}\label{sy.36}
(1-\mathcal{ K})\begin{pmatrix}\phi _1^1\\ \phi
  _2^1\end{pmatrix}=\begin{pmatrix}0\\
\sigma F\widehat{\tau }_\omega \frac{h\overline{q}}{2}\end{pmatrix}.
\end{equation}
Here $F\widehat{\tau }_\omega h\overline{q}/2=\mathcal{O}(1)$ in $\langle
\cdot \rangle^\epsilon L^2$ and Proposition \ref{sy1} gives us a
unique solution in $(\langle \cdot \rangle^\epsilon L^2)^2$ which is
$\mathcal{O}(1)$ in that space. More precisely by (\ref{sy.32.5}),
\begin{equation}\label{sy.41}\begin{split}
\phi _1^1&=(1-AB)^{-1}A\sigma F\widehat{\tau }_\omega
\frac{h\overline{q}}{2}\\ &=A\sigma F\widehat{\tau }_\omega
\frac{h\overline{q}}{2}+\mathcal{ O}\left(h^{s-1}\|AF\widehat{\tau }_\omega
h\overline{q}\|_{\langle \cdot \rangle^\epsilon L^2}\right)\hbox{ in }\langle \cdot
\rangle^\epsilon L^2,
\end{split}
\end{equation}
\begin{equation}\label{sy.42}\begin{split}
\phi _2^1&=(1-BA)^{-1}\sigma F\widehat{\tau }_\omega
\frac{h\overline{q}}{2}\\ &=\sigma F\widehat{\tau }_\omega
\frac{h\overline{q}}{2}
+\mathcal{ O}(h^{s-1}\|F\widehat{\tau }_\omega
h\overline{q}\|_{\langle \cdot
\rangle^\epsilon L^2})\hbox{ in }\langle \cdot
\rangle^\epsilon L^2,
\end{split}
\end{equation}
where $F\widehat{\tau }_\omega
  h\overline{q}/2=B(1)$ and $A\sigma F\widehat{\tau }_\omega
  h\overline{q}/2=AB(1)$ are $\mathcal{ O}(1)$  in
$\langle \cdot \rangle^\epsilon L^2 $ by (\ref{sy.16}).

\par Here we  have used  that $AB,\ BA=\mathcal{ O}(h^{s-1}):\, \langle \cdot
\rangle^\epsilon L^2\to \langle \cdot
\rangle^\epsilon L^2$. 

\par In order to study $\phi _1^1$ in (\ref{sy.41}), we notice that
$A\sigma F\widehat{\tau }_\omega h\overline{q}/2=AB(1)$ and that we have the
same estimates for $AB(1)$ in ${\langle \cdot \rangle^\epsilon L^2}$
as the ones for
$AB:\, \langle \cdot \rangle^\epsilon L^2\to \langle \cdot
\rangle^\epsilon L^2 $. Indeed, in those estimates, starting with the
decomposition (\ref{sy.19}), we just have to replace the fact that
$\widehat{\tau }_\omega \overline{q}:\mathcal{ O}(1):\, \langle \cdot
\rangle^\epsilon L^2\to \langle \cdot \rangle^{\epsilon -2} L^2$ as a
multiplication operator with the fact that
$\widehat{\tau }_\omega \overline{q}$ is $\mathcal{ O}(1)$ in
$\langle \cdot \rangle ^{-2} L^2$ and a fortiori in
$\langle \cdot \rangle^{\epsilon -2} L^2$. We then get
$AF\widehat{\tau }_\omega h\overline{q}/2=\mathcal{ O}(h^{s-1})$ in
$\langle \cdot \rangle^\epsilon L^2$. Using this in (\ref{sy.41}), we get
\begin{equation}\label{sy.44}
\phi _1^1=A\sigma F\widehat{\tau }_\omega \frac{h\overline{q}}{2}+\mathcal{
  O}(h^{2(s-1)})=\mathcal{ O}(h^{s-1})\hbox{ in }\langle \cdot \rangle
^\epsilon L^2.
\end{equation}
  Similarly, from (\ref{sy.42}) we get
  $$
\phi _2^1=\sigma F\widehat{\tau }_\omega \frac{h\overline{q}}{2}+{\mathcal{O}}(h^{s-1})\hbox{ in }\langle \cdot \rangle^\epsilon L^2.
  $$

\begin{prop}\label{sy2} Recall that $q\in\langle \cdot 
	\rangle^{-2}H^{s}(\mathbb{C})$. Let $\epsilon \in ]0,1]$.
Problem (\ref{sy.8}) has a unique solution $\phi _1$, $\phi _2$ of
the form (\ref{sy.33}) with $\phi _1^1$, $\phi _2^1$ in $\langle \cdot
\rangle^\epsilon L^2$, satisfying (\ref{sy.41}), (\ref{sy.42}),
(\ref{sy.44}). For two different $\epsilon \in ]0,1]$
we get the same solutions when $h$ is small enough.
\end{prop}

\section{The case when $q$ is a characteristic function}\label{qcf}
\setcounter{equation}{0}

In this section, we
treat the case when 
\begin{equation}\label{qcf.1}
q=1_{\Omega }
\end{equation}
where $\Omega \Subset \mathbb{C}$ is a strictly convex open set with smooth
boundary. We can repeat the discussion in Section \ref{sy} without any
changes until the study of $\mathcal{K}^2$ in (\ref{sy.18}). We still
make the decomposition in (\ref{sy.19}), (\ref{sy.19.5}) and again
$\mathrm{I},\, \mathrm{II}=\mathcal{O}(h):\, \langle \cdot
\rangle^\epsilon L^2\to \langle \cdot \rangle^\epsilon L^2 $ for
$\epsilon \in ]0,1]$ fixed, while $\mathrm{III}$ will need a new
treatment. In view of (\ref{sy.20}) we have
\begin{equation}\label{qcf.2}
\|\mathrm{III}\|_{\mathcal{L}(\langle \cdot \rangle^\epsilon L^2, \langle \cdot \rangle^\epsilon L^2)}=\mathcal{O}(1)\| \chi ^w\widehat{\tau }_{-\omega }q\chi ^w
\|_{\mathcal{L}(\langle \cdot \rangle^\epsilon L^2, \langle \cdot
  \rangle^{\epsilon -2} L^2)},
\end{equation}
and in this section we shall show that 
\begin{equation}\label{qcf.3}
\chi ^w\widehat{\tau }_{-\omega }q\chi ^w=\mathcal{O}(h^{1/2}):\,
\langle \cdot \rangle^\epsilon L^2\to \langle \cdot \rangle^{\epsilon -2} L^2,
\end{equation}
so $\mathrm{III}=\mathcal{O}(h^{1/2}):\langle \cdot \rangle^\epsilon
L^2\to \langle \cdot \rangle^\epsilon L^2$, implying (cf.\
(\ref{sy.19}),
\begin{equation}\label{qcf.4}
AB=\mathcal{O}(h^{1/2}):\, \langle \cdot \rangle^\epsilon L^2\to \langle
\cdot \rangle L^2.
\end{equation}
Similarly we will
have
\begin{equation}\label{qcf.5}
BA=\mathcal{O}(h^{1/2}):\, \langle \cdot \rangle^\epsilon L^2\to \langle
\cdot \rangle L^2,
\end{equation}
and hence by (\ref{sy.18}),
\begin{equation}\label{qcf.6}
\mathcal{K}^2=\mathcal{O}(h^{1/2}):\ (\langle \cdot \rangle^\epsilon
L^2)^2\to (\langle \cdot \rangle^\epsilon L^2)^2.
\end{equation}

\par We claim that in order to show (\ref{qcf.3}) it suffices to show
that
\begin{equation}\label{qcf.7}
\chi ^w\widehat{\tau }_{-\omega }q\chi ^w=\mathcal{O}(h^{1/2}):\, L^2\to
L^2.
\end{equation}
We shall first show that (\ref{qcf.7}) implies (\ref{qcf.3}) and then
we will establish (\ref{qcf.7}). We basically showed this implication in Section
\ref{sy} (cf.\ (\ref{sy.26})), and here we give a variant of that
argument, exploiting that $q$ now has compact support.

(\ref{qcf.3}) is equivalent to:
\begin{equation}\label{qcf.8}
\langle \cdot \rangle^{2-\epsilon } \chi ^w\widehat{\tau }_{-\omega }
q \chi ^w  \langle \cdot \rangle^\epsilon =\mathcal{O}(h^{1/2}):\
L^2\to L^2. 
\end{equation}
Let $\psi \in C_0^\infty (\mathbb{C})$ be equal to 1 on
$\mathrm{supp\,}q$. We have by $h$-pseudodifferential calculus (\cite{DiSj99})
\begin{multline*}
  \langle \cdot \rangle^{2-\epsilon } \chi ^w\widehat{\tau }_{-\omega
  } q \chi ^w \langle \cdot \rangle^\epsilon = \langle \cdot
  \rangle^{2-\epsilon } \chi ^w\psi \widehat{\tau }_{-\omega }
  q \psi \chi ^w  \langle \cdot \rangle^\epsilon \\
  = \underbrace{\langle \cdot \rangle^{2-\epsilon }\psi }_{\mathcal{
      O}(1)}\underbrace{\chi ^w\widehat{\tau }_{-\omega }q\chi
    ^w}_{\mathcal{O}(h^{1/2})\hbox{\footnotesize by
      (\ref{qcf.7})}}\underbrace{\psi \langle \cdot \rangle^\epsilon
  }_{\mathcal{O}(1)}+\underbrace{\langle \cdot \rangle^{2-\epsilon }[\chi
    ^w,\psi ]\widehat{\tau }_{-\omega}q}_{\mathcal{O}(h^\infty
    )}\underbrace{\psi \chi ^w\langle \cdot \rangle^{\epsilon
    }}_{\mathcal{
      O}(1)}\\
  +\underbrace{\langle \cdot \rangle^{2-\epsilon }\psi \chi ^w}_{\mathcal{
      O}(1)}\underbrace{\widehat{\tau }_{-\omega }q[\psi ,\chi
    ^w]\langle \cdot \rangle^{\epsilon }}_{\mathcal{O}(h^\infty )},
\end{multline*}
where the underbraces indicate bounds on the $\mathcal{L}(L^2,L^2)$
norms (the notation $\mathcal{O}(h^{\infty})$ denotes  a quantity 
which is $\mathcal{O}(h^N )$ 
for every $N\geq 1$). Here we use that the symbol of $[\chi ^w,\psi ]$ is $\mathcal{
  O}(h^\infty )$ over a neighborhood of $\mathrm{supp\,}q$. Thus, we
have seen that (\ref{qcf.3}) follows from (\ref{qcf.7}).

We now turn to the proof of (\ref{qcf.7}). Extend the definition of
$\widehat{u}$ in Section \ref{sy} to the $h$-dependent case and define
$\widehat{u}=\widehat{u}_h$ for $h>0$ by
\begin{equation}\label{qcf.9}
\widehat{u}(\xi )=\mathcal{F}_hu(\xi )=\int_{\mathbb{R}^2}e^{-ix\cdot \xi
  /h}u(x)dx
\end{equation}
as the $h$-Fourier transform of $u$, so that $\mathcal{F}=\mathcal{F}_1$ is the
ordinary Fourier transform. The $h$-Fourier transform of a product of
two functions is equal to the convolution of the $h$-Fourier
transforms with respect to the measure $d\xi /(2\pi h)^2$. As in
Section \ref{sy} we have
$$
\| \chi ^w\widehat{\tau }_{-\omega }q\chi ^w\| =\|\widehat{\tau
}_\omega \chi ^w\widehat{\tau }_{-\omega }q\chi ^w\|
=\| \chi _\omega ^wq\chi ^w\| ,
$$
where
\begin{equation}\label{qcf.10}
\chi _\omega (\xi )=\chi (\xi -\omega ).
\end{equation}
Hence (cf.\ (\ref{sy.28})),
\begin{equation}\label{qcf.10.5}
  \mathcal{F}_h(\chi _\omega ^wq\chi ^wu)(\xi )=
  \int \chi (\xi -\omega )(\mathcal{F}_hq)(\xi -\eta )\chi (\eta
  )\widehat{u}(\eta )\frac{d\eta }{(2\pi h)^2}.
\end{equation}

\par We next evaluate
$$
\mathcal{F}_hq(\xi )=\int_\Omega e^{-ix\cdot \xi /h}dx.
$$
From
$$
-\frac{1}{|\xi |^2}(\xi _1hD_{x_1}+\xi _2hD_{x_2})\left(e^{-ix\cdot
    \xi /h} \right)=e^{-ix\cdot \xi /h}
$$
we get
$$
e^{-ix\cdot \xi /h}dx_1\wedge dx_2=\frac{ih}{|\xi |}d\left(e^{-ix\cdot
    \xi /h}\left(\frac{\xi _1}{|\xi |}dx_2-\frac{\xi _2}{|\xi |}dx_1 \right) \right).
$$
Identifying $dx_1\wedge dx_2$ with the Lebesgue measure $dx$ on 
$\mathbb{R}^2$, we get by Stokes' formula,
$$
\mathcal{F}_h q(\xi )=\frac{ih}{|\xi |}\int_{\partial \Omega }e^{-ix\cdot
  \xi /h}\left(\frac{\xi _1}{|\xi |}dx_2-\frac{\xi _2}{|\xi |}dx_1 \right).
$$
Parametrize $\partial \Omega $ by $x=\gamma (s)$ with positive
orientation and $|\dot{\gamma }(s)|=1$. Then
\begin{equation}\label{qcf.11}
\mathcal{F}_h q(\xi )=\frac{ih}{|\xi |}\int_0^L e^{-i\gamma (s)\cdot \xi
  /h}\frac{\xi }{|\xi |}\cdot n(s)ds,
\end{equation}
where $L$ is the length of $\partial \Omega $ and $n(s)=(\dot{\gamma
}_2(s),-\dot{\gamma } _1(s))$ is the exterior unit normal to $\Omega $ at
$\gamma (s)$. (We could have done the same calculations within the
complex formalism.)

\par Let $w_+(\xi )\in \partial \Omega $ be the point where $\xi /|\xi
|$ is equal to the exterior unit normal (``the North Pole'') and let
$w_- $ be the point where it is equal to the interior unit normal
(``the South Pole''). By stationary phase,
\begin{equation}\label{qcf.12}
\mathcal{F}_h q(\xi )=h^{3/2}\left( c_+(\xi ;h)e^{-iw_+(\xi )\cdot \xi /h}+c_-(\xi ;h)e^{-iw_-(\xi )\cdot \xi /h} \right),
\end{equation}
\begin{equation}\label{qcf.13}
c_\pm (\xi ;h)\sim c_\pm^0(\xi )+hc_\pm^1(\xi )+...,\ \ c^0_\pm\ne 0.
\end{equation}
in $C^\infty (W)$ for any fixed domain $W\Subset \mathbb{R}^2\setminus \{0
\}$.
(At first $c_{\pm}$ are determined up to $\mathcal{O}(h^{\infty})$, and 
there is a similar reminder in (\ref{qcf.12}). This reminder can be 
absorbed by modifying $c_{+}$ or $c_{-}$ by 
$\mathcal{O}(h^{\infty})$.)

Here,
\begin{equation}\label{qcf.14}
H(\xi ):=w_+(\xi )\cdot \xi 
\end{equation}
is the support function of $\Omega $, also given by $H(\xi
)=\sup_{x\in \Omega }x\cdot \xi $. Also,
\begin{equation}\label{qcf.15}
w_-(\xi )\cdot \xi =\inf_{x\in \Omega }x\cdot \xi =-H(-\xi ).
\end{equation}
We notice that
\begin{equation}\label{qcf.16}
\partial _{\xi }H(\xi )=w_+(\xi )\in \partial \Omega .
\end{equation}

\par It follows that the Hessian
\begin{equation}\label{qcf.17}
H''(\xi )=\frac{\partial w_+(\xi )}{\partial \xi }
\end{equation}
can be viewed as a linear map $T_\xi \mathbb{R}^2\to T_{w_+}(\partial \Omega )$ and is of rank 1.

\par Let $p\in C^\infty (\mathbb{R})$ be a real valued function which is
$>0$ in $\Omega $, $<0$ in $\mathbb{R}^2\setminus \Omega $ with $dp\ne 0$
on $\partial \Omega $. Then (\ref{qcf.16}) can be reformulated as the
eikonal equation
\begin{equation}\label{qcf.18}
p(H'(\xi ))=0.
\end{equation}

\par By (\ref{qcf.10.5}), $\chi _\omega ^wq\chi ^w$ is unitarily
equivalent to
$$
\widehat{u}\mapsto \int \chi _\omega (\xi )(\mathcal{F}_hq)(\xi -\eta
)\chi (\eta )\widehat{u}(\eta )\frac{d\eta }{(2\pi h)^2}
$$
and from (\ref{qcf.12})--(\ref{qcf.15}) we see that this operator can
be decomposed as
\begin{equation}\label{qcf.19}
\widehat{u}\mapsto A_+\widehat{u}+A_-\widehat{u},\end{equation}
where
\begin{equation}\label{qcf.20}
  A_\pm \widehat{u}(\xi )=
  \int \chi _\omega (\xi )c_\pm(\xi -\eta ;h)e^{\frac{i}{h}\phi _\pm
    (\xi ,\eta )}\chi (\eta )\widehat{u}(\eta )\frac{d\eta }{(2\pi )^2h^{1/2}},
\end{equation}
\begin{equation}\label{qcf.21}
\phi _\pm (\xi ,\eta )=\mp H(\pm (\xi -\eta )).
\end{equation}
Clearly, the problems of estimating the $L^2$ boundedness of $A_+$ and
of $A_-$ are equivalent and in the following we shall only handle
\begin{equation}\label{qcf.22}
A:=A_+,
\end{equation}
\begin{equation}\label{qcf.23}
Au(\xi )=h^{-1/2}\int \chi _\omega (\xi )c(\xi -\eta ;h)e^{\frac{i}{h}\phi
  (\xi ,\eta )}\chi (\eta )u(\eta  ) d\eta ,
\end{equation}
\begin{equation}\label{qcf.24}
\phi (\xi ,\eta )=-H(\xi -\eta ),\ c=(2\pi )^{-2}c_+.
\end{equation}
Here we write $u$ instead of $\widehat{u}$ since we will work entirely
on the Fourier transform side.

\par We choose the support of $\chi $ contained in a small enough
neighborhood of $0$, so that
\begin{equation}\label{qcf.25}
\mathrm{supp\,}\chi \cap \mathrm{supp\,}\chi _\omega =\emptyset ,
\end{equation}
and hence $\xi -\eta \ne 0$ on $\mathrm{supp\,}(\chi _\omega (\xi
)\chi (\eta ))$.

\par We work in the canonical coordinates $(\xi ,\xi ^*)$ with
symplectic form $\sigma =\sum d\xi _j^*\wedge d\xi _j$, thinking of
$\xi $ as the base variables. With $x=-\xi ^*$, we get
$$
\sigma =\sum -dx_j\wedge d\xi _j =\sum d\xi _j\wedge dx_j 
$$
which is the usual symplectic form on $\mathbb{R}^2_x\times \mathbb{R}^2_\xi
$.

\par We view $A$ as a Fourier integral operator with canonical
relation,
\begin{equation}\label{qcf.26}
(\eta ,-\partial _\eta \phi (\xi ,\eta ))\mapsto (\xi ,\partial _\xi
\phi (\xi ,\eta ))
\end{equation}
for $\eta\in \mathrm{neigh\, supp\,}\chi $, $\xi \in \mathrm{neigh\,
  supp\,}\chi _\omega $. This restriction on $(\xi ,\eta )$ will be
kept below though not constantly recalled. By (\ref{qcf.24}) this
becomes
\begin{equation}\label{qcf.27}
(\eta ,-H'(\xi -\eta ))\mapsto (\xi ,-H'(\xi -\eta )),
\end{equation}
hence by (\ref{qcf.16}), we get
\begin{equation}\label{qcf.28}
(\eta ,-w_+(\xi -\eta ))\mapsto (\xi ,-w_+(\xi -\eta ))
\end{equation}
or equivalently,
\begin{equation}\label{qcf.29}
\mathrm{neigh\,supp\,}\chi \times (-\partial \Omega )\ni (\eta ,\eta
^*)\mapsto (\xi ,\eta ^*)\in \mathrm{neigh\,supp\,}\chi_\omega  \times (-\partial \Omega ),
\end{equation}
with
\begin{equation}\label{qcf.30}
\xi -\eta \in \mathbb{R}_+n(-\eta ^*),
\end{equation}
where $n(-\eta ^*)$ is the exterior unit normal of $\Omega $ at $-\eta
^*$.

\par Identifying $(\xi ,\xi ^*)$ with $(x,\xi )=(-\xi ^*,\xi )$, the
canonical relation becomes
\begin{equation}\label{qcf.31}
  \partial \Omega \times \mathrm{neigh\, supp\,}\chi
\ni (y,\eta )\mapsto (y,\xi )\in
  \partial \Omega \times \mathrm{neigh\, supp\,}\chi _\omega ,
\end{equation}
with
\begin{equation}\label{qcf.32}
\xi -\eta \in \mathbb{R}_+n(y).
\end{equation}
We can also describe the canonical relation by (\ref{qcf.31}) with
\begin{equation}\label{qcf.33}
(y,\xi )\in \exp (\mathbb{ R}_+H_p)(y,\eta )
\end{equation}
which is equivalent to (\ref{qcf.32}).

By means of a finite partition of unity we can decompose $A$ in
(\ref{qcf.23}) into a finite sum of operators
\begin{equation}\label{qcf.34}
A_{\xi _0,\eta_0}u(\xi )=h^{-1/2}\int\chi _{\xi _0}(\xi )c(\xi -\eta
;h)e^{\frac{i}{h}\phi (\xi ,\eta )}\chi _{\eta _0}(\eta )u(\eta )d\eta ,
\end{equation}
where $(\xi _0,\eta _0)$ take finitely many values in
$\mathrm{supp\,}\chi _\omega \times 
\mathrm{supp\,}\chi $ and $\chi _{\xi _0}$, $\chi _{\eta _0}$ are
$C_0^\infty $ cutoffs, supported in small neighborhoods of $\xi _0$
and $\eta _0$ respectively. For a given $(\xi _0,\eta _0)$ we make an
orthogonal change of the $x$-coordinates and the corresponding change
of dual variables, so that
\begin{equation}\label{qcf.35}
\xi _0-\eta _0=t_0(0,1)\hbox{ for some }t_0>0.
\end{equation}
We already know that $\phi ''_{\xi , \eta }(\xi _0,\eta _0)=H''(\xi
_0-\eta _0) $ is of rank
1. In the chosen coordinates, we get more precisely that
\begin{equation}\label{qcf.36}
\phi ''_{\xi , \eta }(\xi _0,\eta _0)=\begin{pmatrix}a_0 &0\\ 0
  &0\end{pmatrix},\ a_0\ne 0,
\end{equation}
or in other terms that
\begin{equation}\label{qcf.37}
\phi ''_{\xi _1,\eta _1}\ne 0,\ \phi ''_{\xi_1 ,\eta_2 }=\phi
''_{\xi_2 ,\eta_1 }=\phi ''_{\xi_2 ,\eta_2 }=0\hbox{ at }(\xi _0,\eta _0). 
\end{equation}

By choosing the supports of $\chi _{\xi _0}$, $\chi _{\eta _0} $
contained in small neighborhoods of $\xi _0$, $\eta _0$ we achieve
that on $\mathrm{supp\,}(\chi _{\xi _0}(\xi )\chi _{\eta _0}(\eta ))$:
\begin{equation}\label{qcf.38}
\phi ''_{\xi _1,\eta _1}-a_0,\ \phi ''_{\xi_1 ,\eta_2 },\ \phi
''_{\xi_2 ,\eta_1 },\ \phi ''_{\xi_2 ,\eta_2 }\hbox{ are small.}
\end{equation} 

In order to shorten the notation, write $A=A_{\xi _0,\eta _0}$,
$$
b(\xi ,\eta ;h)=\chi _{\xi _0}(\xi )c(\xi -\eta ;h)\chi _{\eta
  _0}(\eta ).
$$
Then from (\ref{qcf.34}), we can write
\begin{equation}\label{qcf.39}
  \begin{split}
Au(\xi )=&h^{1/2}\int \int e^{\frac{i}{h}\phi (\xi ,\eta )}b(\xi ,\eta
;h)u(\eta _1,\eta _2)\frac{d\eta _1}{h}d\eta _2\\
=&h^{1/2}\int (A(\xi _2,\eta _2)u(\cdot ,\eta _2))(\xi _1) d\eta _2.
  \end{split}
\end{equation}
From (\ref{qcf.38}) we know that $A(\xi _2,\eta _2)$ is a Fourier
integral operator in one dimension, associated to the canonical {\it
  transformation}
\begin{equation}\label{qcf.40}
(\eta _1,-\partial _{\eta _1}\phi (\xi ,\eta ))\mapsto (\xi
_1,\partial _{\xi _1}\phi (\xi ,\eta )).
\end{equation}
Moreover $A(\xi _2,\eta _2)$ is of order 0, so it follows that
\begin{equation}\label{qcf.41}
\| A(\xi _2,\eta _2)\|\le \mathcal{O}(1),
\end{equation}
where
$$\| A(\xi _2,\eta _2)\|=\| A(\xi _2,\eta _2)\|_{L^2(\mathbb{R}_{\eta_1})\to L^2(\mathbb{R}_{\xi _1})}.$$

By Fubini's theorem, if $v=v(\eta _1,\eta _2)$, then
$$
\| v\|_{L^2(\mathbb{R}^2)}=\|\|v(\cdot ,\cdot \cdot )\|_{L^2_{\eta _1}}\|_{L^2_{\eta _2}}.
$$
By (\ref{qcf.39}) and the triangular inequality for integrals,
\begin{equation}\label{qcf.42}
  \begin{split}
    \| h^{-1/2}Au(\xi _1,\xi _2)\|_{L^2_{\xi _1}}&\le \int \| A(\xi
    _2,\eta _2)u(\cdot ,\eta _2)\|_{L^2_{\xi _1}} d\eta
    _2\\
    &\le \int \underbrace{\| A(\xi _2,\eta _2)\|}_{=:K(\xi _2,\eta
      _2)}\|u(\cdot ,\eta _2)\|_{L^2_{\eta _1}} d\eta _2.
  \end{split}
\end{equation}

\par Let $M$ be the $L^2_{\eta _2}\to L^2_{\xi _2}$ norm of the
integral operator with kernel $K(\xi _2,\eta _2)$. Since $K$ has
compact support and is bounded by (\ref{qcf.41}), we have $M=\mathcal{
  O}(1)$. From (\ref{qcf.42}), we get
\begin{multline*}
\|h^{-1/2}Au\|_{L^2}=\|\| h^{-1/2}Au(\xi _1,\xi _2)\|_{L^2_{\xi
    _1}}\|_{L^2_{\xi _2}}\\
\le M\|\| u(\eta _1,\eta _2)\|_{L^2_{\eta _1}}\|_{L^2_{\eta _2}}=M\| u\|_{L^2}.
\end{multline*}
Thus,
\begin{equation}\label{qcf.43}
h^{-1/2}A=\mathcal{ O}(1):\mathcal{ O}(1):L^2\to L^2,
\end{equation}
and we get
$$
\| A_+\|_{L^2\to L^2}=\mathcal{ O}(h^{1/2}).
$$
Similarly, $\| A_-\|_{L^2\to L^2}=\mathcal{ O}(h^{1/2})$ and putting
things together, we get (\ref{qcf.7}) and hence (\ref{qcf.3}).\hfill{$\Box$}

Summing up, we have proved,
\begin{prop}\label{qcf1}
Let $q$ in (\ref{qcf.1}) be the characteristic function of a strictly
convex open set $\Omega \Subset \mathbb{C}$ with smooth boundary and fix
$\epsilon \in ]0,1]$. Then the conclusions of the propositions
\ref{sy1} and \ref{sy2} hold with $s=3/2$.  
\end{prop}

  \section{Study of an integral in the complex domain}\label{ic}
  \setcounter{equation}{0}
  
%   In the previous section it was shown that the 
%   solution to the Dirac system (\ref{dbarphi}) for small \( h=1/|k| 
%   \) is given for \( q \) being the characteristic function of a 
%   simply connected compact domain in \( \mathbb{C} \) with smooth 
%   boundary by the asymptotic behavior (\ref{Phisasym}) plus 
%   corrections of order \( \mathcal{O}(h^{1/2}) \). Thus 
%   in lowest order of \( h \) the main contribution is given by the 
%   solid Cauchy transform on the domain. The resulting
%   integral is computed in leading order of \( h \) in this section.

\par In this section we keep the assumptions of Section \ref{qcf}, and
we will study the function
\begin{equation}\label{ic0.1}
F\widehat{\tau }_\omega \frac{h\overline{q}}{2}=:\frac{1}{2\pi}\overline{f(z,k)}
\end{equation}
that appears as the leading term in the expansion of $\phi _2^1$ in
(\ref{sy.42}). This will also lead to some additional information on the 
expressions (\ref{sy.41}), (\ref{sy.44}) of $\phi^{1}_{1}$. Here $E$, $F$ are the inverses of
$h\overline{\partial }$, $h\partial $ respectively, acting on $\langle
\cdot \rangle L^2$, and we recall from Proposition \ref{db3} that
\begin{equation}\label{ic0.2}
Ev(z)=\frac{1}{\pi h}\int \frac{1}{z-w}v(w)L(dw).
\end{equation} 
Since $F$ is the complex conjugate, we have
\begin{equation}\label{ic0.3}
Fv(z)=\frac{1}{\pi h}\int \frac{1}{\overline{z}-\overline{w}}v(w)L(dw).
\end{equation} 
The function (\ref{ic0.1}) is therefore equal to 
\begin{equation}\label{ic0.4}
F\widehat{\tau }_\omega \frac{h\overline{q}}{2}(z)=\frac{1}{2\pi
}\int_\Omega \frac{1}{\overline{z}-\overline{w}}e^{kw-\overline{kw}}L(dw)=\frac{1}{2\pi }\overline{f(z,k)},
\end{equation}
where we also used (\ref{sy.6}) and the adjacent discussion. For
notational reasons, we will mainly deal with the complex conjugate,
\begin{equation}\label{ic.0}
f(z,k)=\int_\Omega \frac{1}{z-w}e^{\overline{kw}-kw}L(dw)=\iint_{\Omega
}\frac{e^{\overline{kw}-kw}}{z-w}\frac{d\overline{w}\wedge dw}{2i},
\end{equation}
where $\Omega \Subset \mathbb{C}$ is simply connected with smooth 
boundary. Here $k\in \mathbb{C}$, $|k|\gg 1$. Notice that $\frac{1}{2i}d\overline{w}\wedge dw=d\Re
w\wedge d\Im w$, and we identify this real 2-form with the Lebesgue
measure $L(dw)$.

The main result of this section is given in Theorem \ref{ic2} below. 

\subsection{ Reduction with Stokes' formula.} We observe that in the
sense of differential forms
\[
d_w\left(e^{\overline{kw}-kw} dw \right)=\frac{\partial }{\partial
  \overline{w}}\left(e^{\overline{kw}-kw} \right) d\overline{w}\wedge
dw=\overline{k}e^{\overline{kw}-kw}d\overline{w}\wedge dw.
\]
Inserting the factor $1/(z-w)$, we get
\begin{multline*}
d_w\left(\frac{1}{z-w}e^{\overline{kw}-kw} dw \right)=\frac{\partial
}{\partial \overline{w}}\left(e^{\overline{kw}-kw} \frac{1}{z-w}
\right) d\overline{w}\wedge dw\\
=\left(\frac{\overline{k}}{z-w} -\pi \delta _z(w) \right)
e^{\overline{kw}-kw}d\overline{w}\wedge dw,
\end{multline*}
where $\delta _z(w)$ denotes the delta function at $w=z$.
Hence,
$$
\left(\frac{1}{z-w}-\frac{\pi \delta _z(w)}{\overline{k}} \right)
e^{\overline{kw}-kw}d\overline{w}\wedge dw=
d_w\left(\frac{1}{\overline{k}(z-w)} e^{\overline{kw}-kw} dw
\right),
$$
and Stokes' formula gives after multiplication with $1/(2i)$:
$$
\frac{1}{2i\overline{k}}\int_{\partial \Omega
}\frac{1}{z-w}e^{\overline{kw}-kw}dw=
\iint_\Omega \frac{e^{\overline{kw}-kw}}{z-w}
\frac{d\overline{w}\wedge dw}{2i}-\begin{cases}
  0\hbox{ if }z\not\in \Omega ,\\
  \frac{\pi}{\overline{k}}e^{\overline{kz} -kz},\hbox{ if }z\in \Omega .
\end{cases}
$$
Here we assume that $z\not\in \partial \Omega $. Thus,
\begin{equation}\label{ic.1}
f(z,k)=\frac{1}{2i\overline{k}}\int_{\partial \Omega }\frac{1}{z-w}e^{\overline{kw}-kw}dw
+(\pi/\overline{k})e^{\overline{kz}-kz}1_\Omega(z) .
\end{equation}

\medskip
\subsection{Holomorphic extensions from a real-analytic curve.} The
task is now to study the integral in (\ref{ic.1}), and for that it
will be convenient to add the assumption on $\Omega $ that
\begin{equation}\label{iche.0}
\partial \Omega \hbox{ is a real analytic curve.}
\end{equation}
In other words, we assume that $\Omega \Subset \mathbb{C}$ is given by
$g<0$, where $g$ is real and smooth in a neighborhood of
$\overline{\Omega }$, real-analytic near $\partial \Omega $ and 
$g=0$ and $\nabla g\ne 0$ on $\partial \Omega $. Let $\gamma =\partial \Omega $ be the
oriented boundary of $\Omega $.

Thanks to the analyticity assumption (\ref{iche.0}), we have a
holomorphic extension of the function
$$iu_0(w)=kw-\overline{kw}=i|k|\Re (w\overline{\omega })$$
to
a neighborhood of $\partial \Omega $. (Without the analyticity
assumption on $\partial \Omega $ our study would undoubtedly go
through with minor changes, using an almost holomorphic extension of
$u_0$, i.e.\ a smooth extension whose anti-holomorphic derivative
vanishes to infinite order on $\gamma $.)

\par One way of constructing the holomorphic extension of $u_0$, that
we shall not follow, is to use the antiholomorphic involution
$\iota =\iota _{\partial \Omega }:\, \mathrm{neigh\,}(\partial \Omega
)\to \mathrm{neigh\,}(\partial \Omega )$, characterized by
\begin{equation}\label{iche.1}
\iota \hbox{ is anti-holomorphic},
\end{equation}
\begin{equation}\label{iche.2}
{{\iota }_\vert}_{\partial \Omega }=\mathrm{id}.
\end{equation} $\iota $ can be constructed as follows:
Let $G(z,w)$ be the polarization of $g(z)$, i.e.\ the unique
holomorphic function defined near the anti-diagonal, $\{
(z,\overline{z});\, z\in \mathrm{neigh\,}(\partial \Omega ,\mathbb{C})\}$,
such that
\begin{equation}\label{iche.3}
G(z,\overline{z})=g(z).
\end{equation}
Then $w=\iota (z)$ is given by
\begin{equation}\label{iche.4}
G(z,\overline{w})=0. 
\end{equation}
Now the holomorphic extension of $iu_0$ is given by
\begin{equation}\label{iche.5}
iu(w)=kw-\overline{k\iota (w)},\ w\in \mathrm{neigh\,}(\partial \Omega
,\mathbb{C}).
\end{equation}

\par However, for the practical computations, we choose a more direct method. Let $\partial
\Omega $ be parametrized by
$$
\mathbb{R}/L\mathbb{Z}\ni t\mapsto \gamma (t)\in \mathbb{C},\ L=|\partial
\Omega |,
$$
where $\gamma $ is real-analytic and (for simplicity) $|\dot{\gamma
}(t)|=1$. We parametrize points $w$ in a neighborhood of $\partial
\Omega $ by
\begin{equation}\label{iche.6}
w=\gamma (t)+is\dot{\gamma }(t),\ t\in \mathbb{R}/L\mathbb{Z},\ s\in
\mathrm{neigh\,}(0,\mathbb{R}).
\end{equation}
Notice that $i\dot{\gamma }(t)$ is the interior unit normal to
$\partial \Omega $, since $\gamma $ is positively oriented (so that
$\gamma(t) $ travels along $\partial \Omega $ in the anti-clockwise
direction). We express holomorphicity with respect to $w$ by means of a
$\overline{\partial }$ equation in $t,s$: From (\ref{iche.6}), we get
\[
dw=\underbrace{(\dot{\gamma }(t)+is\ddot{\gamma }(t))}_a
dt+\underbrace{i\dot{\gamma }(t)}_bds.
\]
The conjugate equation is 
$$
d\overline{w}=\overline{a}dt+\overline{b}ds
$$
and inverting this system of two equations, we get with
$c=a\overline{b}-\overline{a}b\ne 0$:
\begin{equation}\label{iche.7}
dt=\frac{1}{c}\left(-i\overline{\dot{\gamma }(t)}dw-i\dot{\gamma }(t) d\overline{w}\right),
\end{equation}
\begin{equation}\label{iche.8}
ds=\frac{1}{c}\left(-(\overline{\dot{\gamma }(t)+is\ddot{\gamma
  }(t)})dw+(\dot{\gamma }(t)+is\ddot{\gamma }(t))d\overline{w}\right).
\end{equation}
By abuse of notation, we write $u(w)=u(s,t)$. Then 
\[\begin{split}\label{}
  cdu=&(\partial _tu)cdt+(\partial _su) cds\\
=&(-i\overline{\dot{\gamma }}\partial _tu-(\overline{\dot{\gamma
  }+is\ddot{\gamma }})\partial _su)dw+\\
&(-i\dot{\gamma }\partial _tu+(\dot{\gamma }+is\ddot{\gamma })\partial _su)d\overline{w}.
\end{split}
\]
From this we see that $u$ is a holomorphic function of $w$ near
$\partial \Omega $ iff
\begin{equation}\label{iche.9}
\partial _tu+i\left(1+is\frac{\ddot{\gamma }}{\dot{\gamma }}
\right)\partial _su=0,\hbox{ near }s=0.\end{equation}

\par If $u_0$ is a given real-analytic function on $\partial \Omega $,
we write $u_0=u_0(t)$ by abuse of notation. Let $u$ be the unique
holomorphic extension to a neighborhood of $\partial \Omega $ and
write
\begin{equation}\label{iche.10}
  u=u(s,t)=u_0(t)+su_1(t)+s^2u_2(t)+\mathcal{O}(s^3).
\end{equation}
The Taylor coefficients $u_1$, $u_2$, ... can be determined from
(\ref{iche.9}), that we first rewrite as
\begin{equation}\label{iche.11}
\partial _su=i\left(1-is\frac{\ddot{\gamma }}{\dot{\gamma }}+\mathcal{
    O}(s^2) \right) \partial _tu.
\end{equation}
Substitution of (\ref{iche.10}) gives
$$
u_1(t)+2su_2(t)=i\left(1-is\frac{\ddot{\gamma }}{\dot{\gamma }} \right)
\left(\partial _tu_0+s\partial _tu_1)+\mathcal{O}(s^2)
\right) ,
$$
leading to
\[
  \begin{split}
    u_1&=i\partial _tu_0,\\
    u_2&=\frac{1}{2}\left(\frac{\ddot{\gamma }}{\dot{\gamma
        }}\partial _tu_0-\partial _t^2 u_0 \right).
  \end{split}
\]
Thus for the holomorphic extension (\ref{iche.10}), we have
\begin{equation}\label{iche.12}
u(s,t)=u_0(t)+is\partial _tu_0+\frac{s^2}{2}\left(\frac{\ddot{\gamma
    }}{\dot{\gamma }}\partial _tu_0-\partial _t^2u_0 \right)+\mathcal{O}(s^3),
\end{equation}
where we recall that $\gamma =\gamma (t)$.

\par Let now $iu_0(w)$ be the restriction to $\partial \Omega $ of
\begin{equation}\label{iche.13}
kw-\overline{kw}=i|k|\Re (w\overline{\omega }),\hbox{ so } u_0(w)=u_0(w,k)=|k|\Re
(w\overline{\omega  }),\ w\in \partial \Omega ,
\end{equation}
and write $u_0(t)=u_0(\gamma (t))$.  Here, we recall that
$\Re (w\overline{\omega })=\langle w,\omega \rangle_{\mathbb{R}^2}$. Write,
\begin{equation}\label{iche.14}
u_0(t)=\frac{|k|}{2}\left(\gamma (t)\overline{\omega }
    +\overline{\gamma }(t) \omega  \right)=|k|\langle \gamma
  (t),\omega \rangle_{\mathbb{R}^2},
\end{equation}
\begin{equation}\label{iche.15}
\partial _tu_0(t)=\frac{|k|}{2}\left(\dot{\gamma }(t)\overline{\omega }
    +\overline{\dot{\gamma }}(t) \omega  \right)=|k|\langle \dot{\gamma
  }(t),\omega \rangle_{\mathbb{R}^2},
\end{equation}
\begin{equation}\label{iche.16}
\partial^2 _tu_0(t)=\frac{|k|}{2}\left(\ddot{\gamma }(t)\overline{\omega }
    +\overline{\ddot{\gamma }}(t) \omega  \right)=|k|\langle \ddot{\gamma
  }(t),\omega \rangle_{\mathbb{R}^2}.
\end{equation}

\par From (\ref{iche.15}) we see that $\gamma (t)$ is a critical point
of $u_0(t)$ iff $\omega $ (which is non-vanishing) is normal to
$\partial \Omega $ at $\gamma (t)$. From $\langle \dot{\gamma }(t),
\dot{\gamma }(t)\rangle =1$ we know that
\begin{equation}\label{iche.17}
\langle \dot{\gamma }(t),\ddot{\gamma }(t)\rangle =0
\end{equation}
and hence $\ddot{\gamma }(t)$ is normal to $\partial \Omega $
everywhere. Thus at a critical point of $u_0$ we have $\ddot{\gamma
}(t)\in \mathbb{R}\omega $. It follows from (\ref{iche.16}) that a critical
point is nondegenerate precisely when $\ddot{\gamma }(t)\ne 0$, i.e.\
when $\partial \Omega $ has non-vanishing curvature there. Such a
point is
\begin{itemize}
  \item a local maximum if $\ddot{\gamma }(t)=c\omega $, $c<0$, and
  \item a local minimum if $\ddot{\gamma }(t)=c\omega $, $c>0$.
\end{itemize}

\par Now recall the assumption that
\begin{equation}\label{iche.18}
\Omega \hbox{ is strictly convex}.
\end{equation}
Then at every point in $\partial \Omega $, $\ddot{\gamma }(t)$ is
non-vanishing and of the form $c(t)\nu (t)$, where $c(t)>0$ and
$\nu (t)=i\dot{\gamma }(t)$ is the interior unit normal. (Recall that
$\gamma $ is positively oriented.)

\par For a fixed $k\ne 0$, we can
decompose
\begin{equation}\label{iche.19}
\partial \Omega =\{ w_-(k )\}\cup \Gamma _+\cup \{ w_+(k )\}\cup
\Gamma _-,
\end{equation}
ordered in the positive direction when starting and ending at
$w_-(k )$. Here
\begin{itemize}
\item $w_-(k )$ is the south pole, where $\nu =c\omega $ for some
  $c<0$. Equivalently this is the global maximum point of $u_0$.
  \item $\Gamma _+$ is the open boundary segment connecting
    $w_-(k )$ to $w_+(k )$ in the positive direction.
\item $w_+(k )$ is the north pole, where $\nu =c\omega $ for some
  $c>0$. Equivalently this is the global minimum point of $u_0$.
  \item $\Gamma _-$ is the open boundary segment connecting
    $w_+(k )$ to $w_-(k )$ in the positive direction.
  \end{itemize}
  We think here of $-u_{0}$ as the latitude, maximal at the north 
  pole and minimal at the south pole.
  
  Notice that
  \begin{equation}\label{iche.20}
\Gamma _\pm =\{ \gamma (t)\in \partial \Omega ;\, \mp \partial
_tu_0(\gamma (t))>0 \},
\end{equation}
where $\Gamma=\Gamma_{+}\cup\Gamma_{-}$. 

\medskip
\par On $\partial \Omega $ we have (\ref{iche.13}):
$$
e^{\overline{kw}-kw}=e^{-iu_0(w)},\ u_0(w)=|k|\Re
(w\overline{\omega  }),\ w=\gamma (t).
$$
The formula (\ref{ic.1}) reads
\begin{equation}\label{icbtf.1}
f(z,k)=\frac{1}{2i\overline{k}}\int_{\partial \Omega
}\frac{1}{z-w}e^{-iu(w,k )}dw+
  (\pi/\overline{k})e^{-i|k|\Re (z\overline{\omega })}1_\Omega(z) ,
\end{equation}
where $u(\cdot ,k)$ is the holomorphic extension of
$u_0(w)=u_0(w,k)$ to a neighborhood of $\partial \Omega $.
\subsection{Contour deformation}
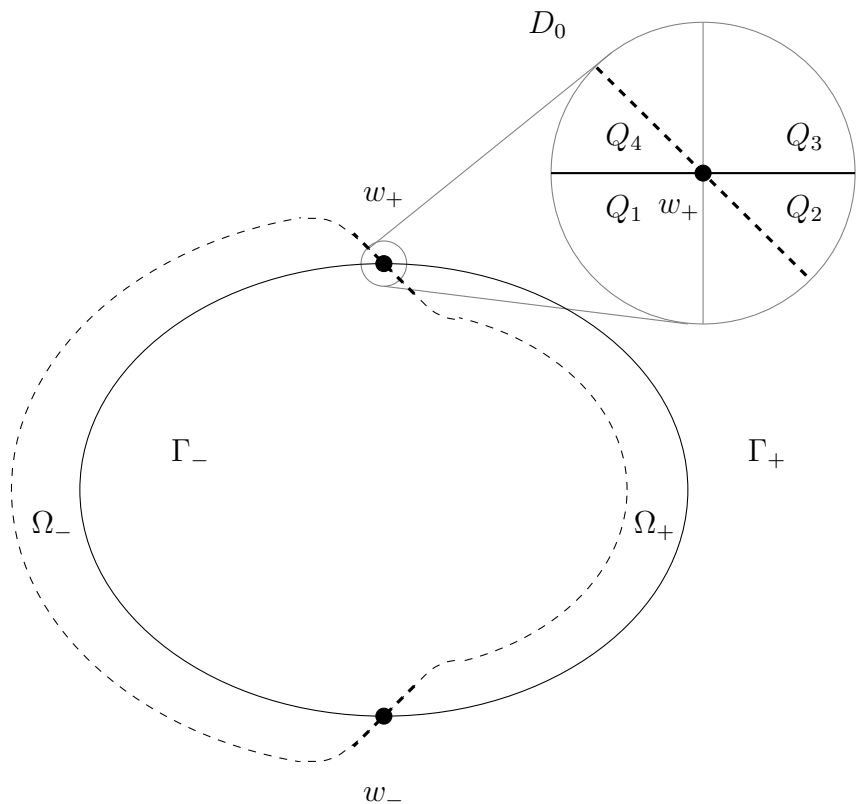
\begin{figure}[htb!]
%\begin{tikzpicture}
%\draw (0,0)  circle [x radius=4 cm, y radius = 3cm];
%
%\draw[very thick, black, dashed] (0.4,2.6)--(-.4,3.4);
%\draw[very thick, black, dashed] (0.4,-2.6)--(-.4,-3.4);
%\node (node) {};
%\node [above left = 0 cm and 2 cm of node](label){$\gamma_-$};
%\node [above right = 0 cm and 4.5 cm of node](label){$\gamma_+$};
%\begin{scope}
%    \clip (-1.1, -3.6) rectangle (-5, 8);
%    \draw(0,0)[dashed] circle   [x radius=4.9 cm, y radius = 3.7cm];
%\end{scope}
%
%\begin{scope}
%    \clip (1.,-3) rectangle (6,6);
%    \draw(-.0,0)[dashed] circle   [x radius=3.2 cm, y radius = 2.4cm];
%\end{scope}
%
%\filldraw 
%(0,+ 3) circle (3pt) node[align=center, above] {$w_+$\\ \vspace{.2cm}}
%(0, -3) circle (3pt) node[align=center, below] {~\\ ~\\ $w_-$} ;
%\draw[dashed,rounded corners=5] (-1,-3.6)--(-.7,-3.6)--(-.3, -3.3)--(0,-3)--(.3,-2.7)--(.7,-2.3) --(1,-2.25);
%\draw[dashed,rounded corners=5] (-1,3.6)--(-.7,3.6)--(-.3, 3.3)--(0,3)--(.3, 2.7)--(.7,2.3) --(1,2.25);
%\end{tikzpicture}
\begin{tikzpicture}

\draw (0,3) [gray] circle[x radius=0.3 cm, y radius = 0.3cm];
\draw [gray] (-.2,3.2) -- (3.0,5.8);
\draw [gray] (0,2.7) -- (4,2.2);
\draw  (4.2, 4.2) [gray]  circle [radius = 2cm];
\draw  [gray](4.2,2.2)--(4.2,6.2);
\draw [very thick, dashed] (2.8,5.6) -- (5.6,2.8);
\draw[thick, black ] (2.2,4.2)--(6.2,4.2);
\node at (4.2, 4.2)   (zoomnode) {};
\node [below left = 0.0 cm and 0.5 cm of zoomnode](label){$Q_1$};
\node [below right = 0.0 cm and 0.8 cm of zoomnode](label){$Q_2$};
%\node [zoomnode](label){$\w_+$};
\node [above left = 0.0 cm and 0.5 cm of zoomnode](label){$Q_4$};
\node [above right = 0.0 cm and 0.8 cm of zoomnode](label){$Q_3$};
\node [above left = 1.5 cm and 1.5 cm of zoomnode](label){$D_0$};

\filldraw (4.2, 4.2) circle (3pt) node[align=center, below left = -0.6 cm and -0.1 cm] {~\\ ~\\ $w_+$} ;

\draw (0,0)  circle [x radius=4 cm, y radius = 3cm];

\draw[very thick, black, dashed] (0.4,2.6)--(-.4,3.4);
\draw[very thick, black, dashed] (0.4,-2.6)--(-.4,-3.4);

\node (node) {};
\node [above left = 0 cm and 2 cm of node](label){$\Gamma_-$};
\node [above right = 0 cm and 4.5 cm of node](label){$\Gamma_+$};
\node [below right = 0 cm and 3 cm of node](label){$\Omega_+$};
\node [below left = 0 cm and 3.8 cm of node](label){$\Omega_-$};

\begin{scope}
    \clip (-1.1, -3.6) rectangle (-5, 8);
    \draw(0,0)[dashed] circle   [x radius=4.9 cm, y radius = 3.7cm];
\end{scope}

\begin{scope}
    \clip (1.,-3) rectangle (6,6);
    \draw(-.0,0)[dashed] circle   [x radius=3.2 cm, y radius = 2.4cm];
\end{scope}

\filldraw 
(0,+ 3) circle (3pt) node[align=center, above] {$w_+$\\ \vspace{.2cm}}

(0, -3) circle (3pt) node[align=center, below] {~\\ ~\\ $w_-$} ;
\draw[dashed,rounded corners=5] (-1,-3.6)--(-.7,-3.6)--(-.3, -3.3)--(0,-3)--(.3,-2.7)--(.7,-2.3) --(1,-2.25);
\draw[dashed,rounded corners=5] (-1,3.6)--(-.7,3.6)--(-.3, 3.3)--(0,3)--(.3, 2.7)--(.7,2.3) --(1,2.25);
\end{tikzpicture}

 \caption{Real analytic strictly convex boundary $\gamma$ (solid) of some domain 
 $\Omega$  and the deformed contour $\Gamma$ (dashed) for this 
 example. A close-up of the region near the north pole is shown. }
 \label{contour}
\end{figure}

From
(\ref{iche.20}) we see that the modulus of the exponential factor in
the integral decreases in the
following two situations:
\begin{itemize}
\item We start from a point in $\Gamma _+$ and move a short distance
  into $\Omega $.
  \item We start from a point in $\Gamma _-$ and move a short distance
    outward to $\mathbb{C}\setminus \overline{\Omega }$. 
  \end{itemize}
  Correspondingly, we seek to deform $\Gamma _+$ to a new contour
  slightly inside $\Omega $ and $\Gamma _-$ to a new contour slightly
  outside $\overline{\Omega }$. Naturally if such a deformation in the
  $w$ plane
  crosses the singularity at $w=z$ we will pick up a residue term. 
  We need to give a more precise description of the deformation near
  the poles and concentrate on the case of $w_+(k)$ for simplicity. 
  We can find a bi-holomorphic map
$K:\mathrm{neigh\,}(0,\mathbb{C})\to \mathrm{neigh\,}(w_+(k 
),\mathbb{C})$ mapping $\mathrm{neigh\,}(0,\mathbb{R})$ with the positive
orientation onto $\mathrm{neigh\,}(w_+(k ),\partial \Omega )$
also with the positive (anti-clockwise) orientation, such
that
\begin{equation}\label{icbtf.7}
u(K(\mu ) ,k )=|k |\left( \frac{u(w_+(k ),k
    )}{|k | }+\frac{\mu ^2}{2} \right). \end{equation} Let $D_0$ be a small
closed disc centered at $w_+(k )$ and let $Q_j$ be the image in
$D_0$ of the closed $j$:th quadrant under $K$. Thus $Q_j$ are
``distorted quadrants'' in $D_0$ with $Q_1$ and $Q_2$ contained in
$\overline{\Omega }$ while $Q_3$ and $Q_4$ are contained in
$D_0\setminus \Omega $. $Q_2\cap Q_3=\overline{\Gamma }_+\cap D_0$
while $Q_4\cap Q_1=\overline{\Gamma } _-\cap D_0$.
We show a   schematic view of the contour deformation for an example 
in Fig.~\ref{contour}.

\par Here $-\Im u(w,k )\ge 0$ in $Q_2\cup Q_4$ with strict
inequality in the interior, while $-\Im u(w,k ) \le 0$ in
$Q_3\cup Q_1$ with strict inequality in the interior. In other words,
$e^{-iu(w,k )}$ is bounded in $Q_2\cup Q_4$ and exponentially
decaying in the interior. It is exponentially large in the interior of
$Q_3\cup Q_1$. Naturally we have a similar description near $w_-(k)$.

\par We deform
$\partial \Omega $ inwards from $\Gamma_+$ and outward from $\Gamma_-$ and
so that the deformed curve $\Gamma$ follows the curve given by
$w=K(\mu )$ where $\mathrm{arg\,}\mu =3\pi /4$ in $D_0\cap Q_2$
and $\mathrm{arg\,}\mu =-\pi /4$ in $D_0\cap Q_4$. (In the
$\mu $-variable, $\Gamma $ here coincides with the oriented line
$e^{-i\pi /4}\mathbb{R}$.)  Thus along this
part of $\Gamma $, we have
$$
-(u(K(\mu ),k )-u(w_+(k ),k ))=i|k ||\mu |^2/2.
$$
Notice that we have just followed the rule of steepest descent. In
$\Gamma _N\setminus \overline{D}_0$, we have
$-\Im (u-u(w_+(k )))\asymp |k |$. Here $\Gamma_{N}$ denotes the 
`northern part of $\Gamma$', defined below after (\ref{icbtf.5}). We do the analogous construction
near $w_-(k)$. Then along $\Gamma$ we have
\begin{equation}\label{icbtf.8}
|e^{-iu(w,k )}|\le e^{-|k |\mathrm{dist\,}(w,\{w_+(k ),w_-(k)\} )^2/C}.
\end{equation}

Let
$\upsilon $ be a smooth vector field, defined near $\partial \Omega $,
transversal to $\partial \Omega $ and pointing outward, then we can
assume that $\Gamma =\{ (\exp t\upsilon (w))_{t=\tau (w)};\, w\in
\partial \Omega  \}$, where $\tau $ is a suitable smooth function on
$\partial \Omega $, $>0$ on $\Gamma _-$ and $<0$ on $\Gamma _-$. Let
$\Omega _+\subset \overline{\Omega } $ and $\Omega _-\subset \mathbb{C}\setminus \Omega $, be the points swept over by the deformation of
$\Gamma _+$ and $\Gamma _-$ respectively,
$$\Omega _+=\{ \exp t\upsilon (w);\, w\in \Gamma _+,\ \tau (w)\le t
\le 0  \},$$
$$\Omega _-=\{ \exp t\upsilon (w);\, w\in \Gamma _-,\ 0\le t
\le \tau (w)  \},$$
and notice that
$$
\Gamma =\partial ((\Omega \setminus \Omega _+)\cup \Omega _-).
$$

Assume for simplicity that $z\not\in \partial \Omega_{+} \cup \partial
\Omega_{-} $ and put
\begin{equation}\label{FGam}
F(z)=F_\Gamma (z)=\int_\Gamma \frac{1}{z-w}e^{-iu(w,k )}dw.
\end{equation}
When deforming the contour $\partial \Omega $ in (\ref{icbtf.1}) into
$\Gamma $ we have to add a residue term if the pole $w=z$ belongs to
$\Omega _+$ or $\Omega _-$. 
By the residue theorem we get (cf.\ (\ref{ic.1}), (\ref{icbtf.1}))
\begin{equation}\label{wn.3}
  \begin{split}
&f(z,k )=\frac{1}{2 i \overline{k}}\int_{\partial \Omega
}\frac{1}{z-w}e^{-iu(w,k )}dw+(\pi /\overline{k})e^{-i|k|\Re
  (z\overline{\omega })}1_\Omega (z)\\
&=\frac{1}{2i\overline{k}}F(z)+(\pi
/\overline{k})\left(e^{-iu(z,k )}(1_{\Omega _-}(z)-1_{\Omega
    _+}(z))+e^{-i|k|\Re (z\overline{\omega })}1_\Omega (z) \right).
  \end{split}
\end{equation}
\begin{remark}\label{ic1}
Away from a neighborhood of $\{w_+(k ),w_-(k ) \}$, we have
flexibility in the choice of $\Gamma $, which determines whether we
should count a residue term or not. This is only an apparent
difficulty because the residue terms become exponentially small when
$z$ approaches $\Gamma $ when $|z-w_\pm(k )|^{2}|k |$ are $\gg
1$. Indeed, if $z\in \mathrm{neigh\,}(w_\pm (k )$ and $|z-w_\pm(k
)|^2|k |\ge 1/\mathcal{ O}(1)$, then we have
\begin{equation}\label{lp.21}
|e^{-iu(z,k )}|\le \mathcal{ O}(1)e^{-\mathrm{dist\,}(z,\partial
  \Omega )|z-w_\pm(k )||k |/\mathcal{ O}(1)},\hbox{ when }z\in
\Omega _+\cup \Omega _-.
\end{equation}
\end{remark}

\subsection{Asymptotics}
It remains to study the asymptotics of $F=F_\Gamma $ in (\ref{wn.3}).
\par We consider different cases depending on $z$:\\

%\medskip
\paragraph{\it Case 1. $\mathrm{dist\,}(z,\{w_+(k),w_-(k)\})\ge
1/\mathcal{ O}(1)$.} After a deformation of $\Gamma $ which does not change
the general properties above and which does not cross the pole at
$w=z$, we may assume that $(z-w)^{-1}=\mathcal{ O}(1)$ along $\Gamma
$. $F_\Gamma $ can then be expanded by the method of stationary phase
-- steepest descent. The asymptotics is determined by the behaviour of
the integrand near the critical points $w_\pm (k)$ and coincides with
the one we would get directly from the corresponding integral over
$\mathrm{neigh\,}(\{w_+(k),w_-(k) \},\partial \Omega )$. $F$ in
(\ref{FGam}) has the asymptotic behaviour:
\begin{multline}\label{icbtf.2}
\sqrt{2\pi }
\left(\frac{1}{z-w_+(k )}e^{-iu(w_+(k ),k )-i\pi /4}
\frac{\dot{\gamma
}(t_+(k ))}{|\partial _t^2u(t_+(k ) )|^{1/2}}\right. \\
\left.
+\frac{1}{z-w_-(k )}e^{-iu(w_-(k ),k )+i\pi /4}
\frac{\dot{\gamma
}(t_-(k ))}{|\partial _t^2u(t_-(k ) )|^{1/2}}
\right)\\
+\mathcal{O}(\langle z\rangle^{-1} k^{-3/2}),
\end{multline}
where we write $w_\pm(k )= \gamma (t_\pm(k ))$ and identify $u(w)$
with $u(\gamma (t))=$ whenever convenient. Since
$u(\gamma (t))=u_0(\gamma (t))$ for real $t$, we get
$\partial _t^2u(\gamma (t_\pm(k )))$ from (\ref{iche.16}). The leading
term is $\mathcal{ O}(k^{-1/2}\langle z\rangle^{-1})$. Notice that the
choice of $\Gamma $ is coherent with the general principle of steepest
descent leading to a new contour $\Gamma $, passing through the saddle
points $w_+$ and $w_-$ of $-\Im u(\cdot ,k )$, so that the integrand
restricted to $\Gamma $ is exponentially small away from those points.\\

\paragraph{\it Case 2. $z$ is close to
  $w_+(k )$ or $w_-(k )$. } 
To fix the ideas, we assume that $z$ is close to the north pole,
\begin{equation}\label{icbtf.4}
  |z-w_+(k )|\ll 1.
\end{equation}

\par
Let $a_\pm\in \Gamma _\pm$ be independent of $z$. We assume that
\begin{equation}\label{icbtf.5}
|z-w_+(k )|\ll |a_\pm -w_+(k )|.
\end{equation}
We decompose $\Gamma $ into the union of two segments;
$\Gamma =\Gamma _N\cup \Gamma _S$, where $\Gamma _N$ is the part of
$\Gamma $ that runs from a point $\widetilde{a}_+\in \Gamma $ near
$a_+$, through the north pole $w _+(k )$, to a point
$\widetilde{a}_-\in \Gamma $ near $a_-$. $\Gamma _S$ is the
remaining part of $\Gamma $, which runs from $\widetilde{a}_-$ through
the south pole $w _-(k )$ to $\widetilde{a}_+$. Define $F_N(z)$,
$F_S(z)$ as in (\ref{FGam}) but with $\Gamma $ replaced by $\Gamma _N$
and $\Gamma _S$ respectively. Clearly,
\begin{equation}\label{icbtf.6}
F(z)=F_N(z)+F_S(z).
\end{equation}

\par For $z$ close to $w_+(k )$ the integral $F_S(z)$
is analyzed by stationary phase -- steepest descent and gives a contribution as in the
second term in (\ref{icbtf.2}).\\

\paragraph{\it Asymptotic expansions when $|z-w_+(k )|^2|k |\gg
  1$.} Recall that we study the case when $|z-w_+(k )|\ll 1$ and
hence when $z\in K(D_0)$. We now add the assumption that $|z-w_+(k
)|\gg |k |^{-1/2}$, i.e.\ 
\begin{equation}\label{icbtf.13}|z-w_+(k )|^2|k |\gg
  1 .\end{equation}

\par If $\mathrm{dist\,}(z,\Gamma _N)\ll |z-w_+(k )|$, we make a
slight deformation of $\Gamma _N$ inside $K(D_0)$ to achieve that
\begin{equation}\label{icbtf'.1}
\mathrm{dist\,}(z,\Gamma _N)\ge \frac{|z-w_+(k )|}{\mathcal{ O}(1)}.
\end{equation}
In fact, we may assume that we are inside a small disc $D_1\subset D_0$
for the $\mu $-variables and that $\mathrm{arg\,}\mu \in \{-\pi
/4,3\pi /4 \}$ on the corresponding part of $\Gamma _N$. In order
to get (\ref{icbtf'.1}) it suffices to rotate this part slightly, so
that we get instead
\begin{equation}\label{ictbf'.2}
\mathrm{arg\,}\mu \in \{c-\pi /4,c+3\pi /4 \},
\end{equation}
where $c\in \mathrm{neigh\,} (0,\mathbb{R})$
is a suitable small constant. (Here it is understood that we deform in
the right direction, avoiding to cross the pole at $w=z$.)

\par Possibly after this additional deformation, we can achieve that
\begin{equation}\label{icbtf'.3}
|w-z|\ge \frac{1}{\mathcal{ O}(1)}\left(|z-w_+(k )| +|k |^{-1/2}
\right),\ w\in \Gamma _N.
\end{equation}
We may still assume that (\ref{icbtf.8}) holds along $\Gamma _N$. It
follows that
$$
\left| e^{-iu(w,k )} \right| = \mathcal{ O}(1)e^{-|k |/C}\hbox{
  on }\Gamma _N\setminus K(D_1),
$$
so after committing a corresponding exponentially small error, we can
replace $\Gamma _N$ with $\Gamma _N\cap K(D_1)$ in the definition of
$F_N(z)$ in (\ref{icbtf.6}).

\par From (\ref{icbtf'.3}) we have
\begin{equation}\label{icbtf'.4}
\frac{1}{z-w}=\frac{\mathcal{ O}(1)}{|z-w_+(k )|+|k |^{-1/2}},\
w\in \Gamma _N
\end{equation}
and this estimate persists for $w\in \mathrm{neigh\,}(\Gamma _N)$ with
$\mathrm{dist\,}(w,\Gamma _N)\le |z-w_+(k )|/\mathcal{ O}(1)$.

Let $\alpha =|z-w_+(k )|$ and make the change of variables
$$
w-w_+(k )=\alpha \widetilde{w},\ \widetilde{w}\in \alpha
^{-1}((\Gamma _N\cap K(D_1))-\{w_+(k ) \})=:\widetilde{\Gamma
}_{N,\alpha }.
$$
Then
$$
u(w,k )-u(w_+(k ),k )=\alpha ^2|k
|\widetilde{u}(\widetilde{w}),\ \partial _w^2u=|k |\partial _{\widetilde{w}}^2\widetilde{u},
$$
where
\begin{equation}\label{icbtf'.5}
  \begin{cases}
    -\Im \widetilde{u}(\widetilde{w})\asymp |\widetilde{w}|^2 , \hbox{
      along } \widetilde{\Gamma
    }_{N,\alpha },\\
    \widetilde{u}(\widetilde{w})=\mathcal{O}(\widetilde{w}^2), \hbox{
      in }\alpha ^{-1}(K(D_1))-\{ w_+(k ) \} .
  \end{cases}
  \end{equation}
  Up to an error $\mathcal{ O}(e^{-|k |/\mathcal{ O}(1)})$, we get
  \begin{equation}\label{icbtf'.6}
    F_N(z)=e^{-iu(w_+(k ),k )}
    \int_{\widetilde{\Gamma }_{N,\alpha }}\frac{1}{\frac{z-w_+(k
        )}{|z-w_+(k )|}-\widetilde{w}}\, e^{-i\alpha ^2|k |\widetilde{u}(\widetilde{w})}d\widetilde{w}.
  \end{equation}
  Here $\alpha ^2|k |\gg 1$ is our new large parameter and from
  (\ref{icbtf'.4}) we see that
  $$
\frac{1}{\frac{z-w_+(k
        )}{|z-w_+(k )|}-\widetilde{w}}=\mathcal{ O}(1)\hbox{ on
    }\widetilde{\Gamma }_{N,\alpha }
    $$
    and even on a larger set $\{\widetilde{w};\,
    \mathrm{dist\,}(\widetilde{w},\widetilde{\Gamma }_{N,\alpha })<
    (1+|\widetilde{w}|)/\mathcal{ O}(1) \}$.

    \par It is then clear that we can apply the method of stationary
    phase (steepest descent) to the integral in (\ref{icbtf'.6}) which
    has a complete asymptotic expansion in powers of $\alpha ^2|k
    |$. Since the choice of $\alpha $ can be modified by
    multiplication with any positive constant of order $\asymp 1$, we
    know in advance that each term in the asymptotic series is
    actually independent of $\alpha $ and in particular the leading
    term in the asymptotic series is independent of $\alpha $ and
    therefore has to coincide with the contribution from $w_+(k
    )$ to the expression in (\ref{icbtf.2}). Let us nevertheless review this in more
    detail.

    \par Recall that $\gamma (t)$ parametrizes $\partial \Omega $,
    that $|\dot{\gamma }(t)|=1$ and that deformations of $\partial
    \Omega $ are then naturally parametrized by $w=\gamma (t)$ where
    we let $t$ follow a deformation of $\mathbb{R}$ in the complex
    $t$-plane. After a translation in $t$, we may assume that
    $w_+(k )=\gamma (0)$.
   In $(\ref{icbtf'.6})$ we can use $\widetilde{w}=(\gamma
    (t)-w_+(k ))/\alpha =:\widetilde{\gamma }(\widetilde{t)}$,
    where $t=\alpha \widetilde{t}$, so that $\dot{\widetilde{\gamma
      }}(\widetilde{t}) = \dot{\gamma }(t)$. We get
    \begin{multline*}
    F_N(z)+\mathcal{ O}((\alpha ^2|k |)^{-3/2})=\\
    e^{-iu(w_+(k ),k )-i\frac{\pi }{4}}\frac{\sqrt{2\pi
      }}{\alpha |k
      |^{1/2}|\partial _{\widetilde{w}}^2\widetilde{u}(0)|^{1/2}}\frac{|z-w_+(k
      )|}{z-w_+(k )}\dot{\gamma }(0).
  \end{multline*}
Recalling that
  $\alpha =|z-w_+(k )|$ and $|k |\partial
  _{\widetilde{w}}^2\widetilde{u}=\partial _w^2u$, we get
  \begin{equation}\label{icbtf'.7}\begin{split}
      F_N(z)=&e^{-iu(w_+(k ),k )-i\frac{\pi }{4}}
\frac{\sqrt{2\pi
  }}{|u''(w_+(k ),k )|^{1/2}}\frac{1}{z-w_+(k
  )}\dot{\gamma }(0)
      \\ &\hskip 1cm  +\mathcal{
  O}((|z-w_+(k )|^2|k |)^{-3/2}),
\end{split}
  \end{equation}
  where the modulus of the first term in the right hand side is
  $$
\asymp (|z-w_+(k )|^2|k |)^{-1/2}.
  $$
~\\

 \paragraph{\it Limiting profile} We allow deformations of $\Gamma _N$
 inside a set
 \begin{equation}\label{lp.1}
V_N=\{w\in \mathrm{neigh\,}(w_+(k ));\, \mathrm{dist\,}(w,\Gamma
_N^0)\le
\mathcal{ O}(1)|k |^{-1/2}+|w-w_+(k )|/\mathcal{ O}(1)
\},
 \end{equation}
where $\Gamma _N^0$ is the contour, given in the $\mu $-variables
by
$$
\{\mu \in \mathrm{neigh\,}(0);\, \mathrm{arg\,}\mu \in \{ 3\pi
/4,-\pi /4 \} \}.
$$
In $V_N$ we have
\begin{equation}\label{lp.2}
|e^{-iu(w,k )}|\le \mathcal{ O}(1)e^{-|k ||w-w_+(k
 ) |^2/\mathcal{ O}(1)}.
\end{equation}
Let $u_2(w-w_+,k )$ be the quadratic Taylor polynomial of
$u(w,k )$ at $w_+=w_+(k )$. Possibly after shrinking the fixed
neighborhood of $w_+$ where we work, (\ref{lp.2})  remains valid
uniformly if we replace $u$ by $u^t=tu+(1-t)u_2(w-w_+,k )$ for
$0\le t\le 1$. By differentiation,
\[
  \begin{split}
    \partial _te^{-iu^t(w,k )}&=i(u_2-u)e^{-iu^t(w,k )}\\
    &=\mathcal{ O}(|k ||w-w_+|^3)e^{-|w-w_+|^2| k |/\mathcal{ O}(1)}\\
    &=\mathcal{ O}(|k |^{-1/2})e^{-|w-w_+|^2|k |/(2\mathcal{ O}(1))},
      \end{split}
\]
and by integration,
\begin{equation}\label{lp.3}
  e^{-iu(w,k )}-e^{iu_2(w-w_+,k )}=
  \mathcal{ O}(|k |^{-1/2})e^{-|w-w_+|^2|k |/(2\mathcal{ O}(1))}.
\end{equation}

\par We can choose the deformation $\Gamma _N$ close to $\Gamma _N^0$
so that
\begin{equation}\label{lp.4}
\frac{1}{z-w}=\frac{\mathcal{ O}(1)}{|z-w_+(k )|+|k |^{-1/2}},\
w\in \Gamma _N.
\end{equation}
Combining (\ref{lp.2}), (\ref{lp.4}) in (\ref{FGam}) with $\Gamma $
replaced by $\Gamma _N$, we get
\begin{equation}\label{lp.5}
F_N(z)=\frac{\mathcal{ O}(1) }{|k | ^{1/2}(|z-w_+(k )|+|k
  |^{-1/2})}=
\frac{\mathcal{ O}(1)}{|z-w_+(k )||k |^{1/2}+1}.
\end{equation}
Using also (\ref{lp.3}), we get
\begin{equation}\label{lp.6}
F_N(z)-F_{N,2}(z)=\frac{\mathcal{ O}(1)}{|k |^{1/2}(|z-w_+(k
  )||k |^{1/2}+1)},
\end{equation}
where
\begin{equation}\label{lp.7}
F_{N,2}(z)=\int_{\Gamma _N}\frac{1}{z-w}e^{-iu_2(w-w_+(k ),k
  )}dw.
\end{equation}
Up to an exponentially small error ($\mathcal{ O}(\exp (-|k |/\mathcal{
  O}(1))$), we may here assume that $\Gamma _N$ is a straight line, 
  whose intersection with $\mbox{neigh}(w_{+}(k))$ is 
contained in the set $V_N$ in (\ref{lp.1}).

Recalling that $u$ is the holomorphic extension from $\partial \Omega
$ of $u_0$, we know from (\ref{iche.16}) that
\begin{equation}\label{lp.8}
(\dot{\gamma }(0)\partial _w )^2u_2(w_+(k ),k
)=|k|\langle \ddot{\gamma }(0),\omega \rangle_{\mathbb{R}^2}=:a>0,
\end{equation}
where for simplicity we assume that $\gamma (0)=w_+(k )$. Here
we also recall that $\omega =2\nu (w_+(k ))$, where
$\nu =\ddot{\gamma }(0)/|\ddot{\gamma }(0)|=i\dot{\gamma }(0)$ is the
interior unit normal at $w_+(k )$ (see the dscussion after 
(\ref{iche.19})). Equivalently,
\begin{equation}\label{lp.9}
u_2(w-w_+(k ),k )=\frac{a}{2}\left(\frac{w-w_+(k
    )}{\dot{\gamma }(0)} \right)^2.
\end{equation}

\par Comparing with (\ref{icbtf.7}), we see that
$$
\frac{w-w_+(k )}{\dot{\gamma }(0)}=c\mu +\mathcal{ O}(\mu ^2)
$$ for some $c=c_k >0$ and we can therefore identify the
quadrants $Q_j$, defined in the $\mu $-plane with those in the
$(w-w_+(k ))/\dot{\gamma }(0)$ plane. Recall that $\Gamma _N$ is
now a straight oriented line close to $\Gamma _{N,0}=e^{-i\pi 
/4}\mathbb{R}$.

Put
\begin{equation}\label{lp.10}
\widetilde{w}=\frac{e^{i\pi /4}\sqrt{a}}{\dot{\gamma
  }(0)}(w-w_+(k )), \hbox{i.e. } w=w_+(k)+\frac{\dot{\gamma
  }(0)}{\sqrt{a}e^{i\pi /4}}\widetilde{w}, 
\end{equation}
so that (cf.\ (\ref{lp.9}),
\begin{equation}\label{lp.11}
iu_2(w-w_+(k ),k )=\frac{\widetilde{w}^2}{2},
\end{equation}
and $\Gamma _N$ becomes a straight line $\widetilde{\Gamma }$ close to
the positively oriented real axis. Define $\widetilde{z}$ similarly by 
\begin{equation}\label{lp.12}
\widetilde{z}=\frac{e^{i\pi /4}\sqrt{a}}{\dot{\gamma
  }(0)}(z-w_+(k )), \hbox{i.e. } z=w_+(k)+\frac{\dot{\gamma
  }(0)}{\sqrt{a}e^{i\pi /4}}\widetilde{z}. 
\end{equation}
Then we get
\begin{equation}\label{lp.13}
F_{N,2}(z)=G(\widetilde{z}),
\end{equation}
where
\begin{equation}\label{lp.14}
G(\widetilde{z})=\int_{\widetilde{\Gamma }}\frac{1}{\widetilde{z}-\widetilde{w}}e^{-\widetilde{w}^2/2}d\widetilde{w}.
\end{equation}

\par Clearly, $G(\widetilde{z})$ does not change if we deform the
contour into a new (straight line) contour with the same properties
(contained in $D(0,\mathcal{ O}(1))+\exp (i]-\pi /4+1/\mathcal{ O}(1), \pi
/4-1/\mathcal{ O}(1)[)\mathbb{R}$, {\it provided that we do not cross
  the pole at $\widetilde{w}=\widetilde{z}$}. Define
$G_r(\widetilde{z})$ as in (\ref{lp.14}) with ``$\widetilde{\Gamma }$
{\it passing below} $\widetilde{z}$'' or equivalently with
``$\widetilde{z}$ always to the left'' when traveling along
$\widetilde{\Gamma }$ in the positive direction. (Thus
``$\widetilde{\Gamma }$ passes to the right of $\widetilde{z}$''.)
Define $G_\ell$ when $\widetilde{z}$ remains to the right when
following $\widetilde{\Gamma }$ with the natural orientation. By the
residue theorem,
\begin{equation}\label{lp.15}
  G_r(\widetilde{z})-G_\ell (\widetilde{z})=-2\pi i e^{-\widetilde{z}^2/2}.
\end{equation}
We define $F_{N,2}^r$, $F_{N,2}^\ell$ similarly. Then
\begin{equation}\label{lp.16}
F_{N,2}^\sigma(z) =G_\sigma (\widetilde{z}),\ \ \sigma =r,\ell .
\end{equation}
If $\Im \widetilde{z}\ge 0$, we can choose a suitable contour
$\widetilde{\Gamma }$, ``passing below'' $\widetilde{z}$, to see that
\begin{equation}\label{lp.17}
G_r(\widetilde{z})=\frac{\mathcal{ O}(1)}{\langle \widetilde{z}\rangle}.
\end{equation}
Similarly, if $\Im \widetilde{z}\le 0$, we have
\begin{equation}\label{lp.18}
G_\ell(\widetilde{z})=\frac{\mathcal{ O}(1)}{\langle \widetilde{z}\rangle}.
\end{equation}
Using (\ref{lp.15}), we then get
\begin{equation}\label{lp.19}
|G_r(\widetilde{z})|,\, |G_\ell (\widetilde{z})|\le \frac{\mathcal{
    O}(1)}{\langle z\rangle}+2\pi |e^{-\widetilde{z}^2/2}|,
\end{equation}
uniformly for $\widetilde{z}\in \mathbb{C}$. \\

%\medskip
\paragraph{\it Summary} 
We assume for the simplicity of the presentation that
$z\not\in \partial \Omega _+\cup \partial \Omega _-$. 

\par If $|z-w_+(k )|^2|k |\le \mathcal{ O}(1)$, the
corresponding contribution to (\ref{icbtf.2}) is no longer pertinent
and we have to replace it by  
\begin{equation}\label{lp.20}
\frac{1}{2i\overline{k}}F_{N,2}^\sigma
(z)=\frac{1}{2i\overline{k}}G^\sigma (\widetilde{z})
\end{equation}
with $\widetilde{z}$ as in (\ref{lp.12}), where we choose $\sigma =r$
when $z$ is inside $\Gamma $ and $\sigma =\ell$ when $z$ is
outside. The same rule  applies for adding a residue term when $z\in \Omega
_+\cup \Omega _-$. Naturally, the same discussion applies near
$w_-(k )$ but we refrain from developing the details.

Putting everything together we get the following long theorem.
\begin{theo}\label{ic2}
Let $\Omega \Subset \mathbb{C}$ be strictly convex with real analytic
boundary and let $f(z,k)$, $z,k\in \mathbb{C}$, $|k|\ge 1$ be the
function appearing in (\ref{ic0.1}), (\ref{ic0.4}), (\ref{ic.0}). Let
$iu_0(w)=kw-\overline{kw}=i|k|\Re (z\overline{\omega })$ (cf.\
(\ref{omega})) and let $u$ be a holomorphic extension of ${{u_0}_\vert}_{\partial \Omega }$  to a neighborhood
of $\partial \Omega $.

\par Assuming
for simplicity that $z\not\in \partial \Omega $, we have
(\ref{icbtf.1}):
$$
f(z,k)=\frac{1}{2i\overline{k}}\int_{\partial \Omega
}\frac{1}{z-w}e^{-iu(w,k )}dw+
  (\pi/\overline{k})e^{-i|k|\Re (z\overline{\omega })}1_\Omega(z) .
  $$
  Let $w_+(k ),\, w_-(k )\in \partial \Omega $ be the points of minimum and
maximum of the (real valued) function $u_0$, also characterized by
$\omega \in \mathbb{R}_{\pm}\nu (w_\pm (k ))$, where $\nu (w)$
denotes the interior normal of $\partial \Omega $ at $w$. Choose a
parametrization $\mathbb{R}/|\partial \Omega |\mathbb{Z}\ni t\mapsto \gamma
(t)\in \partial \Omega $ with positive (anti-clock-wise) orientation
and $|\dot{\gamma }(t)|=1$, so that $\nu (\gamma (t))=\ddot{\gamma
}(t)/|\ddot{\gamma
}(t)|$.

\par Let $\Gamma _+\subset \partial \Omega $ be the open oriented
boundary segment from $w_-(k )$ to $w_+(k )$ and let $\Gamma _-$ be
the similar one from $w_+(k )$ to $w_-(k )$. Let $\Gamma $ be a
deformation of the oriented boundary $\partial \Omega $ as described
prior to Remark \ref{ic1}, see Fig.~\ref{contour}, and recall that $\Gamma $ is obtained by
pushing $\Gamma _+$ inward and $\Gamma _-$ outward, keeping
$w_\pm(k )$ fixed and so that $\Gamma $ coincides near $w_+(k )$ with
the (image of the) oriented line $e^{-i\pi /4}\mathbb{R}$ in the Morse
$\mu $-coordinates in (\ref{icbtf.7}) and similarly near $w_-(k
)$. Let $\Omega _\pm\Subset \mathbb{C}$ be the closed sets swept over,
when deforming $\partial \Omega $ to $\Gamma $. Let
$\hbox{Lead\,(\ref{icbtf.2})}=\hbox{Lead}_+\hbox{(\ref{icbtf.2})}+\hbox{Lead}_-\hbox{(\ref{icbtf.2})}$
denote the leading term in (\ref{icbtf.2}) with the natural
decomposition into contributions from $w_+(k)$ and $w_-(k)$.

\par Let 
\begin{equation}\label{lp.22}
d(z,k )=|k |\min (1, |z-w_+(k )|^2, |z-w_-(k
)|^2).
\end{equation}
We first consider the case when $d(z,k )\ge 1$. Then we have
\begin{equation}\label{lp.23}\begin{split}
&f(z,k )+\mathcal{ O}\left(\frac{d(z,k
    )^{-\frac{3}{2}}}{|k |\langle
    z\rangle} \right)=\hbox{Lead\,(\ref{icbtf.2})}/(2i\overline{k})\\ 
&
+
(\pi /\overline{k})\left( e^{-i\Re (z,\overline{\omega })}1_\Omega(z) 
-e^{-iu(z,k )}1_{\Omega _+}(z)
+e^{-iu(z,k )}1_{\Omega _-}(z)\right).
\end{split}
\end{equation}

\par To cover the remaining case, it suffices to consider the cases
when $|z-w_+(k )|\ll 1$ and $|z-w_-(k )|\ll 1$. Both cases
are similar and we formulate the result only when $|z-w_+(k )|\ll
1$:
\begin{multline}\label{lp.24}
  f(z,k )+\frac{\mathcal{ O}(1)}{|k|^{3/2}(1+d^{1/2})}=
  (F_{N,2}^\sigma (z)+\hbox{Lead}_-\hbox{(\ref{icbtf.2})}) /(2i\overline{k})+\\
(\pi /\overline{k})\left( e^{-i\Re (z,\overline{\omega })}1_\Omega(z)
-e^{-iu(z,k )}1_{\Omega _+}(z)+
e^{-iu(z,k )}1_{\Omega _-}(z) \right) .
\end{multline}
Here, $\sigma =r$ when $z$ is inside $\Gamma $ and $\sigma =\ell $,
when $z$ is outside. $F_{N,2}^\sigma $ is introduced in
(\ref{lp.20}) and the discussion from (\ref{lp.12}) to (\ref{lp.19}). 
In particular 
$F_{N,2}^{\sigma}/\bar{k}=\mathcal{O}(1)/(|k|(1+d^{1/2}))$, 
cf.~(\ref{lp.5}), (\ref{lp.6}). 
\end{theo}

\subsection{Estimates of weighted $L^2$-norms}\label{wn}
Recall the definition of $F_N(z)$ in
(\ref{icbtf.6}). By (\ref{lp.5}) we have
\begin{equation}\label{wn.1}
F_N(z)=\frac{\mathcal{ O}(1)}{(|z-w_+(k )|^2|k |+1)^{1/2}},
\end{equation}
uniformly for $z\in \mathbb{C}$. Similarly,
\begin{equation}\label{wn.2h}
F_S(z)=\frac{\mathcal{ O}(1)}{(|z-w_-(k )|^2|k |+1)^{1/2}}.
\end{equation}

\par We shall estimate $\| f(\cdot ,k )\|_{\langle \cdot
  \rangle^\epsilon L^2}$ for $0<\epsilon \le 1$. We start with the
term
$$
\frac{1}{2i\overline{k}}F=\frac{1}{2i\overline{k}}F_N+\frac{1}{2i\overline{k}}F_S,
$$
appearing in (\ref{wn.3}).
Let $K\Subset\mathbb{C}$ be fixed and fix $r>0$ large enough so that
$$
K\subset D(w_+(k ),r)\cap D(w_-(k ),r).
$$
Then,
\begin{multline*}
\| F_N\|_{L^2(K)}^2\le \mathcal{ O}(1)\int_{D(0,r)}\frac{1}{|k
  ||z|^2+1}L(dz)\\
=\frac{\mathcal{ O}(1)}{|k |}\int_{D(0,r|k
  |^{1/2})}\frac{1}{|\widetilde{z}|^2+1}L(d\widetilde{z})=\frac{\mathcal{
    O}(1)\ln |k |}{|k |},
\end{multline*}
and we have the same estimate for $\| F_S\|_{L^2(K)}^2$ and hence
also for $\| F\|_{L^2(K)}^2$. Thus,
\begin{equation}\label{wn.4}
\left\Vert \frac{1}{2i\overline{k}} F\right \Vert_{L^2(K)}=\mathcal{
  O}(1)\frac{(\ln |k|)^{1/2}}{|k|^{3/2}}.
\end{equation}
For $z\in \mathbb{C}\setminus K$ we have uniformly,
$$
F=\frac{\mathcal{ O}(1)}{|k |^{1/2}|z|},
$$
assuming $K$ large enough so that $|z-w_{\pm}(k )|\asymp
|z|\asymp 1+|z|$ for $z\not\in K$. Then
$$
\|F\|^2_{\langle \cdot \rangle^\epsilon L^2(\mathbb{C}\setminus
  K)}=\frac{\mathcal{ O}(1)}{|k|}\int_{|z|\ge 1}\frac{1}{|z|^{2(1+\epsilon
  )}}L(dz)=\frac{\mathcal{ O}(1)}{|k|}.
$$
Hence,
\begin{equation}\label{wn.5}
\left\Vert \frac{1}{2i\overline{k}}F\right\Vert_{\langle \cdot 
\rangle^\epsilon L^2(\mathbb{C}\setminus
  K)}=\mathcal{ O}(1) \frac{(\ln |k|)^{1/2}}{|k|^{3/2}}
\end{equation}
and we have estimated the norm of the first term in the last
member in (\ref{wn.3}).

\par The estimate of the contribution from the last term in the
parenthesis in (\ref{wn.3}) is obvious:
\begin{equation}\label{wn.6}
\left\| (\pi /\overline{k})e^{i\Re (\cdot \overline{\omega
    })}1_\Omega   \right\|_{\langle \cdot \rangle^\epsilon
  L^2}=\frac{\mathcal{ O}(1)}{|k|}.
\end{equation}

\par We next consider the contribution from the other two terms in the
parenthesis in (\ref{wn.3}), so we look at $\pm (\pi
/\overline{k})e^{-iu(z,k )}$ in $\Omega _\mp$. Away from any fixed
neighborhood of $\{w_+(k ),\, w_-(k ) \}$, these terms are
$$
\frac{\mathcal{ O}(1)}{k}e^{-|k|\mathrm{dist\,}(z,\partial \Omega)/C }
$$
and the corresponding contributions to the squares of the $\langle \cdot
\rangle ^\epsilon L^2$ norms are
$$
\mathcal{ O}(k^{-2})\int _0^1 e^{-|k|s/C}ds=\mathcal{ O}(k^{-3}),
$$
so
\begin{equation}\label{wn.6.5}
  \left\|
1_{\Omega _\pm}e^{-iu}
  \right\|_{\langle \cdot \rangle^\epsilon L^2(\mathbb{C}\setminus
    \mathrm{neigh\,}(\{ w_+,w_-\}))}=\mathcal{ O}(|k|^{-3/2}).
\end{equation}
For the estimate of the contribution from a neighborhood of
$w_+(k )$ we use the $\mu $-variables from (\ref{icbtf.7}) and
get in $\Omega _+\cup \Omega _-$,
$$
(\pi /\overline{k})e^{-iu}1_{\Omega _\pm}=\mathcal{
  O}(k^{-1})e^{-|k||t||s|/C},\hbox{ when }\mu =t+is,\ |t|\le
1/\mathcal{ O}(1),\ |s|\le |t|.
$$
The contribution to the square of the $\langle \cdot \rangle^\epsilon L^2$-norm is
\begin{multline*}
\mathcal{ O}(|k|^{-2})\int_0^1\int_0^t e^{-ts|k|/C}dsdt=
\mathcal{ O}(|k|^{-2})\int_0^1\frac{1}{t|k|}\left(1 - e^{-t^2|k|/C}
\right)dt\\
=\frac{\mathcal{ O}(1) }{|k|^{2}}\left( \int_0^{|k|^{-1/2}}tdt
  +\int_{|k|^{-1/2}}^1\frac{1}{t|k|} dt \right)=
\frac{\mathcal{ O}(1)\ln |k|}{|k|^3}.
\end{multline*}
The same estimate holds for the contribution from a neighborhood of
$w_-(k )$ and we get
\begin{equation}\label{wn.7}
  \left\|
(\pi
/\overline{k})e^{-iu(\cdot ,k )}(1_{\Omega _-}-1_{\Omega
    _+})
  \right\|_{\langle \cdot \rangle^\epsilon L^2}=\frac{\mathcal{ O}(1)(\ln |k|)^{1/2}}{|k|^{3/2}}.
\end{equation}

\par Combining (\ref{wn.3}), (\ref{wn.5}), (\ref{wn.7}), we get
\begin{equation}\label{wn.8}
f(z,k )=(\pi /\overline{k}) e^{-i|k|\Re (z\overline{\omega
  })}1_\Omega(z)+ g(z,k ),\ \ \|g\|_{\langle\cdot  \rangle^\epsilon L^2}= \frac{\mathcal{ O}(1)(\ln |k|)^{1/2}}{|k|^{3/2}},
\end{equation}
where we recall that $kz-\overline{kz}=i|k|\Re (z\overline{\omega })$.

\par (\ref{ic0.4}) now gives,
\begin{equation}\label{wn.9}
F\widehat{\tau }_\omega \frac{h\overline{q}}{2}=\frac{1}{2k}e^{i|k|\Re
  (\cdot \overline{\omega })}1_\Omega +\mathcal{ O}(1)\frac{(\ln
  |k|)^{1/2}}{|k|^{3/2}}\hbox{ in }\langle \cdot \rangle^\epsilon L^2.
\end{equation}
By Proposition \ref{qcf1}, the estimates of Section 3
 are applicable with $s=3/2$ and we have seen after (\ref{sy.16}) that
 $A,B=\mathcal{ O}(1): \langle \cdot \rangle^\epsilon L^2\to \langle
 \cdot \rangle^\epsilon L^2$, so by (\ref{wn.9}),
 \begin{equation}\label{wn.10}
AF\widehat{\tau }_\omega \frac{h\overline{q}}{2}=\mathcal{ O}(h)\hbox{ in
}\langle \cdot \rangle^\epsilon L^2\ \ \ (h=1/|k|),
\end{equation}
and since $AB=\mathcal{ O}(h^{1/2}):\langle \cdot \rangle^\epsilon L^2\to
\langle \cdot \rangle^\epsilon L^2 $, we get from (\ref{sy.41}):
\begin{equation}\label{wn.11}
\phi _1^1=AF\widehat{\tau }_\omega \frac{h\overline{q}}{2}+\mathcal{
  O}(h^{3/2})=
A\left(\frac{1}{2k}e^{i|k|\Re (\cdot ,\overline{\omega })}1_\Omega
\right)+\mathcal{ O}(1)h^{3/2}(\ln (1/h))^{1/2}\hbox{ in }\langle \cdot
\rangle^\epsilon L^2.
\end{equation}
Similarly, since $BA=\mathcal{ O}(h^{1/2}): \langle \cdot \rangle^\epsilon L^2\to
\langle \cdot \rangle^\epsilon L^2 $, we get from (\ref{sy.42}):
\begin{equation}\label{wn.12}
\phi _2^1=\frac{1}{2k}e^{i|k|\Re (\cdot \overline{\omega })}1_\Omega +\mathcal{
  O}(1)h^{3/2}(\ln (1/h))^{1/2}\hbox{ in }\langle \cdot
\rangle^\epsilon L^2.
\end{equation}

\par Recall from (\ref{sy.16}) that
\begin{equation}\label{wn.13}
A=E\widehat{\tau }_{-\omega }\frac{hq}{2},\ \ B=\sigma F\widehat{\tau
}_\omega \frac{h\overline{q}}{2}\ \ \ (q=\overline{q}=1_\Omega ),
\end{equation}
where $q$ is viewed as a multiplication operator. Thus by
(\ref{wn.11}), (\ref{wn.13}):
\begin{multline*}
  \phi _1^1=E\widehat{\tau }_{-\omega }\frac{h}{2}1_\Omega
  F\widehat{\tau }_\omega \frac{h}{2}(1_\Omega) +\mathcal{ O}(h^{3/2})\\
  =E\widehat{\tau }_{-\omega }\frac{h}{2}1_\Omega 
\frac{1}{2k}e^{i|k|\Re (\cdot \overline{\omega })}(1_\Omega )+\mathcal{ O}(h^{3/2}(\ln
(1/h))^{1/2})\hbox{ in }\langle \cdot \rangle^\epsilon L^2.
\end{multline*}
Here the exponential factor corresponds to the action of
$\widehat{\tau }_\omega $ which annihilates the one of $\widehat{\tau
}_{-\omega }$ and we get
\begin{equation}\label{wn.14}
\phi _1^1=\frac{h}{4k}E(1_\Omega )+\mathcal{ O}(1)h^{3/2}(\ln (1/h))^{1/2}
\hbox{ in }\langle \cdot \rangle^\epsilon L^2. 
\end{equation}

\section{Numerical results for the characteristic function of the disk}\label{num}
\setcounter{equation}{0}
In this section we present a detailed numerical study of the system 
(\ref{dbarphi}) for the
characteristic function of the disk. The goal is to compare the 
asymptotic formulae for large \( |k| \) of the previous sections to 
numerical results in this case, and to show that the asymptotic 
formulae allow for a hybrid approach in practice: the asymptotic formulae give 
 a correct description of the solutions with prescribed precision
for values of \( |k| > |k_{c}| \) where \( k_{c} \) is such that the 
numerical solution of the system (\ref{dbarphi}) for \( |k|\leq 
|k_{c}| \) is correct to the same order of accuracy. Thus a 
combination of numerical and semi-classical techniques allows to give 
a solution (with prescribed precision) of the system (\ref{dbarphi}) 
for all values of \( k\in \mathbb{C} \).

\subsection{Numerical approach}
Here we briefly summarize the numerical approach \cite{KS19} for 
potentials with compact support on a disk (for simplicity we only 
consider the unit disk). Note that for the reasons discussed in the 
introduction (possible non-uniqueness of solutions for $\sigma=-1$), 
we only consider the case $\sigma=1$ in (\ref{dbarphi}), i.e., the 
defocusing case for DS II. 

We write 
$z=re^{i\varphi}$ in the disk and $z=e^{i\varphi}/s$ in its complement, 
thus $r\in[0,1]$ and $s\in[0,1]$. System (\ref{dbarphi}) reads in 
polar coordinates
\begin{equation}
\begin{array}{l}
    \mathrm{e}^{\mathrm{i}\varphi}\left(\partial_{r}+\frac{\mathrm{i}}{r}\partial_{\varphi}\right)\phi_{1}  
    =q(r,\varphi)\mathrm{e}^{\bar{k}\bar{z}-kz}\phi_{2},
    \\~\\
     \mathrm{e}^{-\mathrm{i}\varphi}\left(\partial_{r}-\frac{\mathrm{i}}{r}\partial_{\phi}\right)\phi_{2} 
    =\bar{q}(r,\varphi)\mathrm{e}^{kz-\bar{k}\bar{z}}\phi_{1}.
    \label{dbarpolphi}
\end{array}
\end{equation}
In the exterior of the disk, $\phi_{1}$ is a holomorphic function 
tending to 1 at infinity, 
 and $\phi_{2}$ is an anti-holomorphic 
function vanishing at infinity, 
\begin{equation*}
\phi_{1}=1+\sum_{n=1}^{\infty}a_{n}z^{-n}, \quad \phi_{2}=\sum_{n=1}^{\infty}b_{n}\bar{z}^{-n},
\end{equation*}
where $a_{n}$, 
$b_{n}$ are constants for $n=1,2,\ldots$

The system (\ref{dbarpolphi}) is numerically solved in \cite{KS19} by a 
Chebychev-Fourier method. This means that the functions \( \phi_{1} \) 
and \( \phi_{2} \) are approximated by trigonometric polynomials in 
\( \phi \) and by Chebychev polynomials in \( r \), 
\begin{equation}
\begin{array}{l}
    \phi_{1}\approx 
	\sum_{n=-N_{\varphi}/2}^{N_{\varphi}/2-1}\sum_{m=0}^{N_{r}}a_{nm} 
	T_{m}(l)e^{2\pi in/N_\varphi},\\
	\phi_{2}\approx 
	\sum_{n=-N_{\varphi}/2+1}^{N_{\varphi}/2}\sum_{m=0}^{N_{r}}b_{nm} 
	T_{m}(l)e^{2\pi i n/N_\varphi},
	\label{cheb}
\end{array}
\end{equation}
where \( T_{m}=\cos(m\arccos(x)) \), \( m\in \mathbb{N} \) are the 
Chebychev polynomials. Regularity of the solution of 
(\ref{dbarpolphi}) for \( r\to0 \) as 
well as the matching conditions at the rim of the disk uniquely 
determine the solution. The finite dimensional system following with 
(\ref{cheb}) from (\ref{dbarpolphi}) for the coefficients \( a_{nm} 
\), \( b_{nm} \) is solved via a fixed point iteration, see 
\cite{KS19} for details. Since it is known that the coefficients of a 
Chebychev and a Fourier series are exponentially decreasing for an 
analytic function, the decrease of the coefficients \( a_{nm} \), \( 
b_{nm} \) for large \( |n|,m\ \) indicates the numerical resolution 
of the problem and allows to estimate the numerical error, see again 
\cite{KS19}. 

\begin{figure}[htb!]
  \includegraphics[width=0.32\textwidth]{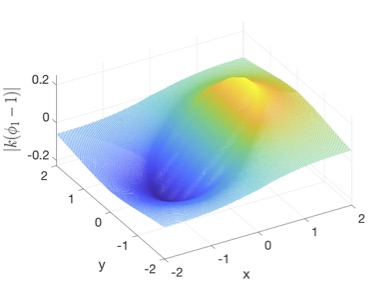}
  \includegraphics[width=0.32\textwidth]{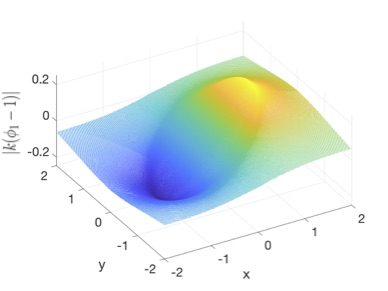}
  \includegraphics[width=0.32\textwidth]{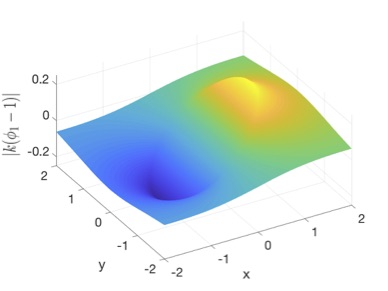}
 \caption{Difference between the solution $\phi_{1}$ for the 
 characteristic function of the disk and 
1 multiplied by 
 $k$ for $k=10,100,1000$ from left to right.}
 \label{Phi1}
\end{figure}

We now apply this numerical approach to the case of \( q \) being the 
characteristic function of the unit disk.
Because of the radial symmetry of $q$, we can concentrate on values 
of $k>0$ without loss of generality. 
In Fig.~\ref{Phi1} we show the results of a 
numerical computation of the modulus of $\phi_{1}-1$ for three 
values of $k$. It can be seen that this difference decreases as 
$1/k$ in agreement with (\ref{wn.14}). 

\begin{figure}[htb!]
  \includegraphics[width=0.32\textwidth]{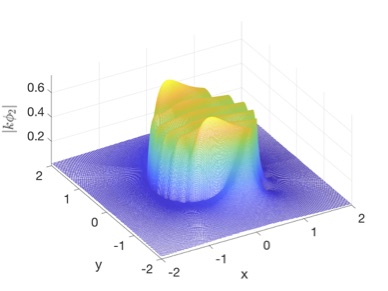}
  \includegraphics[width=0.32\textwidth]{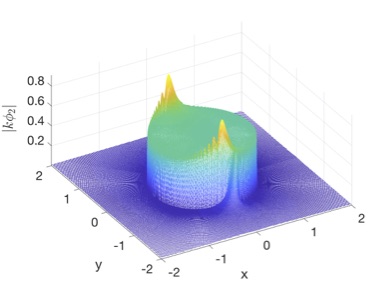}
  \includegraphics[width=0.32\textwidth]{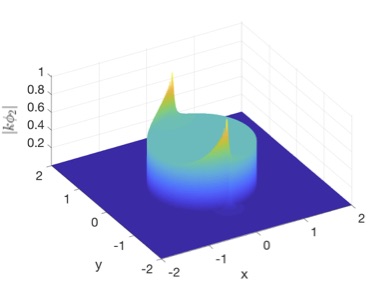}
 \caption{The solution $\phi_{2}$ for the 
 characteristic function of the disk multiplied by 
 $k$ for $k=10,100,1000$ from left to right.}
 \label{Phi2}
\end{figure}

In Fig.~\ref{Phi2} we show the corresponding plots for $\phi_{2}$. A 
scaling proportional to $1/k$ as in (\ref{wn.12}) 
is not obvious for all shown values of 
$k$. It appears to be realized for the higher values $k=100,1000$.

\begin{remark}
	In this section we will always show the pointwise difference 
	between the 
	solution to the d-bar system and various asymptotic formulae for 
	the latter. Note, however, that the asymptotic formulae have been 
	derived for some weighted $L^{2}$ norms. Thus the found 
	differences near the maxima in Fig.~\ref{Phi2} will contribute 
	much less 
	in the $L^{2}$ spaces than shown here, where 
	the agreement is already very good.
\end{remark}

\subsection{Asymptotic formulae}

\begin{figure}[htb!]
  \includegraphics[width=0.32\textwidth]{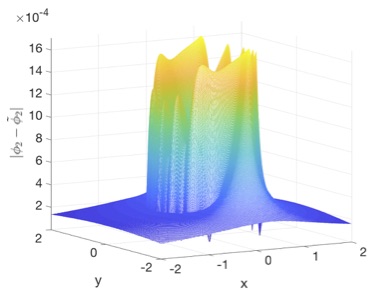}
  \includegraphics[width=0.32\textwidth]{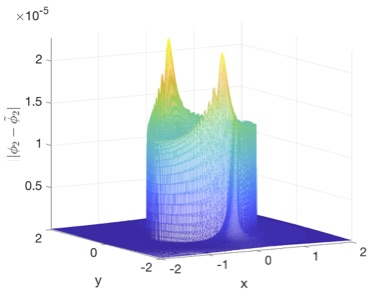}
  \includegraphics[width=0.32\textwidth]{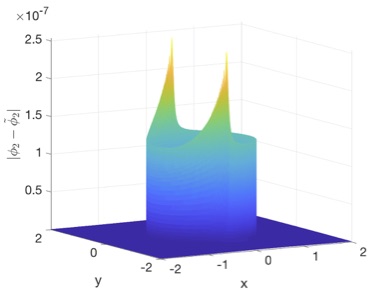}
 \caption{Difference between the solution $\phi_{2}$ for the 
 characteristic function of the disk and 
 $\tilde{\phi}_{2}$ (\ref{Phi2int}) for $k=10,100,1000$ from left to right.}
 \label{Phi2diff}
\end{figure}

The results of Section \ref{ic} imply that 
$\phi_{1}=1+\mathcal{O}(1/|k|)$ for $|k|\to\infty$, see (\ref{wn.14}). 
Thus in leading order of $1/k$ 
the second equation in (\ref{dbarphi}) has the approximate solution
\begin{align}
    \tilde{\phi}_{2}&=\frac{1}{2\pi}\int_{|w|\leq1}^{}\frac{\mathrm{e}^{kw-\bar{k}\bar{w}}}{\bar{z}-\bar{w}}d^{2}w
    =\frac{1}{4\pi 
	k}\int_{0}^{2\pi}\frac{\mathrm{e}^{ke^{-i\varphi}-\bar{k}e^{i\varphi}}-\mathrm{e}^{\bar{k}\bar{z}-kz}}{\bar{z}-e^{-i\varphi}}e^{-i\varphi}d\varphi
	\nonumber\\
	& =\frac{1}{2\pi}\bar{f}(z,k)
    \label{Phi2int},
\end{align}
i.e.\ up to a factor the complex conjugate of the integral in 
\ref{ic.0}, see also (\ref{ic.1}). We first check how well the function \( \tilde{\phi}_{2} \) of
(\ref{Phi2int}) approximates \( \phi_{2} \) for large \( k \). 
Since the integral (\ref{Phi2int}) is singular 
near the boundary and highly oscillatory, it is numerically 
challenging to evaluate. Therefore we compute it  by numerically inverting the $\partial$ 
operator, the same way as when solving
the system (\ref{dbarphi}). The difference between $\phi_{2}$ and 
$\tilde{\phi}_{2}$ 
can be seen in Fig.~\ref{Phi2diff}. It appears to scale as $1/k^{2}$ 
(see also (\ref{wn.12})). 

The task is thus to compute the function $f(z,k)$ of (\ref{ic.0}) 
in (\ref{Phi2int}) to leading order in $1/|k|$  as in Section 
\ref{ic}. We briefly recall the main  steps for the example of the characteristic function of the unit disk: 
We consider the holomorphic extension in $w$ of 
the integrand by noting that $\bar{w}=1/w$ on the unit circle. Thus 
we have 
\begin{equation}
    \frac{1}{2\pi}\bar{f}(z,k)=\frac{1}{4\pi 
	ki}\int_{|w|=1}^{}\frac{\mathrm{e}^{k/w-\bar{k}w}}{\bar{z}-w}dw
    \label{int2},
\end{equation}
which was computed in Section \ref{ic} asymptotically via a contour deformation and 
steepest descent techniques. We illustrate the various steps to obtain the 
asymptotic formula (\ref{wn.14}) for the unit disk below.  The interior of the disk $r\leq 1$ 
 and its complement in 
the complex plane given by $s:=1/r<1$  will be always shown 
separately in dependence of polar coordinates. 
The exponent in the integrand of (\ref{int2}) has the 
stationary points $w_{\pm}=\pm i$. The fact that stationary phase 
approximations are essentially quadratic approximations of the phase 
means that length scales of order $1/\sqrt{|k|}$ are important. This  
leads to a natural decomposition of the disk and its complement in the 
complex plane into zones. Let $r_{k}$ be such that 
$1-r_{k}=\mathcal{O}(1/\sqrt{|k|})$. 
We consider the following cases:\\
I.1 $r<r_{k}$ respectively $s<r_{k}$;\\
I.2 $r_{k}<r<1$ respectively $r_{k}<s<1$ and 
$\pi/2+\delta<\varphi<3\pi/2-\delta$ where $\delta>0$ such that 
$\delta=O(1/\sqrt{|k|})$ or
$0\leq \varphi<\pi/2-\delta$ or $3\pi/2+\delta<\varphi<2\pi$;\\
II. $r_{k}<r<1$ respectively $r_{k}<s<1$ and 
$\pi/2-\delta \varphi<\pi/2+\delta$ or $3\pi/2-\delta 
\varphi<3\pi/2+\delta$.

\begin{figure}[htb!]
  \includegraphics[width=0.32\textwidth]{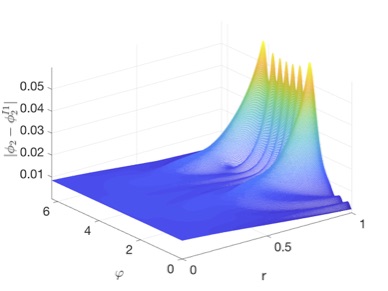}
  \includegraphics[width=0.32\textwidth]{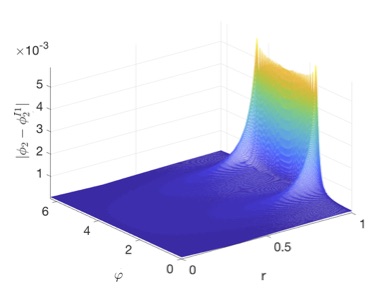}
  \includegraphics[width=0.32\textwidth]{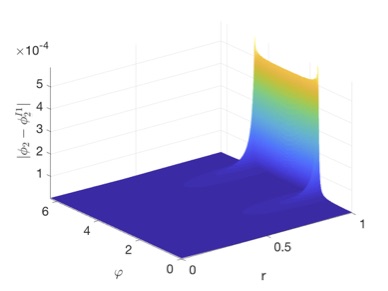}\\
  \includegraphics[width=0.32\textwidth]{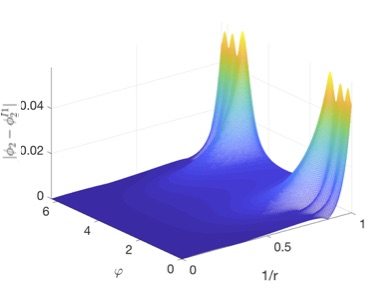}
  \includegraphics[width=0.32\textwidth]{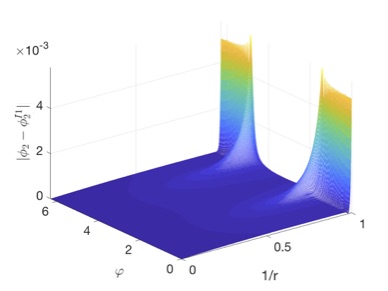}
  \includegraphics[width=0.32\textwidth]{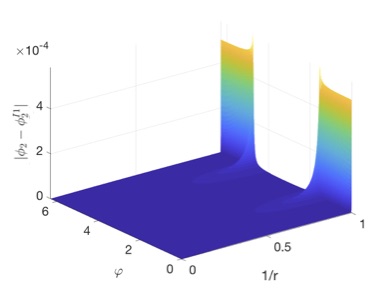}
 \caption{Difference between the solution $\phi_{2}$ for the 
 characteristic function of the disk and 
 $\phi_{2}^{I1}$, in the 
 upper row the interior of the disk, in the lower row the exterior of 
 the disk, both   for $k=10,100,1000$ from left to right.}
 \label{Phi2diff2}
\end{figure}

\paragraph{\textbf{Case I.1}}
In this case the integral (\ref{int2}) can be evaluated with a 
standard stationary phase approximation, see Section \ref{ic}. Its leading order 
contribution in $1/|k|$  to $\phi_{2}$ is with (\ref{wn.3}) and 
(\ref{icbtf.2})
\begin{equation}
	\begin{cases}
    \phi_{2}^{I1}&=\frac{e^{kz-\bar{k}\bar{z}}}{4k},\quad  |z|\leq 
	1,\\
    \phi_{2}^{I1}&=0,\quad |z|>1,
	\end{cases}	
    \label{I}
\end{equation}
plus corrections of  order $O(1/|k|^{3/2})$.
In Fig.~\ref{Phi2diff2}, we show the difference of $\phi_{2}$ and the 
$\phi_{2}^{I1}$ of (\ref{I}), in 
the upper row in the interior of the disk,  in the lower row in the 
complement of the disk in the complex plane, both in polar coordinates.
It can be seen that the approximation is as expected not of the 
wanted order near the rim of the disk. The precise region of 
applicability of the approximation will be studied below.

\paragraph{\textbf{Case I.2}}
Writing $w = \omega\exp(i\psi)$, $\omega>0$, $\psi\in\mathbb{R}$,  we get for the exponent in the 
integral (\ref{int2})
\begin{equation}
    k\left(w-\frac{1}{w}\right) = k\cos\psi 
    \left(\omega-\frac{1}{\omega}\right)
    +ik\sin\psi \left(\omega+\frac{1}{\omega}\right)
    \label{exp}.
\end{equation}
This means  we can deform the integration contour in (\ref{int2}) as 
in Fig.~\ref{contour} for 
$\pi/2<\psi<3\pi/2$ to a semicircle with $\omega<1$ to get an 
exponentially small contribution to the integral there for $k$ large. 
Similarly we deform the contour for $0\leq \psi<\pi/2$ or 
$3\pi/2<\psi<2\pi$ to a semicircle with $\omega>1$ to get again an 
exponentially small contribution to the integral. In both cases it is 
possible to pick up a contribution due to a residue if the pole of 
the integrand is crossed by the deformation. This gives the 
following leading order contributions to $\phi_{2}$ in these cases: For 
$1-r_{k}<r<1$  and $\pi/2<\varphi<3\pi/2$, one has 
\begin{equation}
    \phi_{2}^{I2}=\frac{1}{2k}\left(\mathrm{e}^{kz-\bar{k}\bar{z}}
    -e^{ke^{i\varphi}/r-\bar{k}e^{-i\varphi}r}\right)
    \label{case2}.
\end{equation}
The asymptotic formula in the complement of the disk does not change 
in this case since no residue can be picked up.

\begin{figure}[htb!]
  \includegraphics[width=0.32\textwidth]{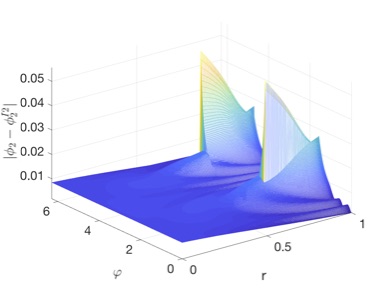}
  \includegraphics[width=0.32\textwidth]{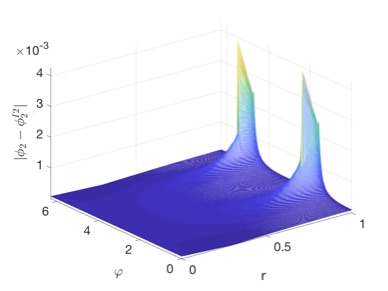}
  \includegraphics[width=0.32\textwidth]{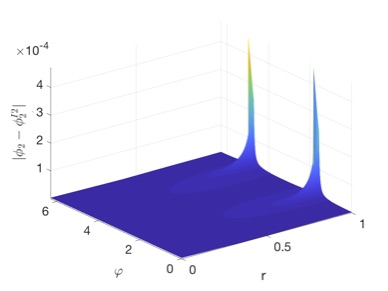}
 \caption{Difference between the solution $\phi_{2}$ for the 
 characteristic function of the disk and 
 $\phi_{2}^{I2}$ from  (\ref{case2}) at the disk 
 for $k=10,100,1000$ from left to right.}
 \label{Phi2diffcase2}
\end{figure}

If we take this enhanced asymptotic description into account, the upper row of figures in Fig.~\ref{Phi2diff2} is now replaced by the figures in 
Fig.~\ref{Phi2diffcase2}. The expected behavior can be seen except in 
the vicinity of the points $w_{\pm}$. There the difference is still 
linear in $1/k$. 

\begin{figure}[htb!]
  \includegraphics[width=0.32\textwidth]{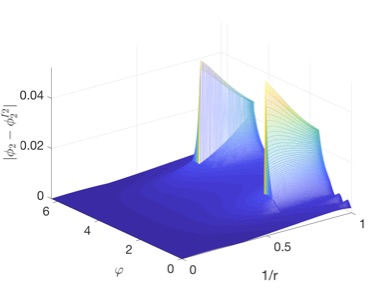}
  \includegraphics[width=0.32\textwidth]{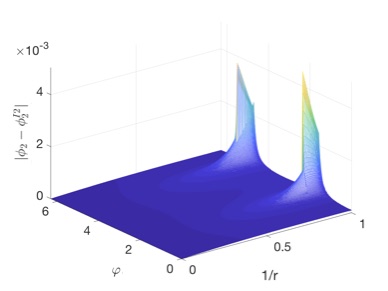}
  \includegraphics[width=0.32\textwidth]{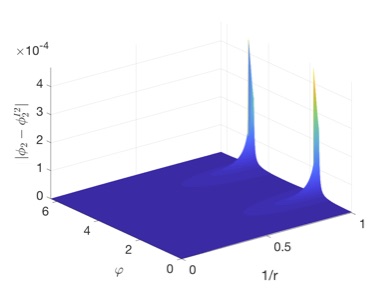}
 \caption{Difference between the solution $\phi_{2}$ for the 
 characteristic function of the disk and 
 $\phi_{2}^{I2}$ from  (\ref{case3}) for the 
 exterior of the disk 
 for $k=10,100,1000$ from left to right.}
 \label{Phi2diffcase3}
\end{figure}

For $\varphi<\pi/2$ or $\varphi>3\pi/2$ there is no contribution due to a 
residue in the interior of the disk. But in its complement in 
$\mathbb{C}$, we get
\begin{equation}
    \phi_{2}^{I2}=\frac{1}{2k}
    e^{ke^{i\varphi}s-\bar{k}e^{-i\varphi}/s}
    \label{case3}.
\end{equation}
This leads to Fig.~\ref{Phi2diffcase3} which shows the same behavior as 
Fig.~\ref{Phi2diffcase2} for the interior of the disk.

\paragraph{\textbf{Case II}}
We now address the case that $z$ is close to the stationary points 
$w_{\pm}$ of 
the exponent in (\ref{int2}), $|z-w_{\pm}|=\mathcal{O}(1/\sqrt{|k|})$. 
In Section \ref{ic}  a quadratic 
approximation to the exponent was considered. We put $w = 
w_{\pm}+\xi$ and get for $|\xi|\ll1$ for the exponent
$w-1/w = \pm2i(1-\xi^{2}/2)+\mathcal{O}(|\xi|^{3})$. 
As the integration path we use the line 
$\xi=a_{\pm}\eta$, where $\eta\in\mathbb{R}$ and where  
$a_{\pm}=\sqrt{w_{\pm}/(2k)}$ in order to get an integrand exponentially 
decaying on the integration path. 
We consider the function (\ref{lp.14})
\begin{equation}
    G(z) := 
    \int_{-\infty}^{\infty}\frac{e^{-t^{2}/2}}{z-t}dt
    \label{G}.
\end{equation}
Note that the function $G$ is not uniquely defined by (\ref{G}) 
because of the pole on the real axis which is also the integration contour. 
We denote by $G_{r}(z)$ the analytical continuation to the whole 
complex plane of 
the function obtained by computing $G$ in standard way for $\Im z>0$, 
and $G_{l}(z)$ for $\Im z<0$. 

\begin{figure}[htb!]
  \includegraphics[width=0.49\textwidth]{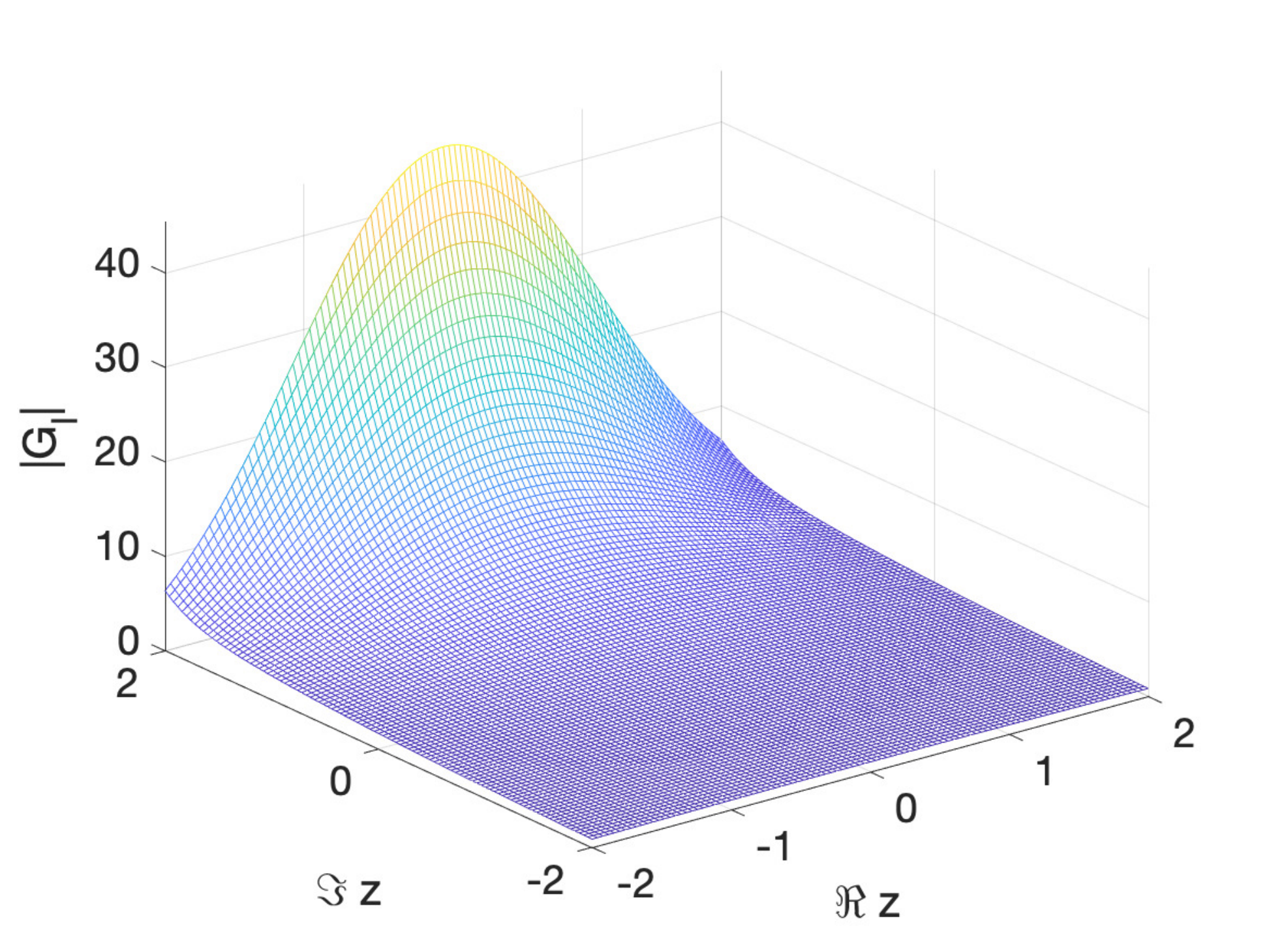}
  \includegraphics[width=0.49\textwidth]{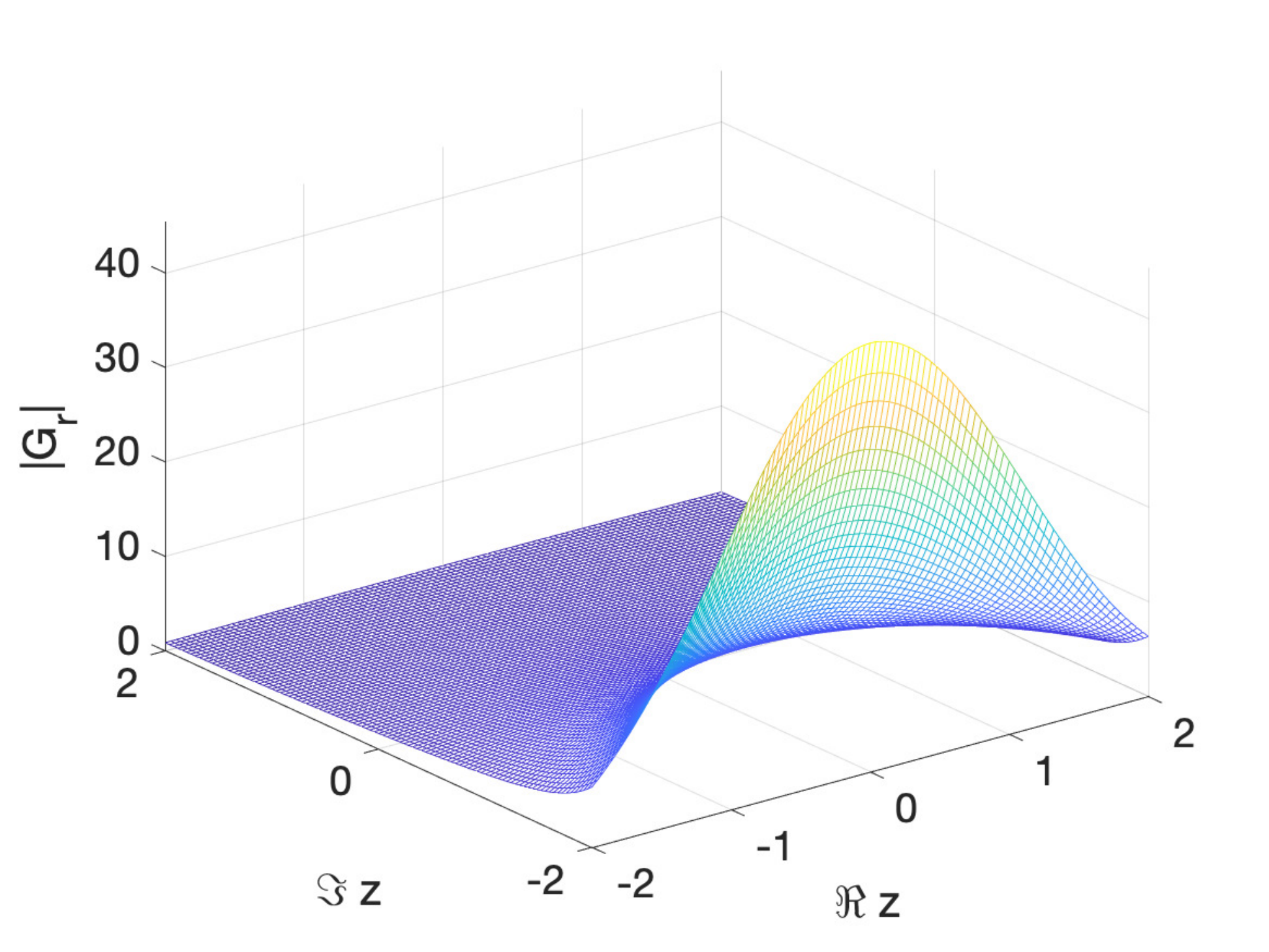}
 \caption{Moduli of the function $G_{l}$ on the left and $G_{r}$ on 
 the right.}
 \label{figG}
\end{figure}

Numerically these functions are computed on the parallels to the real 
axis going through $-i(\Im z +3)$ for $G_{l}$ and $i(\Im z +3)$ for 
$G_{r}$. On these lines, the integrand is approximated via a 
truncated Fourier series on a sufficiently large period, \( t\in 
L[-\pi,\pi] \) (we use \( L=10 \) in the following). We compute in  
\( t \) the standard  discrete Fourier transform, i.e., sample the 
integrand on
\( t_{n}=L(-\pi+nh) \), \( n=1,\ldots,N \), where \( N\in \mathbb{N} 
\) is the number of collocation points and where \( h=2\pi/N \). The 
integral in (\ref{G}) is the Fourier coefficient with index 0 of the 
discrete Fourier transform of this function, i.e., simply the sum over 
\( n \) of 
the integrand in (\ref{G}) sampled at the collocation points \( t_{n} 
\). Since this is one of the coefficients of the discrete Fourier 
transform, the resulting numerical method is a so-called spectral 
method. 
This means the numerical error in approximating the integrand (which 
is analytic on the chosen integration path) decreases exponentially 
with \( N \). The
numerical 
accuracy is controlled via the decay of the discrete Fourier 
coefficients which can be computed with a fast Fourier transform.
We show both functions $G_{l,r}$ in 
Fig.~\ref{figG}. The difference between the functions on the real axis 
is according to (\ref{lp.15}) equal to $2\pi i \exp(-x^{2}/2)$. The functions satisfy the symmetry relation
\begin{equation}
    G_{l}(\bar{z}) = \overline{G_{r}(z)}
    \label{sym}.
\end{equation}

Thus we get for the integral (\ref{int2})
\begin{equation}
	\begin{split}
		   \frac{\bar{f}(k,z)}{2\pi}\approx &\frac{e^{k(1/w_{\pm}-w_{\pm})}}{4\pi ik}\int_{\xi_{1}}^{\xi_{2}}\frac{
    e^{-k\xi^{2}/w_{\pm}}}{\bar{z}-w_{\pm}-\xi}d\xi\\
	&=
     \frac{e^{k(1/w_{\pm}-w_{\pm})}}{4\pi ik}\int_{\eta_{1}}^{\eta_{2}}\frac{
    e^{-\eta^{2}/2}}{(\bar{z}-w_{\pm})/a_{\pm}-\eta}d\eta
	\end{split}
    \label{square}.
\end{equation}
Since the integrand is exponentially 
decaying, we finally arrive for $|\bar{z}+i|\leq C/|k|$, where $C$ is 
some positive constant, at the approximations for the interior of the 
disk 
\begin{equation}
    \phi_{2}^{II}=
    \begin{cases}
    	\frac{e^{2ik}}{4\pi ik}G_{r}\left(
    \frac{\bar{z}+i}{a_{-}}\right)+\frac{1}{2k}\left(\mathrm{e}^{kz-\bar{k}\bar{z}}
    -e^{ke^{i\varphi}/r-\bar{k}e^{-i\varphi}r}\right), & \Im z +1\leq \Re z \\
    	\frac{e^{2ik}}{4\pi ik}G_{l}\left(
    \frac{\bar{z}+i}{a_{-}}\right)+\frac{e^{kz-\bar{k}\bar{z}}}{2k}, & \Im z +1> \Re z
    \end{cases}
    \label{square2}
\end{equation}
and
\begin{equation}
  \phi_{2}^{II}=
    \begin{cases}
    	\frac{e^{2ik}}{4\pi ik}G_{l}\left(
    \frac{\bar{z}+i}{a_{-}}\right), & \Im z +1\leq \Re z \\
    	\frac{e^{2ik}}{4\pi ik}G_{r}\left(
    \frac{\bar{z}+i}{a_{-}}\right)-\frac{1}{2k}
    e^{ke^{i\varphi}/r-\bar{k}e^{-i\varphi}r}, & \Im z +1> \Re z
    \end{cases}
    \label{square2e}
\end{equation}
in the exterior of the disk. Analogous formulae hold for 
$|\bar{z}-i|\leq C/|k|$. 

\begin{figure}[htb!]
  \includegraphics[width=0.32\textwidth]{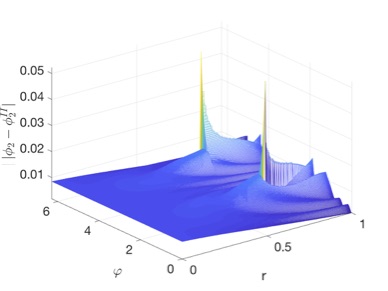}
  \includegraphics[width=0.32\textwidth]{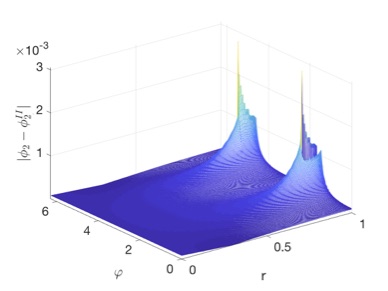}
  \includegraphics[width=0.32\textwidth]{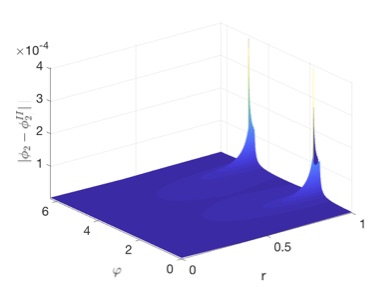}
 \caption{Difference between the solution $\phi_{2}$ for the 
 characteristic function of the disk and 
 $\phi_{2}^{II}$ from (\ref{square2}) for the 
 interior of the disk 
 for $k=10,100,1000$ from left to right ($C=1$).}
 \label{Phi2diffstaint}
\end{figure}

In Fig,~\ref{Phi2diffstaint} we show the effect of all above 
asymptotic descriptions for the interior of the disk. Approximation 
(\ref{square2}) is applied for $|\bar{z}-w_{\pm}|<C/\sqrt{|k|} $ for 
$C=1$, see (\ref{wn.3}). It can be seen that the approximation is excellent near the 
points $w_{\pm}$, the error is largest near these points where the 
approximation of case I is applied.

\begin{figure}[htb!]
  \includegraphics[width=0.32\textwidth]{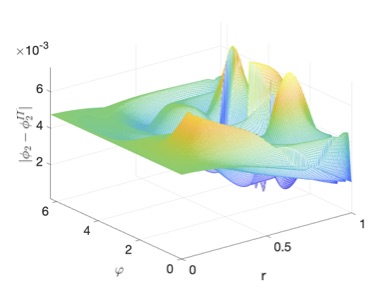}
  \includegraphics[width=0.32\textwidth]{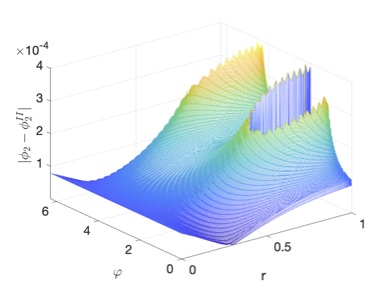}
  \includegraphics[width=0.32\textwidth]{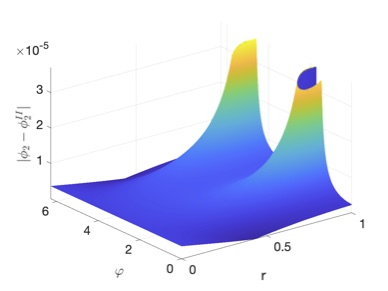}
 \caption{Difference between the solution $\phi_{2}$ for the 
 characteristic function of the disk and 
 $\phi_{2}^{II}$ from (\ref{square2}) for the 
 interior of the disk 
 for $k=10,100,1000$ from left to right ($C=4$).}
 \label{Phi2diffstaintC4}
\end{figure}

Since the error near the points $w_{\pm}$ in 
Fig.~\ref{Phi2diffstaint} is largest where the approximation 
(\ref{square2}) is not applied, it appears reasonable that larger 
values of the constant $C$ should be considered. 
The asymptotic formulae of section \ref{ic} do not fix this 
constant. Since we consider values of $k$ as low as 10, we cannot 
choose $C$ too large since otherwise the regions in the vicinity of 
$w_{\pm}$ would overlap. In Fig. \ref{Phi2diffstaintC4} we show the 
same differences as in Fig. \ref{Phi2diffstaint}, but this time for 
$C=4$. The overall error is considerably lower in this case than for 
$C=1$ and is still dominated by the regions $|z-w_{\pm}|>C/\sqrt{k}$. 
For large $|k|$ one could optimize the choice of $C$, but this is 
beyond the goal of this paper. 

\begin{figure}[htb!]
  \includegraphics[width=0.32\textwidth]{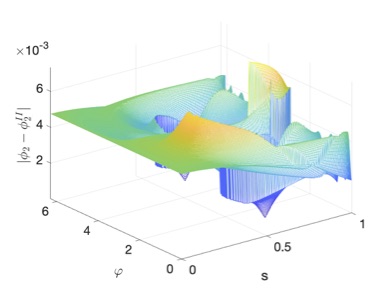}
  \includegraphics[width=0.32\textwidth]{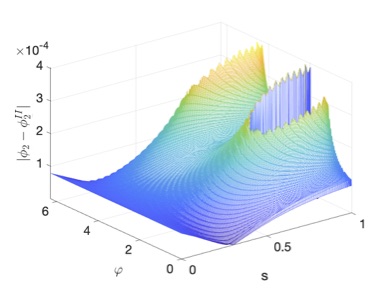}
  \includegraphics[width=0.32\textwidth]{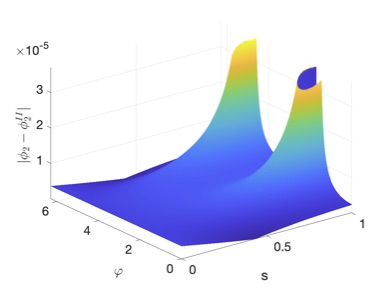}
 \caption{Difference between the solution $\phi_{2}$ for the 
 characteristic function of the disk and 
 $\phi_{2,e}^{III}$ from (\ref{square2}) for the 
 exterior of the disk 
 for $k=10,100,1000$ from left to right.}
 \label{Phi2diffstaext}
\end{figure}

A similar behavior can be seen in the complement of the disk in the 
complex plane in Fig.~\ref{Phi2diffstaext}.

\subsection{Reflection coefficient}
As stated in the introduction, the main quantity of interest in an 
inverse scattering approach to DS II is the reflection coefficient 
(\ref{eq:r-def}) which plays here the role of the Fourier transform for 
linear equations and is the angle in terms of action-angle variables.

The analysis in Section \ref{ic} has shown that the function 
$\phi_{1}$ is given in leading order by the expresssion 
(\ref{wn.14}). In the case of the unit disk we are interested in 
here, this takes the form  
\begin{equation}
	\tilde{\phi}_{1}=1+\frac{\bar{z}}{4k} 
    \label{Phi1t}
\end{equation}
for $|k|\to\infty$ 
($\phi_{1}=\tilde{\phi}_{1}+\mathcal{O}(|k|^{-3/2}|k|)$). 
In the exterior of the disk, the function is holomorphic and tends to 
1 at infinity. Since it is continuous at the disk, we have 
$\tilde{\phi}_{1}=1+\frac{1}{4kz}$ for $|z|>1$. 

\begin{figure}[htb!]
  \includegraphics[width=0.32\textwidth]{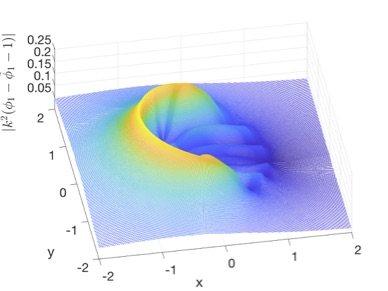}
  \includegraphics[width=0.32\textwidth]{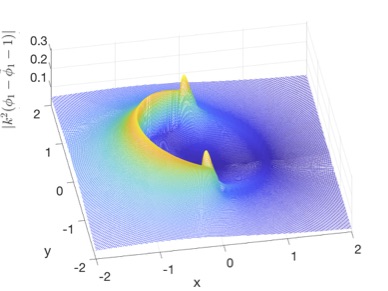}
  \includegraphics[width=0.32\textwidth]{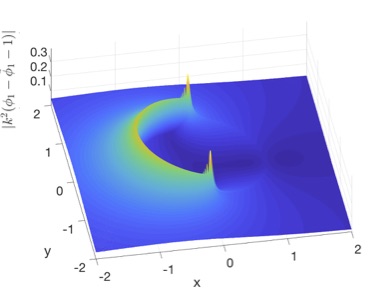}
 \caption{Difference between the solution $\phi_{1}$ for the 
 characteristic function of the disk and 
 $1+\frac{\bar{z}}{4k}$ multiplied by 
 $k^{2}$ for $k=10,100,1000$ from left to right.}
 \label{Phi1diff}
\end{figure}

We show in Fig.~\ref{Phi1diff} the difference 
between $\phi_{1}$ and $\tilde{\phi}_{1}$. This difference is largest near 
the rim of the disk, but appears to be of order $1/k^{2}$. 

\begin{figure}[htb!]
  \includegraphics[width=0.49\textwidth]{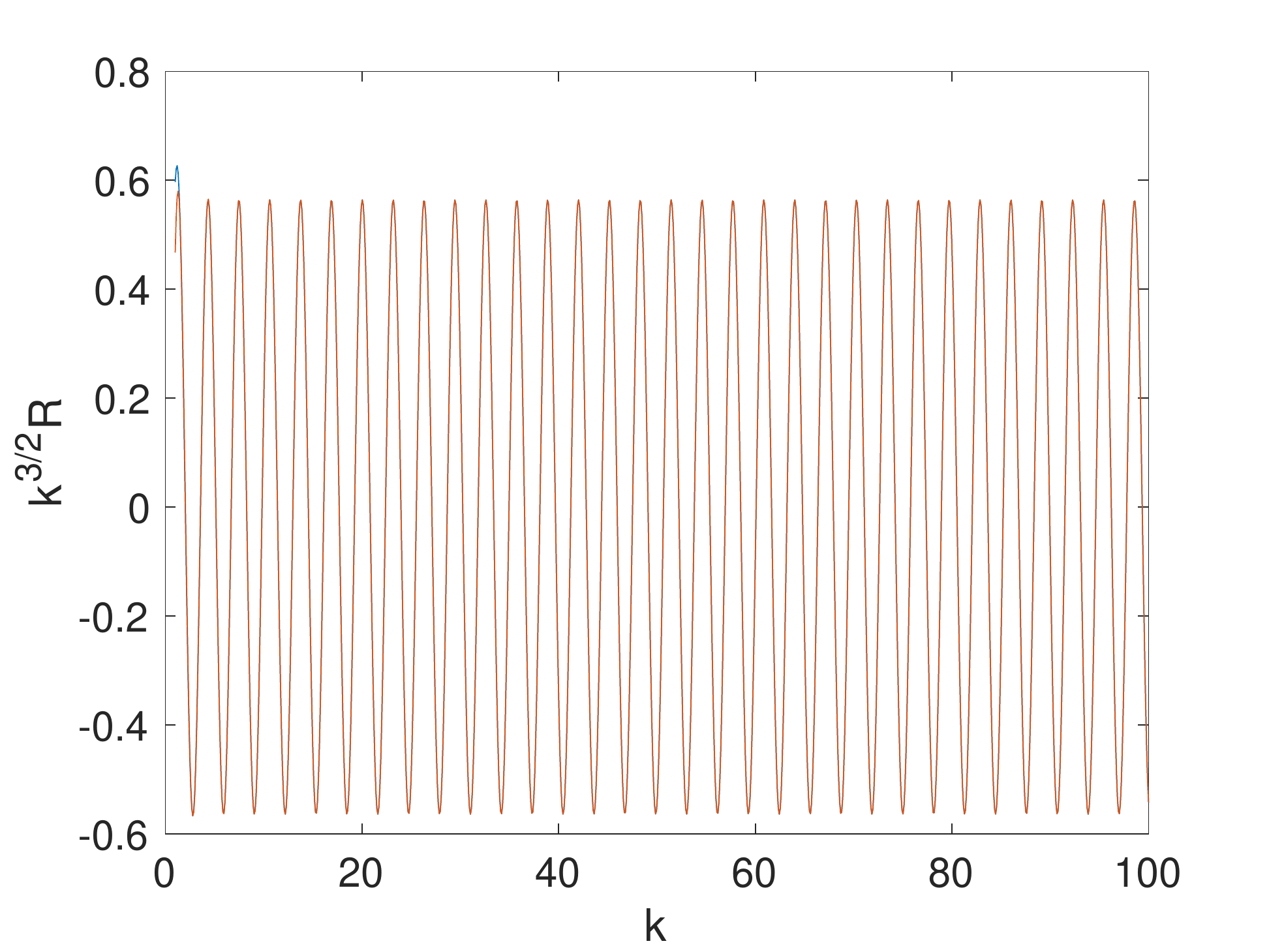}
  \includegraphics[width=0.49\textwidth]{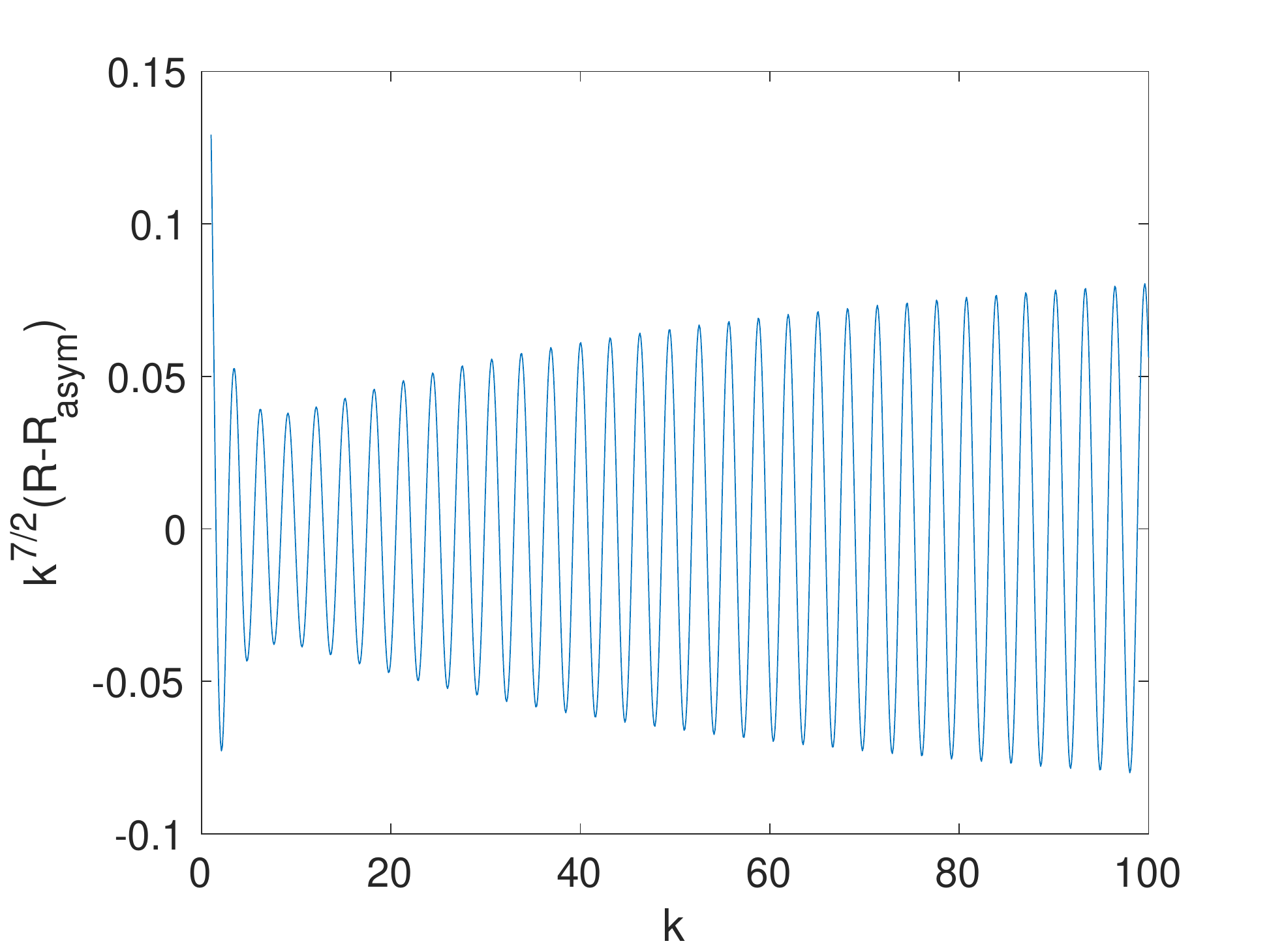}
 \caption{Reflection coefficient for the characteristic function of 
 the disk, on the left $R$ in blue and $R_{asym}$ from (\ref{Rintsp}) 
 in red, both multiplied with $k^{3/2}$, on the right the difference 
 between both multiplied with $k^{7/2}$.}
 \label{figref}
\end{figure}

The reflection coefficient is given via $\bar{R}(k) = 2\lim 
_{z\to\infty}\bar{z}\phi_{2}$, i.e.,
\begin{equation}
    \bar{R} = 
    \frac{2}{\pi}\int_{|z|\leq1}^{}\mathrm{e}^{kw-\bar{k}\bar{w}}\phi_{1}d^{2}w
    \approx     \frac{2}{\pi}\int_{|z|\leq1}^{}\mathrm{e}^{kw-\bar{k}\bar{w}}
    \left(1+\frac{\bar{w}}{4k}\right)d^{2}w
    \label{Rint}.
\end{equation}
The reflection coefficient is real in this case. 
Note that the error term in (\ref{wn.14}) can 
contribute in the oscillatory integral (\ref{Rint}) in the order we 
would like to study. If we conjecture that this is not the case, then 
the integral (\ref{Rint}) can be computed once more with a stationary phase 
approximation (of higher order), which allows us to study higher order terms. After some calculation we get
\begin{equation}
    R    \approx   R_{asym}:=  \frac{1}{\sqrt{\pi k^{3}}}\left(
    \sin(2k-\pi/4)-\frac{5}{16k}\cos(2k-\pi/4)\right)
    \label{Rintsp}.
\end{equation}
We show the reflection coefficient in Fig.~\ref{figref} on the left 
in blue. 
The asymptotic formula for the coefficient shown in the same figure 
in red
agrees so well with the coefficient even for values of $k$ of the 
order 10 that we show on the right of the same figure the difference between $R$ and the 
asymptotic formula (\ref{Rintsp}) multiplied with $k^{7/2}$. This 
indicates that the corrections to formula (\ref{Rintsp}) are of the 
order $k^{-7/2}$, but that this asymptotic regime is only reached for 
larger values of $k$ than shown. 

An error term of the order of $|k|^{-7/2}$ implies that for $k=1000$, 
which can be reached numerically at least with an accuracy of the 
order of $10^{-11}$ as discussed in \cite{KS}, will be of the order 
$10^{-11}$. This means that the asymptotic formula for the reflection 
coefficient is applicable already for values of the spectral 
parameter $k$ where the numerical approach can be applied. In other 
words, in the example of the characteristic function of the disk, a 
hybrid approach of a numerical approach for $|k|<k_{0}$ with 
$k_{0}\sim 10^{3}$ combined with the asymptotic formula 
(\ref{Rintsp}) for $|k|>k_{0}$ allows for a computation of the 
reflection coefficient for all $k\in \mathbb{C}$ with an accuracy of 
$10^{-11}$ and better.

\end{document}